\documentclass[sn-mathphys,Numbered]{sn-jnl}

\usepackage{graphicx}%
\usepackage{multirow}%
\usepackage{amsmath,amssymb,amsfonts}%
\usepackage{amsthm}%
\usepackage{siunitx}
\usepackage{dutchcal}

\usepackage{xcolor}%
\usepackage{textcomp}%
\usepackage{manyfoot}%
\usepackage{booktabs}%
\usepackage{algorithm}%
\usepackage{algorithmicx}%
\usepackage{algpseudocode}%
\usepackage{listings}%
\usepackage{chemfig}
\usepackage[export]{adjustbox}
\usepackage{bm}

\theoremstyle{thmstyleone}%
\newtheorem{theorem}{Theorem}

\theoremstyle{thmstyletwo}%
\newtheorem{remark}{Remark}%

\theoremstyle{thmstylethree}%
\newtheorem{definition}{Definition}%

\usepackage{calrsfs}
\DeclareMathAlphabet{\pazocal}{OMS}{zplm}{m}{n}
\newtheorem{lemma}{Lemma}
\newtheorem{Corollary}{Corollary}
\usepackage{gensymb}
\usepackage{tikz}
\usetikzlibrary{matrix,backgrounds,shapes}
\usetikzlibrary{positioning}
\usetikzlibrary{fit}
\usetikzlibrary{plotmarks}
\usetikzlibrary{decorations.pathreplacing}
\usepackage{pgfplots}
\usetikzlibrary{shapes}
\usetikzlibrary{positioning,arrows}
\usetikzlibrary{fadings,shapes.arrows,shadows}
\usetikzlibrary{arrows.meta}
\usepackage{bm}
\usepackage{xspace}
\usepackage{afterpage}
\usepackage{color}
\usepackage{float}
\definecolor{royalblue}{rgb}{0.254,0.41,0.88}
\usepackage{xcolor}

\DeclareMathAlphabet\mathbfcal{OMS}{cmsy}{b}{n}

\newcommand{\ie}[0]{{i.e.\@}\xspace}

\newtheorem{assume}{Assumption}

\newcommand{\Tr}{\ensuremath{^{\mr{T}}}}

\newcommand{\mr}[1]{\ensuremath{\mathrm{#1}}}
\newcommand{\fnc}[1]{\ensuremath{#1}}
\newcommand{\bfnc}[1]{\ensuremath{\bm{#1}}}
\newcommand{\mat}[1]{\ensuremath{\mathsf{#1}}}

\newcommand{\xm}[1]{\ensuremath{x_{#1}}}
\newcommand{\xil}[1]{\ensuremath{\xi_{#1}}}
\newcommand{\nxm}[1]{\ensuremath{n_{\xm{#1}}}}
\newcommand{\nxil}[1]{\ensuremath{n_{\xil{#1}}}}
\newcommand{\M}[0]{\ensuremath{\mat{P}}}
\newcommand{\Mbar}[0]{\ensuremath{\overline{\mat{P}}}}

\newcommand{\Imat}[1]{\ensuremath{\mat{I}_{#1}}}
\newcommand{\Rxiln}[2]{\ensuremath{\mat{R}_{\xil{#1}}^{(#2)}}}
\newcommand{\Rxilnbar}[2]{\ensuremath{\overline{\mat{R}}_{\xil{#1}}^{(#2)}}}
\newcommand{\Porthoxil}[1]{\ensuremath{\mat{P}_{\xil{#1}}^{\perp}}}
\newcommand{\Porthoxilbar}[1]{\ensuremath{\overline{\mat{P}}_{\xil{#1}}^{\perp}}}

\newcommand{\W}[0]{\ensuremath{\bm{\fnc{W}}}}

\newcommand{\bmU}[0]{\ensuremath{\bm{U}}}
\newcommand{\bmf}[0]{\ensuremath{\bm{R}}}

\newcommand{\bmW}[0]{\ensuremath{\bm{W}}}
\newcommand{\Vmn}[2]{\ensuremath{\mat{C}_{#1,#2}}}

\newcommand{\Vhatmn}[2]{\ensuremath{\widehat{\mat{C}}_{#1,#2}}}

\newcommand{\Vtildemn}[2]{\ensuremath{\widetilde{\mat{C}}_{#1,#2}}}
\newcommand{\Lyap}[0]{\ensuremath{\fnc{V}}}
\newcommand{\bmF}[1]{\ensuremath{\mathfrak{F}_{x_{#1}}}}
\newcommand{\scp}[2]{\left\langle{#1,\, #2}\right\rangle}

\newcommand{\Jxilxm}[2]{\ensuremath{\fnc{J}\frac{\partial\xil{#1}}{\partial\xm{#2}}}}

\newcommand{\bmfI}[1]{\ensuremath{\bm{F}^{(c)}_{x_{#1}}}}
\newcommand{\bmfV}[1]{\ensuremath{\bm{F}^{(d)}_{x_{#1}}}}

\newcommand{\Dxil}[1]{\ensuremath{\mat{D}_{\xil{#1}}}}
\newcommand{\Dxilbar}[1]{\ensuremath{\overline{\mat{D}}_{\xil{#1}}}}
\newcommand{\Qxil}[1]{\ensuremath{\mat{Q}_{\xil{#1}}}}
\newcommand{\Qxilbar}[1]{\ensuremath{\overline{\mat{Q}}_{\xil{#1}}}}
\newcommand{\Exil}[1]{\ensuremath{\mat{E}_{\xil{#1}}}}
\newcommand{\Exilbar}[1]{\ensuremath{\overline{\mat{E}}_{\xil{#1}}}}
\newcommand{\Sxil}[1]{\ensuremath{\mat{S}_{\xil{#1}}}}
\newcommand{\el}[1]{\ensuremath{\bm{e}_{#1}}}
\newcommand{\Morthoxil}[1]{\ensuremath{\mat{P}_{\xil{#1}}^{\perp}}}

\newcommand{\uk}[1]{\ensuremath{\bm{u}^{#1}}}
\newcommand{\wk}[1]{\ensuremath{\bm{w}^{#1}}}
\newcommand{\fk}[1]{\ensuremath{\bm{r}^{#1}}}
\newcommand{\vk}[1]{\ensuremath{\bm{v}^{#1}}}
\newcommand{\psik}[1]{\ensuremath{\bm{\psi}^{#1}}}
\newcommand{\matJk}[1]{\ensuremath{\left[\fnc{J}^{#1}\right]}}
\newcommand{\matJkbar}[1]{\ensuremath{\left[\overline{\fnc{J}}^{#1}\right]}}
\newcommand{\matJxilxmk}[3]{\ensuremath{\left[\Jxilxm{#1}{#2}\right]_{#3}}}
\newcommand{\matJxilxmkbar}[3]{\ensuremath{\overline{\left[\Jxilxm{#1}{#2}\right]}_{#3}}}
\newcommand{\matVtildelmk}[3]{\ensuremath{\left[\Vtildemn{#1}{#2}\right]_{#3}}}
\newcommand{\Thetajk}[2]{\ensuremath{\bm{\Theta}_{#1}^{#2}}}

\newcommand{\wki}[2]{\ensuremath{\bm{w}^{#1}\left(#2,:\right)}}
\newcommand{\uki}[2]{\ensuremath{\bm{u}^{#1}\left(#2,:\right)}}
\newcommand{\vki}[2]{\ensuremath{\bm{v}^{#1}\left(#2\right)}}
\newcommand{\psiki}[2]{\ensuremath{\bm{\psi}^{#1}\left(#2\right)}}
\newcommand{\USC}[4]{\ensuremath{\bfnc{F}^{lc}\left(\uki{#1}{#2},\uki{#3}{#4}\right)}}
\newcommand{\matUSC}[2]{\ensuremath{\mat{F}^{lc}\left(\uk{#1},\uk{#2}\right)}}
\newcommand{\ones}[0]{\ensuremath{\bm{1}}}
\newcommand{\onesbar}[0]{\ensuremath{\overline{\bm{1}}}}
\newcommand{\matdutildedwjk}[2]{\ensuremath{\left[\left|\frac{\partial\tilde{\bmU_{#1}}}{\partial\W}\right|\right]_{#2}}}
\newcommand{\matFSC}[3]{\ensuremath{\mat{F}^{lc}_{x_{#1}}\left(\uk{#2},\uk{#3}\right)}}
\newcommand{\FSCempty}[0]{\ensuremath{\bm{F}^{lc}}}

\newcommand{\URRK}[0]{\ensuremath{\bm{U}}}

\newcommand{\yRRK}[0]{\ensuremath{\bm{Y}}}
\newcommand{\fRRK}[0]{\ensuremath{\bm{F}}}
\newcommand{\dRRK}[0]{\ensuremath{\bm{d}}}
\newcommand{\PhiRRK}[0]{\ensuremath{\bm{\Phi}}}

\newcommand{\MJ}[0]{\ensuremath{\M_{\fnc{J}}}}
\newcommand{\MJbar}[0]{\ensuremath{\overline{\M}_{\fnc{J}}}}
\newcommand{\Dxms}[1]{\ensuremath{\mat{D}_{\xm{#1}}^{s}}}
\newcommand{\Dxmsk}[2]{\ensuremath{\mat{D}_{\xm{#1}}^{s,#2}}}
\newcommand{\Qxmsk}[2]{\ensuremath{\mat{Q}_{\xm{#1}}^{s,#2}}}
\newcommand{\Qxmskbar}[2]{\ensuremath{\overline{\mat{Q}}_{\xm{#1}}^{s,#2}}}
\newcommand{\Sxmsk}[2]{\ensuremath{\mat{S}_{\xm{#1}}^{s,#2}}}
\newcommand{\Qxms}[1]{\ensuremath{\mat{Q}_{\xm{#1}}^{s}}}
\newcommand{\Sxms}[1]{\ensuremath{\mat{S}_{\xm{#1}}^{s}}}
\newcommand{\Exm}[1]{\ensuremath{\mat{E}_{\xm{#1}}}}
\newcommand{\Dxmsc}[1]{\ensuremath{\mat{D}_{\xm{#1}}^{sc}}}
\newcommand{\Dxmsck}[2]{\ensuremath{\mat{D}_{\xm{#1}}^{sc,#2}}}

\newcommand{\Sxmsck}[2]{\ensuremath{\mat{S}_{\xm{#1}}^{sc,#2}}}
\newcommand{\Dxmc}[1]{\ensuremath{\mat{D}_{\xm{#1}}^{c}}}
\newcommand{\Dxmck}[2]{\ensuremath{\mat{D}_{\xm{#1}}^{c,#2}}}
\newcommand{\Exmk}[2]{\ensuremath{\mat{E}_{\xm{#1}}^{#2}}}
\newcommand{\Exmkbar}[2]{\ensuremath{\overline{\mat{E}}_{\xm{#1}}^{#2}}}
\newcommand{\matVhatmnk}[3]{\ensuremath{\left[\Vhatmn{#1}{#2}^{#3}\right]}}
\newcommand{\Rxis}[1]{\ensuremath{\mat{R}_{\xil{#1}}}}

\newcommand{\onek}[0]{\ensuremath{\bm{1}}}
\newcommand{\Mk}[0]{\ensuremath{\mat{P}}_{j}}
\newcommand{\matJkmod}[0]{\ensuremath{\mat{J}}_{j}}

\UseRawInputEncoding
\usepackage{float}
\usepackage{booktabs} 
\newsavebox{\measurebox}
\usepackage{hyperref}
\hypersetup{
 colorlinks=true,
 citecolor=blue,
 linkcolor=blue,
 urlcolor=blue}
\usepackage{tikz}
\usetikzlibrary{matrix,backgrounds,shapes}
\usetikzlibrary{positioning}
\usetikzlibrary{fit}
\usetikzlibrary{plotmarks}
\usetikzlibrary{decorations.pathreplacing}
\usepackage{pgfplots}
\usetikzlibrary{shapes}
\usetikzlibrary{positioning,arrows}
\usetikzlibrary{fadings,shapes.arrows,shadows}
\usetikzlibrary{arrows.meta}
\usepackage{tikz,graphics,color,fullpage,float,epsf,caption,subcaption}
\usepackage{appendix}

\begin{document}

\title[Article Title]{Fully-discrete provably Lyapunov consistent discretizations
for convection-diffusion-reaction PDE systems} 


\author*[1,]{\fnm{Rasha} \sur{Al Jahdali}}\email{rasha.aljahdali@kaust.edu.sa}

\author[2]{\fnm{David C.~Del Rey} \sur{Fern{\'a}ndez}}\email{ddelreyfernandez@uwaterloo.ca}

\author[1]{\fnm{Lisandro} \sur{Dalcin}}\email{dalcinl@gmail.com}

\author[1,3]{\fnm{Matteo} \sur{Parsani}}\email{matteo.parsani@kaust.edu.sa}

\affil[1]{\orgdiv{Computer Electrical and Mathematical Science and Engineering Division (CEMSE)}, \orgname{ King Abdullah University of Science and Technology (KAUST)}, \orgaddress{ \city{Thuwal}, \postcode{23955-6900},  \country{Saudi Arabia}}}

\affil[2]{\orgdiv{ Department of Applied Mathematics}, \orgname{University of Waterloo}, \orgaddress{ \city{Waterloo}
, \country{Canada}}}

\affil[3]{\orgdiv{Physical Science and Engineering Division (PSE)}, \orgname{ King Abdullah University of Science and Technology (KAUST)}, \orgaddress{ \city{Thuwal}, \postcode{23955-6900},  \country{Saudi Arabia}}}

\abstract{
Convection-diffusion-reaction equations are a class of second-order partial 
differential equations widely used to model phenomena involving the 
change of concentration/population of one or more substances/species distributed
in space.
Understanding and preserving their stability properties in numerical simulation is crucial for accurate predictions, system analysis, and decision-making.
This work presents a comprehensive framework for constructing fully discrete 
Lyapunov-consistent discretizations of any order for convection-diffusion-reaction models. 
We introduce a systematic methodology for constructing discretizations that mimic 
the stability analysis of the continuous model using Lyapunov's direct method. 
The spatial algorithms are based on collocated discontinuous Galerkin methods with the summation-by-parts property and the 
simultaneous approximation terms approach for imposing interface coupling and boundary conditions.
Relaxation Runge-Kutta schemes are used to integrate in time and achieve fully
discrete Lyapunov consistency. To
verify the properties of the new schemes, we numerically solve a system of 
convection-diffusion-reaction partial differential equations governing the dynamic evolution of monomer
and dimer concentrations during the dimerization process. Numerical results
demonstrated the accuracy and consistency of the proposed discretizations. 
The new framework can enable further 
advancements in the analysis, control, and understanding of
general convection-diffusion-reaction systems.
}




\maketitle

\section{Introduction}\label{sec1}
Convection-diffusion-reaction partial differential equations (PDEs) are a class of time-dependent second-order PDEs 
that can describe the dynamics of substance concentrations in a medium under the influence of convection, 
diffusion, and reaction processes. These PDEs are often used to describe different phenomena, 
including fluid dynamics, plasma physics, biological
population genetics, neurology, heat transfer, combustion, and reaction chemistry.

Stability analysis is fundamental to understanding the dynamics governing convection-diffusion-reaction models (e.g., stability of equilibrium points and input-output stability).  
Lyapunov's second method, also known as Lyapunov's direct method, is a rigorous approach to analyzing the asymptotic behavior and stability of dynamical systems modeled using ordinary differential equations (ODEs).
This approach uses a generalized notion of energy functions (i.e., Lyapunov functions) to address stability.

Laypunov's direct method has been successfully used in the PDE context to, 
for example, study the stability of fixed points and the stabilization of 
distributed parameter PDEs
\cite{KRSTIC2008750,burton2013stability,GREENBERG198466,komornik1994exact,krstic2008output,smyshlyaev2010boundary,BUONOMO20104010}.
In the PDEs context, the conventional Lyapunov function is replaced by a
Lyapunov functional which is typically constructed by integrating an appropriate
Lyapunov function over the spatial domain. However, generating Lyapunov functionals
is nontrivial. Nevertheless, many studies have made significant progress in
identifying and constructing Lyapunov functionals for specific classes of PDEs,
thereby contributing to understanding stability properties and enabling further 
advancements in control and system analysis.

Recently, Al Jahdali and collaborators \cite{aljahdali2024brain} have 
developed discretizations for (parabolic) reaction-diffusion PDEs that mimic the 
Lyapunov stability properties of the continuous models. The construction of those 
discretizations is rooted in Lyapunov's direct method. Hence, they are referred to as  
``Lyapunov consistent" discretizations. 
In this work, starting from the framework presented in
\cite{aljahdali2024brain}, we develop fully discrete Lyapunov consistent discretizations
that preserve the stability properties of continuous
convection-reaction-diffusion PDE models. 
To design Lyapunov consistent schemes for this class of PDEs, the analysis of the 
convective term requires additional considerations. Here, we develop appropriate 
two-point flux functions \cite{tadmor2003entropy} that render convective 
terms and result in Lyapunov consistent algorithms.

In this work, we present a general framework that enables the construction of 
fully discrete Lyapunov consistent discretizations that are high-order in space and time. 
The algorithms are devised to ensure provable stability, emulating the continuous 
stability proofs used in Lyapunov's direct method. In the stability analysis 
of PDE models, integration by parts (IBP) often plays a vital role. Thus, it is 
natural to leverage a class of spatial operators that discretely mimic IBP. 
Within this class, known as the summation-by-parts (SBP) methods \cite{Svard2014,Fernandez2014},
the discrete analog of IBP serves as a fundamental design principle.
At the spatial level, the proposed methods are constructed using the SBP approach. 
Specifically, the one-dimensional baseline spatial discretizations are collocated
discontinuous Galerkin (DG) methods with the SBP property built at the 
Legendre--Gauss--Lobatto (LGL) quadrature points  \cite{Carpenter2014, Parsani2016}. 
The DG-SBP operators are then extended to multiple dimensions through tensor products. 
The semi-discrete Lyapunov consistency across the
computational domain is achieved by combining the SBP operators with the
simultaneous approximation terms (SAT) approach for the inter-element coupling 
and the enforcement of boundary conditions. The SAT methodology was originally 
developed in the context of linear stability 
\cite{Carpenter1994, Carpenter1999,  Nordstrom1999, mattsson2003boundary} 
and subsequently extended to address nonlinear (entropy) stability  
\cite{Parsani2015b, fernandez_entropy_stable_hp_ref_snpdea_2019}. 
Additionally, the methodology described in \cite{fernandez_entropy_stable_hp_ref_snpdea_2019} 
is used to apply these DG-SBP-SAT discretizations to curvilinear grids with 
geometric ($h$) and polynomial ($p$) 
refinements. To ensure fully discrete Lyapunov consistency and avoid spurious numerical solutions, 
the time integration procedure is performed using relaxation Runge--Kutta (RK) schemes
\cite{ketcheson2019relaxation,ranocha2019relaxation}.
These temporal integration techniques can incorporate properties such as 
conservation, dissipation, or other solution properties associated with convex 
functionals, such as  Lyapunov functionals. This remarkable feature is 
accomplished by multiplying the Runge--Kutta update at each time step by a 
relaxation parameter. 

To demonstrate the properties of the new numerical framework, we solve a system
of convection-diffusion-reaction PDEs that describe the dynamic evolution of
monomer and dimer concentrations during the dimerization process. The algorithms
are implemented in the high-performance computing SSDC framework \cite{PARSANI2021109844},
specifically designed and optimized to exploit spatial and temporal discretization properties and address the scalability and robustness requirements for complex simulations.

The paper is organized as follows. In Section \ref{sec:lyapunov-consistency}, in
the framework of Lyapunov's direct method, we introduce and discuss the concept
of ``Lyapunov consistency." Section \ref{sec:equilibrium-ode-pde} briefly discusses that the PDE extensions of ODE models inherit the ODE model's equilibrium points. 
In Section \ref{sec:dg-sbp-sat-global}, the spatial discretization framework for developing Lyapunov consistent collocated discontinuous Galerkin schemes with the summation-by-part property is presented. Section \ref{section-rrk-lyapunov} summarizes the results
of the relaxation Runge--Kutta schemes, which allow the design of fully
discrete Lyapunov consistent algorithms. In Section \ref{sec:numerical-results}, extensive numerical results, including convergence studies, are presented to demonstrate the accuracy and stability properties of the proposed algorithms. Finally, conclusions are
drawn in Section \ref{sec:discussion-conclusions}.  

\section{Lyapunov Consistency} \label{sec:lyapunov-consistency}

In this work, we consider a system of nonlinear convection-diffusion-reaction
equations in $dim$ spatial dimensions having $r$ species. This system of PDEs
written in conservation (or divergence) form reads 
\begin{equation}\label{eq:SEIR_con_diff}
\frac{\partial \bmU}{\partial t}  =  \bfnc{R}
    - \sum\limits_{l=1}^{dim} \frac{\partial  \bmfI{l}
}{\partial\xm{l}} +\sum\limits_{l=1}^{dim}
\frac{\partial  \bmfV{l}}{\partial
    \xm{l}},\quad
    \bm{x}\in \Omega, \quad t \in (t_0, T_f], 
\end{equation}
where $\bmU \in \mathbb{R}^{r}$, 
$\bmf = \bmf\left(\bmU\right)$ contains the reaction terms, and
$\bmfI{l} =
\bmfI{l}\left(\bmU\right) \in \mathbb{R}^{r \times dim}$ and $\bmfV{l} =
\bmfV{l}\left(\bmU, \frac{\partial \bmU}{\partial \xm{l}}\right)  \in \mathbb{R}^{r \times dim}$ are the convective and diffusive
fluxes in the $l^{\text{th}}$ coordinate direction, respectively.
The position vector is defined as
$\bm{x}=\left(\xm{1},\dots,\xm{dim}\right)\Tr\in \Omega$,
where $\Omega \in \mathbb{R}^{dim}$ is the spatial domain with boundary
$\Gamma$. The symbol $T_f$ indicates the final time, i.e., the upper bound of the
time interval. 
The viscous flux components in \eqref{eq:SEIR_con_diff} are defined as 
\begin{equation}
\bmfV{l} = \Vmn{l}{m}\frac{\partial\bmU}{\partial\xm{m}},
\end{equation}
where $ \Vmn{l}{m}$ is the $r \times r$ matrix of the diffusion coefficients for the
conservative variable, $\bmU$.  
Finally, because system \eqref{eq:SEIR_con_diff} corresponds to an
initial boundary value problem, we have to close it with an
initial condition, 
$\bmU(\bm{x},t_0)=\bmU^{0}$, and appropriate boundary conditions.  

In what follows, we derive the time rate of change of 
the Lyapunov functional, $\tilde{\Lyap}$, which represents a generalized energy function for system \eqref{eq:SEIR_con_diff} 
when it can be constructed. As shown later in the section, the time derivative of $\tilde{\Lyap}$ can be used in 
conjunction with LaSalle's invariance principle  \cite{lasalle1961stability,lasalle1976stability}
to investigate the stability of the fixed points of several PDE models, which can be rewritten as \eqref{eq:SEIR_con_diff}. 
Many of the steps used in the procedure described next
follow the landmark works of Harten \cite{HARTEN1983151} and Tadmor
\cite{TADMOR1984428,tadmor2003entropy} for systems of conservation laws.
These steps also resemble those used in the entropy analysis of the compressible Euler and 
Navier--Stokes equations.
It is important to highlight that the advancement of mathematical theory for 
nonlinear conservation law systems has
been documented over the years in various monographs and books. This journey 
began with the seminal work of Courant and Hilbert in 1962 \cite{courant2008methods}, 
followed by contributions from Lax in 1973 \cite{lax1973hyperbolic}, Smoller in 1983 
\cite{smoller2012shock}, Whitham in 1999 \cite{whitham2011linear}, 
Serre in 2000 \cite{serre1999systems}, Bressan in 2000 \cite{bressan2000hyperbolic}, 
and Dafermos in 2016 \cite{dafermos2016hyperbolic}.

We begin by assuming that a
convex Lyapunov function \cite{khalil_book_2002}, $\Lyap = \Lyap(\bmU)$, for
system \eqref{eq:SEIR_con_diff} exists and can be explicitly generated. The
convexity of the Lyapunov function means that its Hessian is positive definite, i.e.,
\begin{equation}
    {\bm\zeta}^T \, \frac{\partial^2 \Lyap}{\partial \bmU^2} \, {\bm\zeta} > 0, \,
    \, \text{for all nonzero vectors} \, \bm\zeta .
\end{equation}
Often, the 
Lyapunov functional, $\tilde{\Lyap}$, is constructed as
\begin{equation}
\tilde{\Lyap}\equiv\int_{\Omega}\Lyap\mr{d}\Omega.
\end{equation}

Next, we define the Lyapunov
variables (i.e., the entropy variables in \cite{HARTEN1983151,TADMOR1984428}) as
\begin{equation}
    \bm{\fnc{W}}\equiv\left(\frac{\partial\fnc{\Lyap}}{\partial\bmU}\right)\Tr,
\end{equation}
where the convexity of $\Lyap$ 
ensures that there is a one-to-one correspondence between the (state) conservative variables,
$\bmU$, and the 
Lyapunov variables, $\W$. 
Because of the one-to-one mapping, the variables $\bmU$ can be written as a
function of the Lyapunov variables, i.e., $\bmU = \bmU\left(\W\right)$.

At this point, it is useful to
recast the diffusion contribution in system \eqref{eq:SEIR_con_diff}, $\bmfV{l}$, in terms of the Lyapunov
variables, i.e., 
\begin{equation}
\bmfV{l} = \Vhatmn{l}{m}\frac{\partial\bmW}{\partial\xm{m}},
\end{equation}
where $\Vhatmn{l}{m}\equiv\Vmn{l}{m}\frac{\partial\bmU}{\partial \W}$ is the $r
\times r$ matrix of the diffusion coefficients in the new set of variables, $\W$.
Thus, system \eqref{eq:SEIR_con_diff} can be rewritten as 
\begin{equation}\label{eq:SEIR_con_diffW}
\frac{\partial \bmU}{\partial t}  =  \bfnc{R}
- \sum\limits_{l=1}^{dim} \frac{\partial  \bmfI{l}}{\partial\xm{l}} 
    +\sum\limits_{l,m=1}^{dim}
\frac{\partial}{\partial
    \xm{l}}\left(\Vhatmn{l}{m}\frac{\partial\W}{\partial\xm{m}}\right).
\end{equation}
Furthermore, we assume that the following matrix:
\begin{equation}\label{eq:mhat}
\widehat{\mat{M}}\equiv\left[
  \begin{array}{ccc}
    \Vhatmn{1}{1}&\dots&\Vhatmn{1}{dim}\\
    \vdots&&\vdots\\
\Vhatmn{dim}{1}&\dots&\Vhatmn{dim}{dim}\\
  \end{array}
\right]
\end{equation}
is such that $\widehat{\mat{M}}+\widehat{\mat{M}}\Tr$ is positive semi-definite,
\ie, $\widehat{\mat{M}}+\widehat{\mat{M}}\Tr\ge 0$.

\begin{remark}
  The assumption that $\widehat{\mat{M}}+\widehat{\mat{M}}\Tr\ge 0$ is the major assumption in this paper and allows 
  us to completely decouple the analysis of the diffusion from that
  of the convection and reaction terms.
\end{remark}

To obtain the time derivative of
the Lyapunov functional, $d\tilde{\Lyap}/dt$, we proceed as follows.  
We multiply from the left Equation \eqref{eq:SEIR_con_diffW} by $\W\Tr$, which results in
\begin{equation}\label{eq:SEIR_con_diffW-e1}
\W\Tr\frac{\partial \bmU}{\partial t}  =  \W\Tr \bfnc{R}
- \sum\limits_{l=1}^{dim} \W\Tr \frac{\partial \bmfI{l}}{\partial\xm{l}} 
+  \sum\limits_{l,m=1}^{dim} \W\Tr
\frac{\partial}{\partial
    \xm{l}}\left(\Vhatmn{l}{m}\frac{\partial\W}{\partial\xm{m}}\right).
\end{equation}
Then, we integrate over the domain $\Omega$:
\begin{equation}\label{eq:SEIR_con_diffW-e2}
\int_{\Omega}\W\Tr\frac{\partial \bmU}{\partial t} \, \mr{d}\Omega  =
    \int_{\Omega} \W\Tr \bfnc{R} \, \mr{d}\Omega
- \int_{\Omega}\sum\limits_{l=1}^{dim} \W\Tr \frac{\partial
    \bmfI{l}}{\partial\xm{l}} \, \mr{d}\Omega 
+  \int_{\Omega}\sum\limits_{l,m=1}^{dim} \W\Tr
\frac{\partial}{\partial
    \xm{l}}\left(\Vhatmn{l}{m}\frac{\partial\W}{\partial\xm{m}}\right) \, \mr{d}\Omega.
\end{equation}
The temporal term reduces as follows:
\begin{equation}\label{eq:step2_diff}
  \int_{\Omega} \W\Tr \, \frac{\partial \bmU}{\partial t} \, \mr{d}\Omega = 
  \int_{\Omega}\frac{\partial\Lyap}{\partial\bmU}\frac{\partial\bmU}{\partial t}
    \, \mr{d}\Omega = 
  \int_{\Omega}\frac{\partial\Lyap}{\partial t} \, \mr{d}\Omega =  \frac{\mr{d}}{\mr{d}t}
    \int_{\Omega} \Lyap \, \mr{d}\Omega,
\end{equation}
where in the last step we used the Leibniz's rule.
The viscous terms is rewritten using the multidimensional IBP rule as follows:
\begin{equation}\label{eq:step3_diff}
\begin{split}    
  \int_{\Omega}
    \sum\limits_{l,m=1}^{dim}\W\Tr\frac{\partial}{\partial\xm{l}}\left(\Vhatmn{l}{m}\frac{\partial\W}{\partial\xm{m}}\right)
    \, \mr{d}\Omega
    & = \oint_{\Gamma}
     \sum\limits_{l,m=1}^{dim}
     \W\Tr \, \Vhatmn{l}{m}\frac{\partial\W}{\partial\xm{m}} \, \nxm{l} \,
     \mr{d}\Gamma - 
     \int_{\Omega} \sum\limits_{l,m=1}^{dim} \frac{\partial\W}{\partial\xm{l}}\Tr \, \Vhatmn{l}{m} \,
    \frac{\partial\W}{\partial\xm{m}} \, \mr{d}\Omega \\ 
    &= \oint_{\Gamma}
     \sum\limits_{l,m=1}^{dim}
     \W\Tr \, \Vhatmn{l}{m} \, \frac{\partial\W}{\partial\xm{m}} \, \nxm{l} \,
     \mr{d}\Gamma - 
     \int_{\Omega} \bfnc{Z}\Tr \, \hat{\mat{M}}\bfnc{Z} \, \mr{d}\Omega,
\end{split}    
\end{equation}
where $\nxm{l}$ is the component of the outward facing normal on the
boundary $\Gamma$ in the direction
$\xm{l}$,
$\bfnc{Z}\equiv\left[\partial\W/\partial\xm{1},\dots\partial\W/\partial\xm{dim}\right]\Tr$,
and $\hat{\mat{M}}$ is the matrix defined in Equation \eqref{eq:mhat}. Note that one
could obtain expression \eqref{eq:step3_diff} by using the product rule first and then the divergence theorem
\cite{aljahdali2024brain}. 
The convective contribution is manipulated as follows:
\begin{equation}\label{eq:step4_diff}
    \begin{split}
        \int_{\Omega}\sum\limits_{l=1}^{dim} \W\Tr \frac{\partial  \bmfI{l}}{\partial\xm{l}} \, \mr{d}\Omega  &= 
        \int_{\Omega}\sum\limits_{l=1}^{dim} \frac{\partial\fnc{\Lyap}}{\partial\bmU}
        \frac{\partial \bmfI{l}}{\partial\bmU}  \frac{\partial
        \bmU}{\partial\xm{l}} \, \mr{d}\Omega, \\ 
        &=
        \int_{\Omega}\sum\limits_{l=1}^{dim} \frac{\partial \bmF{l}}{\partial \bmU}    \frac{\partial
        \bmU}{\partial\xm{l}} \, \mr{d}\Omega, \\
        &= 
        \int_{\Omega}\sum\limits_{l=1}^{dim} \frac{\partial
        \bmF{l}}{\partial\xm{l}} \, \mr{d}\Omega, \\
        &= \oint_{\Gamma} \sum\limits_{l=1}^{dim} \bmF{l} \, \nxm{l} \,
        \mr{d}\Gamma,
    \end{split}
\end{equation}
with
\begin{equation}
\frac{\partial \Lyap}{\partial \bmU} \frac{\partial \bmfI{l}}{\partial \bmU} =
\frac{\partial \bmF{l}}{\partial \bmU},
\end{equation}
where the scalar function $\bmF{l} = \bmF{l}(\bmU)$ is the Lyapunov flux
(i.e., the entropy flux in
\cite{HARTEN1983151,tadmor2003entropy}) in the
$\xm{l}$ coordinate direction, corresponding to the Lyapunov function
$\Lyap=\Lyap(\bmU)$. The last expression in \eqref{eq:step4_diff} shows that the
contribution of convective terms to the time rate of change of the Lyapunov
functional is ``only" a boundary term. This feature/property is often referred
to as the ``telescoping" property; see \cite{Fisher2013} where this
terminology was introduced in the context of fluid dynamics.

Using Equations \eqref{eq:step2_diff}-\eqref{eq:step4_diff}, expression
\eqref{eq:SEIR_con_diffW-e2} becomes
\begin{equation}\label{eq:SEIR_con_diffW-final}
 \frac{d\widetilde{V}}{d t}=\frac{\mr{d}}{\mr{d}t}
    \int_{\Omega} \Lyap \, \mr{d}\Omega = \int_{\Omega} \W\Tr \bfnc{R} \,
    \mr{d}\Omega +  \oint_{\Gamma} \left( - \sum\limits_{l=1}^{dim} \bmF{l} +
    \sum\limits_{l,m=1}^{dim}
     \W\Tr \, \Vhatmn{l}{m} \, \frac{\partial\W}{\partial\xm{m}} \right) \nxm{l}
 \mr{d}\Gamma - 
     \int_{\Omega} \bfnc{Z}\Tr \, \hat{\mat{M}}\bfnc{Z} \, \mr{d}\Omega.
\end{equation}
Equation \eqref{eq:SEIR_con_diffW-final} provides the expression of the time rate of
change of the Lyapunov functional, $\widetilde{\Lyap}$, for the PDE model
\eqref{eq:SEIR_con_diff}.

Lyapunov's direct method, when combined with LaSalle's invariance principle, can
be used to analyze the stability of infinite dimensional systems, such as PDE systems,
under the following conditions: 
i) a Lyapunov functional, $\widetilde{\Lyap}$, exists and can be generated for the dynamical system of interest, and ii)
the trajectories of the dynamical system are contained within a compact set \cite{curtain2020introduction,luo2012stability}.
When these two conditions are satisfied, the following result holds:

\vspace{0.5cm}
\begin{theorem}\label{thm:lyap_direct_pde}
  Consider the existence of the equilibrium point $\bmU_{eq}$ of
    system~\eqref{eq:SEIR_con_diffW}, (\ie, a point in phase space such that the right-hand side 
    of system \eqref{eq:SEIR_con_diffW} is zero) where each of its components is
    bounded. Furthermore, assume that the following conditions are met: 
  \begin{enumerate}
    \item $\Lyap$ is convex and locally positive definite,
    \item $\W\Tr\bfnc{R}\leq 0$,
    \item $\widehat{\mat{M}}+\widehat{\mat{M}}\Tr\ge 0$,
    \item Appropriate boundary conditions can be found such that 
    \begin{equation*}
    \oint_{\Gamma} \left( - \sum\limits_{l=1}^{dim} \bmF{l} +
    \sum\limits_{l,m=1}^{dim}
     \W\Tr \, \Vhatmn{l}{m} \, \frac{\partial\W}{\partial\xm{m}} \right) \nxm{l}
 \mr{d}\Gamma \leq 0 \, .
    \end{equation*}
  \end{enumerate}
 Then, defining the Lyapunov functional 
  \begin{equation}\label{eq:v_functional}
    \tilde{\Lyap}\equiv\int_{\Omega}\Lyap\mr{d}\Omega,
  \end{equation}
  we have the following stability result:
  \begin{itemize}
    \item  If $\tilde{\Lyap}$ is locally positive definite,  $d\tilde{\Lyap}/dt$ is negative
      definite, and $||\bmU|| \rightarrow +\infty \Rightarrow \tilde{\Lyap}(\bmU) \rightarrow + \infty$, then 
	$\bmU_{eq}=0$ is globally asymptotically stable.
     \end{itemize}
\end{theorem}

\begin{proof} 
 Under the stated conditions equation~\eqref{eq:SEIR_con_diffW-final} reduces to
  \begin{equation*}
    \frac{\mr{d}\widetilde{\Lyap}}{\mr{d}t}\leq 0.
    \end{equation*}
    The Lyapunov functional, $\widetilde{\Lyap}$, is locally positive since $\Lyap$ is locally positive 
    and $\frac{\mr{d}\widetilde{\Lyap}}{\mr{d}t}=0$ at the equilibrium point, $\bmU_{eq}$. 
  Thus, the Lyapunov function, $\Lyap$, can be converted into a 
  Lyapunov functional, $\widetilde{\Lyap}$, with the same properties. Then, 
    we apply the LaSalle's invariance principle 
  to prove the stability statement for the fix point $\bmU_{eq}$. 
\end{proof}

In this work, we develop a systematic framework to design fully discrete algorithms for system \eqref{eq:SEIR_con_diffW}, which mimic at the discrete level the properties
of the Lyapunov functional, $\tilde{\Lyap}$, and its time derivative,
$d\tilde{\Lyap}/dt$, listed in Theorem \ref{thm:lyap_direct_pde}. These algorithms are referred to as
Lyapunov-consistent discretizations and mimic term-by-term each contribution appearing in  
equation \eqref{eq:SEIR_con_diffW-final}.

\section{From ODEs to PDEs: Equilibrium points}\label{sec:equilibrium-ode-pde}
In some modeling approaches, spatially varying models are built from ODE models 
and a natural questions that arises is if the PDE extension inherits the equilibirium points of the ODE model. 
In the context of convection-diffusion-reaction models considered here, we have the following theorem:
\vspace{0.5cm} 
\begin{theorem}\label{thrm:inherit}
  If the convection-diffusion-reaction model~\eqref{eq:SEIR_con_diff} is the PDE
    extension of the ODE model
\begin{equation*}
    \frac{d \bmU}{d t}  =  \bfnc{R}\left(\bmU\right)
 \quad t \in (t_0, T_f],
\end{equation*}
   then, it inherits the 
  equilibrium points of the ODE model. 
\end{theorem}
\begin{proof}
Since the inherited equilibrium points in the PDE context are constant in space,
the spatial terms in Equation \eqref{eq:SEIR_con_diff} are zero. Since $\bm{R}$ must also be zero at an equilibirum point of the ODE model, we see that the ODE's equilibrium points must be the PDE's equilibrium points.
\end{proof}

\section{Spatial discretization: diffusion and convection}
\label{sec:dg-sbp-sat-global}

In this section, we present our spatial discretization framework. The goal is to develop schemes
that are Lyapunov consistent, i.e., they mimic term-by-term each contribution appearing in  
Equation \eqref{eq:SEIR_con_diffW-final} which governs at the continuous level the evolution of
the spatial integral of the Lyapunov function. To construct such schemes, we rely on the SBP framework~\cite{KREISS1974195,gustafsson2013time,Fernandez2014,Svard2014} to derive discrete operators which approximate derivatives 
and SATs to enforce inter-element
coupling and boundary conditions weakly \cite{Carpenter1994, Nordstrom1999, Nordstrom2001b, 
mattsson2003boundary, svard2008stable,Mattsson2009,Parsani2015,Fernandez2018SAT,Yang2018,Duru2019}.
SBP operators are built to discretely mimic the IBP rule, while SATs are used to extend this property over the entire mesh,
including the 
imposition of boundary conditions. SBP operators are typically constructed in a fixed 
computational domain and construct approximations to 
Cartesian derivatives over an element. Here, to simplify the presentation of the semi-discretization 
of~\eqref{eq:SEIR_con_diffW} and the ensuing stability analysis, we use the resultant global SBP 
operators from the procedure mentioned above. Their construction is detailed in Appendix~\ref{app:globalSBP}.
We refer the reader to \cite{Fernandez2019_dense} for a discussion on global SBP operators, and \cite{Crean2018,Alund2019} 
for a related discussion on constructing Cartesian SBP operators from computational SBP operators. While 
the global SBP perspective is helpful for presentation and analysis, 
implementing these schemes is generally done using the local element description.
We provide details on the curvilinear element-wise discretization in Appendix~\ref{app:discelement}. 

\subsection{Preliminaries}

When solving PDEs numerically, the physical domain $\Omega$ with boundary
$\Gamma$, characterized by Cartesian coordinates (e.g., in three
dimensions, $\bm{x} = \left(\xm{1},\xm{2},\xm{3}\right)\subset\mathbb{R}^{3}$) is partitioned into a set of $K$ 
non-overlapping elements. The domain of the element with index $\kappa$ is denoted
by $\Omega^{\kappa}$, and its boundary
is denoted by $\Gamma^{\kappa}$. Thus, $\Omega = \bigcup_{\kappa=1}^{K} \Omega^{\kappa}$. 
We assume that $\Omega^{\kappa}$ is piecewise smooth for all
elements $\kappa$. Without loss of generality, this work considers a partition (or
tessellation) of
the domain $\Omega$ accomplished with tensor-product elements, i.e.,
quadrilateral in two dimensions (2D) or hexahedral elements in in three dimensions (3D).
Numerically, we solve PDEs in computational 
coordinates 
$\bm{\xi} = \left(\xil{1},\xil{2},\xil{3}\right)\Tr\subset\mathbb{R}^{3}$, 
where each $\Omega^{\kappa}$ is logically transformed to the reference element
$\widehat{\Omega}^{\kappa}$, with boundary ${\widehat\Gamma}^{\kappa}$,
using a pull-back curvilinear coordinate transformation which satisfies the following assumption:
\vspace{0.5cm}
\begin{assume}\label{assume:curv}
Each element in physical space is transformed using 
a local and invertible curvilinear coordinate transformation that is compatible at 
shared interfaces, meaning that the push-forward element-wise mappings are continuous across physical element interfaces.
Note that this is the standard assumption requiring that the curvilinear coordinate transformation is water-tight.
\end{assume}

\begin{figure}
 \begin{center}
   \includegraphics[width=0.6\textwidth]{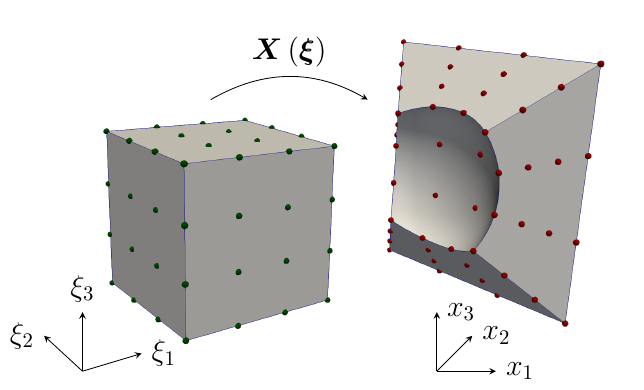}
   \caption{The reference element and its mapping to an element in a
     3D mesh.}
   \label{fig:ref_cell}
  \end{center}
\end{figure}

Precisely, one maps from the reference coordinates
$\bm{\xi} \in [-1,1]^{dim}$
to the physical element (see Figure \ref{fig:ref_cell} for the 3D case) by the push-forward transformation 
\begin{equation}\label{eq:3d_mapping}
    \bm{x} = X \left(\bm{\xi}\right),
\end{equation}
which, in the presence of curved elements, is usually a high-order degree
polynomial.
In our algorithm, we use a cell-wise isoparametric approach based on Lagrange basis
functions. However, any set of polynomial basis functions and corresponding 
dual coefficients could be used to characterize the push-forward mapping \eqref{eq:3d_mapping}. 
As usual in unstructured mesh schemes, the procedure described in what follows does not require
explicit knowledge nor construction of the pull-back mappings.

Since derivatives are approximated with differentiation operators defined in 
computational space, we use the Jacobian of the push-forward mapping and the chain rule
\begin{equation*}
  \frac{\partial}{\partial\xm{l}}=
  \sum\limits_{m=1}^{dim}\frac{\partial\xil{m}}{\partial\xm{l}}\frac{\partial}{\partial\xil{m}},\quad
  \frac{\partial^{2}}{\partial\xm{l}^{2}}=
  \sum\limits_{m,a=1}^{dim}\frac{\partial\xil{m}}{\partial\xm{l}}
  \frac{\partial}{\partial\xil{m}}\left(
  \frac{\partial\xil{a}}{\partial \xm{l}}\frac{\partial}{\partial\xil{a}}  
  \right),
\end{equation*} 
to transform the PDEs from physical, $\bm{x}$ , to computational space, $\bm{\xi}$. From here onward, we denote
the determinant of the metric Jacobian by the symbol $J$.

Herein, we discretize functions using their nodal values.
Suppose we have a set of $N$ nodes placed in the domain $\Omega$ whose
coordinates are stored in $\bm{x}^{g} = \{(x_{i,1}, \ldots, x_{i,dim}
)\}^{N}_{i=1}$. Then, a generic continuous function $U \in L^2
(\Omega)$ evaluated at the nodes is represented using the column vector
\begin{equation*}
    \bm{u}^{g} = \left(U(x_{1,1}, \ldots, x_{1,dim}),
    U(x_{2,1},
    \ldots, x_{2,dim}), \ldots, U(x_{N,1}, \ldots, x_{N,dim})\right)
    \in \mathbb{R}^{rN}.
\end{equation*}
Note that the superscript ``$g$" in $\bm{x}^{g}$ and $\bm{u}^{g}$ stands for
``global". This superscript is used for all the quantities that span the domain
$\Omega$, as opposed to the local quantities and operators associated with each element $\kappa$
with volume $\Omega^{\kappa}$.
We can now present the global SBP operators for the first derivative and their properties. The details 
of each component of these operators are presented in Appendix~\ref{app:globalSBP}.
\vspace{0.5cm}
 \begin{definition}\label{def:SBPD1}
    Consider a domain $\Omega$ partitioned by $K$ non-overlapping elements resulting in a tessellation of the 
    boundary surface $\Gamma$ with $M$ surfaces, where each face belongs
     to a uniquely defined element. Then, the order $p$ global, diagonal norm, skew-symmetric SBP operator
     $\Dxmsk{l}{g}\in\mathbb{R}^{(N\,\times r)\times (N\,\times r)}$ 
    approximating the first derivative $\partial/\partial\xm{l}$ for a system of
     $r$ equations (i.e., $\partial\bmU/\partial\xm{l}$, $\bmU \in \mathbb{R}^{r}$) 
    on the $N$ node nodal distribution $\bm{x}^{g}$ has the following properties.
    \begin{enumerate}
        \item
            $\Dxmsk{l}{g}\, \uk{g}\approx\frac{\partial\bmU}{\partial\xm{l}}\left(\bm{x}^{g}\right)+\mathcal{O}\left(\Delta
            x^{p}\right)$,
        \item $\Dxmsk{l}{g}\equiv\left(\MJ^{g}\right)^{-1}\Qxmsk{l}{g}$, where
            the norm matrix, $\MJ^{g} = \M^{g}$\matJk{}, is diagonal and positive definite,
        \item $\Qxmsk{l}{g}\equiv\Sxmsk{l}{g}+\frac{1}{2}\Exmk{l}{g}$, where $\Sxmsk{l}{g}=-\left(\Sxmsk{l}{g}\right)\Tr$, 
        $\Exmk{l}{g}=\left(\Exmk{l}{g}\right)\Tr$,
        \item the matrix $\Exmk{l}{g}$ is diagonal
      and is constructed as, 
  \begin{equation}\label{eq:Exmgapp}
    \Exmk{l}{g}\equiv\sum\limits_{i=1}^{M}n_{\xi_{\Gamma_i}}\left(\Rxis{\Gamma_i}\right)\Tr\Porthoxil{\Gamma_i}\Rxis{\Gamma_i}
    \matJxilxmk{\Gamma_i}{m}{},
  \end{equation}
  where $\Gamma_{i}$ specifies the computational coordinate orthogonal to the $i^{\text{th}}$ surface (i.e., in three dimensions 
 $\Gamma_{i}$ can assume the value $1$ or $2$ or $3$). The diagonal matrices $\matJxilxmk{i}{l}{i}$ contain the the metric terms 
      along their diagonal. Furthermore, the action of the matrices $\Rxis{\Gamma_i}$ is to 
      pick off the appropriate portions of the solution vector at the $i^{\text{th}}$ 
            boundary surface which belongs to $\Gamma$. 
     Thus, $\Rxis{\Gamma_i}\uk{g}$ results in a vector whose entries are the portion of the global 
      solution vector $\uk{g}$ at the boundary surface with index $i$, which belongs to $\Gamma$. 

            The diagonal matrices $\Porthoxil{\Gamma_i}$ contain 
      cubature weights for the same $i^{\text{th}}$ surface. The diagonal matrix
            $[J]$ used to construct $\MJ^{g}$ is a diagonal matrix with the metric Jacobian, $J$, of each node along its diagonal.
    \end{enumerate}
\end{definition}

The matrices that appear in Definition \ref{def:SBPD1} are closely related to integral bilinear forms;
see \cite{Hicken2016} for the case of multidimensional SBP operators. For
completeness, we report their relationships. The diagonal
matrix $\MJ^{g}$ can be interpreted as a lumped mass matrix and can be used to approximate an inner
product of two functions, e.g., their $L^{2}$ inner product. For instance, 
given two column vectors $\bm{U}, \bm{V}  \in \mathbb{R}^{r}$ whose $r$
components, $U_i, V_i$, are
given by square integrable functions over the domain $\Omega$ (i.e.,
$U_i = U_i(\bm{x}), V_i = V_i(\bm{x})
\in L^2(\Omega)$, $i=1, \ldots, r$) we have:
\begin{equation}\label{eq:L2}
    \left(\vk{g}\right)\Tr \, \MJ^{g} \, \uk{g}\approx\int_{\Omega}\bm{V}\Tr \,
    \bm{U} \, \mr{d}\Omega. 
\end{equation}
Furthermore, for the matrix $\Qxmsk{l}{g}$ we have that
\begin{equation}\label{eq:Qxmsk}
    \left(\vk{g}\right)\Tr \, \Qxmsk{l}{g} \, \uk{g}\approx\int_{\Omega}\bm{V}\Tr \,
    \frac{\partial \bm{U}}{\partial \xm{l}} \, \mr{d}\Omega 
\end{equation}
and
\begin{equation}\label{eq:QxmskT}
    \left(\vk{g}\right)\Tr \, \left(\Qxmsk{l}{g}\right)\Tr \, \uk{g}\approx\int_{\Omega} 
    \left(\frac{\partial \bm{V}}{\partial \xm{l}}\right)\Tr
    \bm{U} \,
     \, \mr{d}\Omega. 
\end{equation}
Finally, we note that the matrix $\Exmk{l}{g}$ can be interpreted as an approximation the 
the following surface integral: 
\begin{equation}\label{eq:Exmg-2}
    \left(\vk{g}\right)\Tr \, \Exmk{l}{g} \, \uk{g} \, \approx\oint_{\Gamma}\bm{V}\Tr \,
    \bm{U} \, \nxm{l} \, \mr{d}\Gamma. 
\end{equation}

Now, we demonstrate how the global SBP operator $\Dxmsk{l}{g}$ mimics the IBP formula, i.e.,
\begin{equation}\label{eq:IBPmain}
  \int_{\Omega}\left(\bm{V}\Tr\frac{\partial\bm{U}}{\partial\xm{l}}+
  \bfnc{\bm{U}}\Tr\frac{\partial\bm{V}}{\partial\xm{l}}
  \right)\mr{d}\Omega=\oint_{\Gamma}\bm{V}\Tr \, \bm{U} \,
    \nxm{l} \,\mr{d}\Gamma.
\end{equation}
Discretizing the left-hand side of~\eqref{eq:IBPmain} using the global SBP
operator $\Dxmsk{l}{g}$ and 
the norm matrix $\MJ^{g}$, both of which are given in Definition \ref{def:SBPD1}, results in the equality 
\begin{equation}\label{eq:SBPmain}
    \left(\vk{g}\right)\Tr \, \MJ^{g} \, \Dxmsk{l}{g} \, \uk{g} +
    \left(\uk{g}\right)\Tr \, \MJ^{g} \,
    \Dxmsk{l}{g}\, \vk{g}=
    \left(\vk{g}\right)\Tr \, \Exmk{l}{g} \, \uk{g} ,
\end{equation}
where each term is an approximation to the corresponding term in the IBP 
formula~\eqref{eq:IBPmain}.

In analogy to developing nonlinearly-stable (entropy-stable) numerical techniques for hyprbolic conservation laws, the convective term in \eqref{eq:SEIR_con_diff} requires special care.
In fact, while $\Dxmsk{l}{g}$ mimics the IBP rule, it does not mimic the 
telescoping property that is needed for the Lyapunov consistency analysis,
namely (see Equation \eqref{eq:step4_diff})
\begin{equation}\label{eq:step4_diff_last}
    \begin{split}
        \int_{\Omega}\sum\limits_{l=1}^{dim} \W\Tr \, \frac{\partial
        \bmfI{l}}{\partial\xm{l}} \, \mr{d}\Omega &= \oint_{\Gamma}
        \sum\limits_{l=1}^{dim} \bmF{l} \, \nxm{l} \,
        \mr{d}\Gamma.
    \end{split}
\end{equation}
To mimic at the semi-discrete level \eqref{eq:step4_diff_last}, we can combine
the SBP operator with Tadmor's two-point flux function \cite{tadmor2003entropy} 
and their extension to SBP discretizations via the Hadammard formalism \cite{Fisher2013b,fernandez_entropy_stable_hp_ref_snpdea_2019}. The discrete operator that we require is given as 
\begin{equation}\label{eq:DU}
    2\Dxmsk{l}{g}\circ\matFSC{l}{g}{g}\bm{1}\approx\frac{\partial\bmfI{l}}{\partial
    \xm{l}}\left(\bm{x}^{g}\right),
\end{equation}
where $\circ$ is the Hadammard product (entry-wise multiplication) between two matrices. 
Moreover $\matFSC{l}{g}{g}$ is a symmetric matrix constructed from two point flux functions \cite{tadmor2003entropy} and 
$\bm{1}$ is an appropriately sized vector of ones. The superscript ``$lc$" in $\matFSC{l}{g}{g}$ stands
for ``Lyapunov conservative".
The details of this discrete operator on the left-hand side of \eqref{eq:DU} are presented in Appendix~\ref{app:Hadammard}.
In the following 
theorem, we state that this operator mimics formula~\eqref{eq:step4_diff_last}.
\vspace{0.5cm}
\begin{theorem}\label{thrm:SBPnon}
    If the matrix $\matFSC{l}{g}{g}$ is constructed from an entropy-consistent two-point flux function 
  that is symmetric and satisfies the Tadmor shuffle condition \cite{tadmor2003entropy} and $\Dxmsk{l}{g}$ is 
    constructed such that $\Dxmsk{l}{g} \, \bm{1}=\bm{0}$ (see Appendix~\ref{app:Hadammard} for more 
  details on both of these points), then for 
  the class of diagonal norm SBP operators with diagonal $\Exmk{l}{g}$ the following holds:
  \begin{equation}
      2\left(\wk{g}\right)\Tr \, \MJ^{g} \, \Dxmsk{l}{g}\circ\matFSC{l}{g}{g} \, \bm{1}=
      \bar{\bm{1}}\Tr \, \Exmkbar{l}{g} \, \mathfrak{f}_{x_l}^{g} ,
  \end{equation}
where the overbar notation means the scalar version of a matrix or a vector, \ie, 
$\Exmk{l}{g}=\Exmkbar{l}{g}\otimes\Imat{r}$, where $r$ is the number of equations. Furthermore, 
the vector $\mathfrak{f}_{x_l}^{g}$ is the vector constructed from evaluating
    the Lyapunov flux $\bmF{l}$ at the mesh nodes.
\end{theorem}
\begin{proof}
  See Appendix~\ref{app:Hadammard} for the proof.
\end{proof}

Next, we present approximations to the variable coefficient second derivative that are SBP. 
To do so, we use two different global matrix difference operators that are not SBP but combined, 
result in an SBP discretization of the second derivative. 
\vspace{0.5cm}
 \begin{definition}\label{def:SBPD2}
  Under the same conditions as in Definition~\ref{def:SBPD1} the SBP approximation 
  to the second derivative used in this paper is given as
 \begin{equation}\label{eq:D2gmain}
    \Dxmsck{l}{g}\matVhatmnk{l}{m}{g}\Dxmck{m}{g}=\left(\MJ^{g}\right)^{-1}
    \left(-\left(\Dxmck{l}{g}\right)\Tr\MJ^{g}\matVhatmnk{l}{m}{g}\Dxmck{m}{g}+
    \Exmk{l}{g}\matVhatmnk{l}{m}{g}\Dxmck{m}{g}\right),
  \end{equation}
  where $\matVhatmnk{l}{m}{g}$ is a block diagonal matrix with the $\Vhatmn{l}{m}$ 
  evaluated at its mesh nodes along its diagonal. Moreover, the matrix difference 
  operator $\Dxmsck{l}{g}$ is called ``the strong-conservation form" approximation 
  to the derivative $\frac{\partial}{\partial\xm{l}}$ and while the matrix 
  difference operator $\Dxmck{l}{g}$ is called ``the chain rule form". Neither of 
  these operators are SBP, and their construction is detailed in Appendix~\ref{app:globalSBP}.
 \end{definition}
 The specific IBP formula that we wish to mimic with our discretization of the second 
 derivative is (see also Equation \eqref{eq:step3_diff})
 \begin{equation}\label{eq:IBPD2main}
  \int_{\Omega}\bfnc{V}\Tr \, \frac{\partial}{\partial\xm{l}} \,
     \left(\Vhatmn{l}{m} \, \frac{\partial\bmU}{\partial\xm{m}}\right)
     \, \mr{d}\Omega = 
  \oint_{\Gamma}\bfnc{V}\Tr \, \Vhatmn{l}{m} \,
     \frac{\partial\bmU}{\partial\xm{m}} \, \nxm{l} \, \mr{d}\Gamma
  -\int_{\Omega}\frac{\partial\bfnc{V}}{\partial\xm{l}}\Tr \, \Vhatmn{l}{m} \, \frac{\partial\bmU}{\partial\xm{l}}
     \, \mr{d}\Omega.
 \end{equation}
 Discretizing the left-hand side of~\eqref{eq:IBPD2main} using our SBP operator and the norm matrix 
 results in the following equality:

 \begin{equation}\label{eq:SBPD2main}
     \left(\vk{g}\right)\Tr \, \Dxmsck{l}{g} \, \MJ^{g} \, \matVhatmnk{l}{m}{g} \, \Dxmck{m}{g}\uk{g}=
 \left(\vk{g}\right)\Tr \, \Exmk{l}{g} \, \matVhatmnk{l}{m}{g} \,
     \Dxmsck{m}{g}\uk{g} \, 
-\left(\vk{g}\right)\Tr \, \left(\Dxmck{l}{g}\right)\Tr \, \MJ^{g} \,
     \matVhatmnk{l}{m}{g} \, \Dxmck{m}{g}\uk{g},
 \end{equation}
 where we see that each term is an approximation to the corresponding term 
 in the IBP formula~\eqref{eq:IBPD2main}. 

We now have all the relevant pieces to construct our semi-discrete form of~\eqref{eq:SEIR_con_diffW},
\begin{equation}\label{eq:SEIR_con_diffWdisc}
  \begin{split}
      \frac{\mr{d} \uk{g}}{\mr{d} t}  =  \, &\fk{g}
      -\sum\limits_{l=1}^{d}2\Dxmsk{l}{g} \circ \matFSC{l}{g}{g} \, \bm{1}
+\sum\limits_{l,m=1}^{d}\Dxmsck{l}{g} \, \matVhatmnk{l}{m}{g} \, \Dxmck{m}{g} \, \uk{g}\\
      &+\bm{SAT}^{(bc)}+\bm{diss}^{(c)}+\bm{diss}^{(d)},
  \end{split}
\end{equation}
where $\bm{SAT}^{(bc)}$ contains the boundary SATs to weakly impose the boundary condition 
and $\bm{diss}^{(c)}$ and $\bm{diss}^{(d)}$ add interface dissipation for the
convective 
and the diffusion terms, respectively; their construction is detailed in Appendix~\ref{app:diss}. 
For our analysis, we note that $\bm{diss}^{(c)}$ and $\bm{diss}^{(d)}$
are constructed such that
\begin{equation}\label{eq:dissprop}
    \left(\wk{g}\right)\Tr \, \MJ^{g} \, \bm{diss}^{(c)}\leq 0,\quad
    \left(\wk{g}\right)\Tr \, \MJ^{g} \, \bm{diss}^{(d)}\leq 0.
\end{equation}
We observe that system \eqref{eq:SEIR_con_diffWdisc} can be recast in the 
following compact and general ODE notation:
\begin{equation}\label{eq:system_odes}
    \frac{d \bmU(t) }{dt} = \bm{F}\left(\bmU(t)\right).
\end{equation}

The proof of the stability of the equilibrium points of the semi-discrete
equations defined by equation~\eqref{eq:SEIR_con_diffWdisc} follows identically 
to that of the continuous equations. Thus, we need to derive
the expression of the time rate of change of the discrete Lyapunov functional. 
The first step is multiplying through by the Lyapunov variables and discretely integrating 
in space, followed by the discrete counterparts to the steps in
equations~\eqref{eq:SEIR_con_diffW-e1} through~\eqref{eq:SEIR_con_diffW-final}. We summarize the result 
in the following theorem, i.e., the analog of Theorem \ref{thm:lyap_direct_pde} presented for the
analysis at the continuous level.

\vspace{0.5cm}
\begin{theorem}\label{thrm:boundLyapDisc}
    Consider the  equilibrium point $\bm{u}_{eq}$ of the semi-discrete equations
    defined by the ODE system~\eqref{eq:SEIR_con_diffWdisc}. 
  Furthermore, assume that the following conditions are met: 
  \begin{enumerate}
    \item $\bmW\Tr\bfnc{R}\leq 0$,
    \item $\Lyap$ is convex and locally positive definite,
    \item $\widehat{\mat{M}}+\widehat{\mat{M}}\Tr\ge 0$,
    \item Appropriate boundary conditions and associated $\bm{SAT}^{(bc)}$ can be found such that
      \begin{equation*}
-\sum\limits_{l=1}^{dim}\bar{\bm{1}}\Tr \, \Exmkbar{l}{g} \,
          \mathfrak{f}_{x_l}^{g} \,
          +\sum\limits_{l,m=1}^{dim}\left(\wk{g}\right)\Tr \, \Exmk{l}{g} \,
          \matVhatmnk{l}{m}{g} \, \Dxmsck{m}{g}\wk{g} \, + \left(\wk{g}\right)\Tr \,
          \MJ^{g} \, \bm{SAT}^{(bc)}
\leq 0.
      \end{equation*} 
  \end{enumerate}

      Then, defining the discrete Lyapunov functional 
\begin{equation}\label{eq:discrete_lyap_functional}
    \tilde{v}^{g}\equiv\sum\limits_{k=1}^{K}\onesbar\Tr \, \Mbar^{k} \, \matJkbar{k} \, \vk{k},
  \end{equation}
    where the vector $\vk{k}$ is the discrete version of the Lyapunov function at the nodes of the
    $k^\text{th}$ element, we have the following stability result:
  \begin{itemize}
    \item  If $\tilde{v}^{g}$ is locally positive definite and $d\tilde{v}^{g}/dt : \mathbb{R}^{r} \rightarrow \mathbb{R}$ is negative
      definite, and $||\uk{g}|| \rightarrow +\infty \Rightarrow \tilde{v}^{g}(\uk{}) \rightarrow + \infty$, then 
	$\bm{u}_{eq}=0$ is globally asymptotically stable.
     \end{itemize}
  \end{theorem}

  \begin{proof}
We mimic the steps at the continuous level and seek to construct a bound on the time-derivative of $\vk{g}$. 
Multiplying~\eqref{eq:SEIR_con_diffWdisc} by $\left(\wk{g}\right)\Tr\MJ^{g}$ results in
\begin{equation}\label{eq:discentropy1}
  \begin{split}
\left(\wk{g}\right)\Tr\MJ^{g}\frac{\mr{d} \uk{g}}{\mr{d} t}  =&  \left(\wk{g}\right)\Tr\MJ^{g}\fk{g}
-\sum\limits_{l=1}^{dim}2\left(\wk{g}\right)\Tr\MJ^{g}\Dxmsk{l}{g}\circ\matFSC{l}{g}{g}\bm{1}\\
&+\sum\limits_{l,m=1}^{dim}\left(\wk{g}\right)\Tr\MJ^{g}\Dxmsck{l}{g}\matVhatmnk{l}{m}{g}\Dxmck{m}{g}\uk{g}\\
&+\left(\wk{g}\right)\Tr\MJ^{g} \, \bm{SAT}^{(bc)}+\left(\wk{g}\right)\Tr\MJ^{g} \, \bm{diss}^{(c)}+\left(\wk{g}\right)\Tr\MJ^{g} \, \bm{diss}^{(d)}.
\end{split}
\end{equation}
The temporal term reduces as follows:
\begin{equation}\label{eq:tempdisc}
  \begin{split}
\left(\wk{g}\right)\Tr\MJ^{g}\frac{\mr{d}\uk{g}}{\mr{d}t} =& 
\sum\limits_{j=1}^{N^{dim}}\MJbar^{g}(j,j)\wki{g}{i}\Tr\frac{\mr{d}\uki{g}{i}}{\mr{d}t}
=\sum\limits_{j=1}^{N^{dim}}\MJbar^{g}(j,j)\frac{\mr{d}\vki{g}{i}}{\mr{d}t}\\
=&\onesbar\Tr \, \MJbar^{g} \, \frac{\mr{d}\vk{g}}{\mr{d}t},
  \end{split}
\end{equation}
      where $\onesbar$ is a vector of ones of size $N^{dim}$ (i.e., the number of nodes in the domain, $N$, to the
      power $dim$) and the bar over $\MJbar^{g}$ 
signifies the $N^{dim}\times N^{dim}$ version of this matrix. Furthermore, the notation 
$\uki{k}{i}$ means the $r\times 1$ portion of $\uk{g}$ at the $i^\text{th}$ node.
The inviscid terms are recast with the help of Theorem~\ref{thrm:SBPnon}, while the 
diffusion terms are recast using~\eqref{eq:D2gmain}. Thus, Equation~\eqref{eq:discentropy1} 
becomes 

\begin{equation}\label{eq:discentropy2}
  \begin{split}
 \onesbar\Tr \, \MJbar^{g} \, \frac{\mr{d}\vk{g}}{\mr{d}t} =&  \left(\wk{g}\right)\Tr \,
      \MJ^{g} \, \fk{g}
-\sum\limits_{l=1}^{dim}\bar{\bm{1}}\Tr \, \Exmkbar{l}{g} \, \mathfrak{f}_{x_l}^{g} \, \\
&+\sum\limits_{l,m=1}^{dim}\left(\wk{g}\right)\Tr
\left(-\left(\Dxmck{l}{g}\right)\Tr \, \MJ^{g} \, \matVhatmnk{l}{m}{g} \, \Dxmck{m}{g}+
    \Exmk{l}{g} \, \matVhatmnk{l}{m}{k} \, \Dxmck{m}{g} \, \right)\wk{g}\\
      &+\left(\wk{g}\right)\Tr \, \MJ^{g} \, \bm{SAT}^{(bc)}+\left(\wk{g}\right)\Tr \,
      \MJ^{g} \, \bm{diss}^{(c)}+\left(\wk{g}\right)\Tr \, \MJ^{g} \,
      \bm{diss}^{(d)}.
\end{split}
\end{equation}
  By assumption $\left(\wk{g}\right)\Tr \, \MJ^{g} \, \fk{g}\leq 0$. Moreover,
      we can recast the firs part of the diffusion term as 
\begin{equation*}
-\sum\limits_{l,m=1}^{d}\left(\wk{g}\right)\Tr\left(\Dxmck{l}{g}\right)\Tr \, \MJ^{g}
    \, \matVhatmnk{l}{m}{g}
    \, \Dxmck{m}{g}
    \, \wk{g} = 
\bm{z}\Tr \, \hat{\mat{M}} \, \bm{z},
\end{equation*}
where
      $\bm{z}\equiv\left[\left(\Dxmck{1}{g}\wk{g}\right)\Tr,\dots,\left(\Dxmck{d}{g}\wk{g}\right)\Tr\right]\Tr$.
      Thus, via the assumptions, we have $\bm{z}\Tr\hat{\mat{M}}\bm{z}\leq 0$.
      In addition, 
      $\left(\wk{g}\right)\Tr\MJ^{g}\bm{diss}^{(c)}\leq 0$ and
      $\left(\wk{g}\right)\Tr\MJ^{g}\bm{diss}^{(d)}\leq 0$. Therefore, dropping 
all negative semi-definite terms, expression \eqref{eq:discentropy2} reduces to 
\begin{equation}\label{eq:discentropy3}
  \begin{split}
 \frac{\mr{d}\tilde{v}^{g}}{\mr{d}t}=\onesbar\Tr\MJbar^{g}\frac{\mr{d}\vk{g}}{\mr{d}t} \leq& 
-\sum\limits_{l=1}^{dim}\bar{\bm{1}}\Tr \, \Exmkbar{l}{g} \,\mathfrak{f}_{x_l}^{g} \,
+\sum\limits_{l,m=1}^{dim}\left(\wk{g}\right)\Tr \, \Exmk{l}{g} \,
      \matVhatmnk{l}{m}{k} \, \Dxmsck{m}{g} \, \wk{g} \, + \left(\wk{g}\right)\Tr \,
      \MJ^{g} \, \bm{SAT}^{(bc)}.
\end{split}
\end{equation}

The final results arises from the assumption on the interplay between terms on the right-hand 
side of~\eqref{eq:discentropy3}. Moreover, $\tilde{v}^{g}$ is locally positive since ${v}^{g}$ is locally positive 
    and $\frac{\mr{d}\tilde{v}^{g}}{\mr{d}t}=0$ at the equilibrium point.    
  Thus, the discrete Lyapunov function, ${v}^{g}$, can be converted into a 
  Lyapunov functional, $\tilde{v}^{g}$, with the same properties. Therefore, we can
  apply Lyapunov's direct method and LaSalle's invariance principle \cite{lasalle1961stability,lasalle1976stability} 
  to prove the stability statement for $\bm{u}_{eq}$.
  \end{proof}

\section{Fully discrete schemes} \label{section-rrk-lyapunov}
The semi-discrete scheme \eqref{eq:discentropy3} retains the stability
properties of the original system of PDEs; however, these can be easily lost when 
discretizing the temporal term.
Here, we consider a broader class of temporal discretization schemes that 
allow us to perform numerical time integration while retaining the stability 
properties of \eqref{eq:discentropy3}.
Specifically, we use the explicit class of relaxation Runge--Kutta schemes 
\cite{ketcheson2019relaxation,ranocha2019relaxation}. 
Relaxation Runge--Kutta schemes can be used to enforce conservation, dissipation, or other
solution properties with respect to any convex function by adding a relaxation parameter 
that multiplies the Runge--Kutta update at each time step. The computational overhead is the 
solution of one additional nonlinear scalar algebraic equation for which a good initial guess is available. The analysis of
these schemes in the context of Lyapunov consistent discretizations was first
presented in \cite{aljahdali2024brain}. Here, for completeness, we report the main results.

%
%
\subsection{Relaxation Runge--Kutta methods}

A general (explicit or implicit) Runge--Kutta method with $s$ stages can be represented
by its Butcher tableau \cite{butcher2008numerical},
\begin{equation}
\label{eq:butcher}
\begin{array}{c|c}
\bm{c} & \mat{A} \\ \hline & \bm{b}\Tr
\end{array}\, ,
\end{equation}
where $\mat{A} \in \mathbb{R}^{s \times s}$ and $\bm{b}, \bm{c} \in \mathbb{R}^s$.
For the system of ODEs \eqref{eq:SEIR_con_diffWdisc}, a step from $\URRK^n \approx \URRK(t_n)$ to
$\URRK^{n+1} \approx \URRK(t_{n+1})$, where $t_{n+1} = t_n + \Delta t$, is given by
\begin{subequations}
\label{eq:RK-step}
\begin{align}
\label{eq:RK-stages}
\yRRK_i &= \URRK^n + \Delta t \sum_{j=1}^{s} a_{ij} \, \fRRK(t_n + c_j \Delta t, \yRRK_j), \qquad i = 1, \dots, s,
\\
\label{eq:RK-final}
\URRK^{n+1}&=\URRK^n + \Delta t \sum_{i=1}^{s} b_{i} \, \fRRK(t_n + c_i \Delta t, \yRRK_i).
\end{align}
\end{subequations}
Here, $\yRRK_i$ are the stage values of the Runge--Kutta method, and $\fRRK$ is
the right-hand side of the ODE system \eqref{eq:SEIR_con_diffWdisc}.
To prevent confusion regarding the notation for these vectors of continuous functions, we note that 
uppercase $\URRK$, and $\fRRK$ are used for Runge--Kutta method in this section,
as often used in the 
literature in our previous work \cite{aljahdali2024brain}. In addition, for simplicity, we will make use of the shorthand
\begin{align}
\fRRK_i & \equiv \fRRK(t_n + c_i \Delta t, \yRRK_i), \quad \fRRK_0=\fRRK(t_n, \URRK^n).
\end{align}

The basic idea to make a given Runge--Kutta method Lyapunov
consistent (\ie, entropy stable in the framework of the compressible Euler 
and Navier--Stokes equations \cite{ranocha2019relaxation}) 
is to replace the update \eqref{eq:RK-final}
with 
\begin{equation}
\label{eq:RK-final-gamma}
\URRK^{n+1}_\gamma= \URRK^n + \gamma_n \Delta t \sum_{i=1}^{s} b_{i} \, \fRRK_i,
\end{equation}
where $\gamma_n$ is a real number called the relaxation parameter. The
generalization to Lyapunov consistency is to enforce the condition
\begin{equation}
  \label{eq:rrk-condition}
\Lyap(\URRK^{n+1}_\gamma) - \Lyap(\URRK^n)=\gamma_n \Delta t \sum_{i=1}^s b_i \scp{\Lyap'(\yRRK_i)}{\fRRK_i},
\end{equation}
by finding a root $\gamma_n$ of
\begin{equation}
\label{eq:r}
q(\gamma)=\Lyap\biggl( \URRK^n + \gamma \Delta t \sum_{i=1}^{s} b_{i} \fRRK_i \biggr)-\Lyap(\URRK^n)
-\gamma \Delta t \sum_{i=1}^s b_i \scp{\Lyap'(\yRRK_i)}{\fRRK_i},
\end{equation}
where $\Lyap$ is the Lyapunov function or functional, and $\scp{\cdot}{\cdot}$ denote the inner product
inducing the norm $\| \cdot \|$. Note that $\Lyap'(\yRRK_i)=\bmW(\yRRK_i)$. Moreover, the direction
\begin{equation}
\label{eq:d}
\dRRK^n = \sum_{i=1}^{s} b_{i} \fRRK_i\, ,
\end{equation}
and the estimate of the Lyapunov function change
\begin{equation}
\label{eq:e}
e = \Delta t \sum_{i=1}^s b_i \scp{\Lyap'(\yRRK_i)}{\fRRK_i}
\end{equation}
can be computed on the fly during the calculations required by the  Runge--Kutta method and are not influenced
by the procedure used to compute $\gamma_n$ \cite{ranocha2019relaxation}.
Hence, existing low-storage Runge--Kutta implementations can be used.

In equation \eqref{eq:r}, $q(\gamma = 1)$ can be interpreted as the energy 
production of the unmodified Runge--Kutta method. Indeed,
$\Lyap\left(\URRK^{n+1}\right) - \Lyap\left(\URRK^{n}\right)$ is the Lyapunov 
change and $e$ in \eqref{eq:e} is the energy function change. 
If the weights 
of the underlying Runge--Kutta method are $b_i \geq0$, then $e$
has the same sign as the 
true time derivative of the Lyapunov function or functional, $d\Lyap/dt$.
Hence, $q$ can be seen as the temporal energy production. Thus, finding a root 
of $q$ yields a scheme that is Lyapunov consistent.
This is formally
summarized in the following theorem that was introduced in \cite{ranocha2019relaxation}  
for convex functionals in the context
of the compressible Euler and Navier--Stokes equations.
\vspace{0.5cm}
\begin{theorem}
The method defined by \eqref{eq:RK-stages} and \eqref{eq:RK-final-gamma}, 
where $\gamma_n$ is a root of \eqref{eq:r}, is Lyapunov conservative. 
If the weights $b_i$ are non-negative and $\gamma_n \geq 0$, then the method 
is Lyapunov dissipation preserving.
\end{theorem}

The following corollary provides a crucial element for preserving 
any form of Lyapunov consistency when a system of ODEs, endowed with a Lyapunov
function or functional,
is integrated with a RRK method.
\vspace{0.5cm}
\begin{Corollary}\label{cor:v_rrk}
  Any Lyapunov function of system \eqref{eq:system_odes} is also a Lyapunov function of the relaxation
  Runge--Kutta method defined by \eqref{eq:RK-stages} and \eqref{eq:RK-final-gamma} if
  the weights $b_i$ are non-negative and $\gamma_n$ is a root of \eqref{eq:r} or $\gamma_n \geq 0$.
\end{Corollary}

The following theorem provides the formal order of accuracy of relaxation Runge--Kutta 
methods.

\vspace{0.5cm}
\begin{theorem}\label{thrm:error-estimate}
  The local error for relaxation Runge--Kutta method satisfies
  \begin{subequations}
    \label{eq:local-error-augmented}
    \begin{align}
      \| \Lyap(\URRK_\gamma^{n+1})-\Lyap(\URRK(t_n+\gamma \Delta t)) \| = \mathcal{O}(\Delta t^{p+1}),\\
      \| \URRK_\gamma^{n+1} - \URRK(t_n+\gamma\Delta t) \| = \mathcal{O}(\Delta t^{p+1}),
    \end{align}
  \end{subequations}
  where $p$ is the order of the underlying Runge--Kutta method.
\end{theorem}

\begin{proof}
  The error estimates are proven in \cite{ketcheson2019relaxation,ranocha2019relaxation}.
\end{proof}

It is known that Runge--Kutta schemes can 
produce spurious solutions which, despite being non-oscillatory, are not 
solutions to the original differential equation \cite{hairer1990equilibria,griffiths1992spurious}. 
In the terminology of \citet{Iserles1990} such methods are not regular. In some 
circumstances, \citet{newell1977finite}, \citet{brezzi1984real}, and \citet{griffiths1992spurious}
found spurious solutions for certain values of the time step below 
the linearised stability limit where spurious invariant curves appear. The presence
of spurious solutions may affect computations in that,
although the fixed point of the continuous problem is globally asymptotically
stable, they generally restrict the range of initial data that are attracted to
it \cite{griffiths1992spurious}. In the following, we state the main
results that assure that 
the relaxation Runge--Kutta schemes do not produce spurious solutions.
If we define a generic time marching method as
\begin{equation*}
  \URRK^{n+1} = \PhiRRK\left(\URRK^{n}\right),
\end{equation*}
then a fixed point is a discrete
solution, $\URRK^{n}$, such that $\PhiRRK\left(\URRK^{n}\right)=\URRK^{n}$, and therefore, $\URRK^{n+1}=\URRK^{n}$. 
We define the set of fixed points of a time marching method as
\begin{equation*}
 E_{\Delta t}\equiv\left\{\URRK^{n}\in\mathbb{R}^{l}\,|\,\Phi\left(\URRK^{n}\right)=\URRK^{n}\right\},
\end{equation*}
where $l$ is the number of equations in the ODE system. Similarly, the set of equilibrium points
of the ODE is defined as
\begin{equation*}
  E\equiv\left\{\URRK\in\mathbb{R}^{l}\,|\,\bfnc{F}\left(\URRK\right)=0\right\}.
\end{equation*}
Note that for the set of fixed point of the time marching scheme we use the
subscript $\Delta t $ because it may be a function of the time step
\cite{hairer1990equilibria}.

For immediate use, we define a strict Lyapunov function for the ODE system 
as one whereby the equalities in
\begin{equation*}
  \frac{dV}{dt} = \W\Tr\bfnc{F}\leq 0
\end{equation*}
only holds on the set $E$. This implies that for $t_2>t_1$
\begin{equation}
\Lyap\left(\URRK(t_2)\right)\leq\Lyap\left(\URRK(t_{1})\right),
\end{equation}
where the equality holds only if $\URRK\in E$.

Similarly, we define a strict discrete Lyapunov function as one that satisfies
\begin{equation}\label{eq:discreteLyap}
  \Lyap\left(\URRK_{\gamma}^{n+1}\right)\leq\Lyap\left(\URRK^{n}_{\gamma}\right),
\end{equation}
where the equality holds for $\URRK^{n}\in E_{\Delta t}$.

Our first task is to determine the conditions such that $E_{\Delta t}=E$ for the
relaxation Runge--Kutta method.
The following result is due to Iserles \cite{Iserles1990} (see Theorem 3).
\vspace{0.2cm}
\begin{lemma}
  $E\subseteq E_{\Delta t}$.
\end{lemma}

In the subsequent analysis we require that $\gamma_n$ is positive; the following theorem 
(proven in~\cite{ranocha2019relaxation}) states the conditions required for $\gamma_n$ to 
be positive.
\vspace{0.5cm}
\begin{theorem}\label{thrm:positivegamma}
For all, at least second order accurate, Runge--Kutta schemes and sufficiently small $\Delta t$, 
  if $\Lyap''(\URRK^n)\scp{\fRRK_0}{\fRRK_0}>0$, where $\fRRK_0=\fRRK(\URRK^n)$, 
$q(\gamma_n)$ (\ref{eq:r}) has a positive root, \ie, $\gamma_n>0$.
\end{theorem}

\vspace{0.2cm}
\begin{remark}\label{rem:hessian}
  The statement $\Lyap''(\URRK^n)\scp{\fRRK_0}{\fRRK_0}>0$ requires that $\fRRK_0\neq 0$, which is automatically satisfied
  for models away from their equilibrium. This makes $\Lyap''(\URRK^n)>0$ the only 
  requirement we need to satisfy.
\end{remark}
\vspace{0.2cm}
\begin{remark}\label{rem:zero-rhs}
  At equilibrium points, $\fRRK_0 = 0$. Thus, we require that 
  $\Lyap''(\URRK^n)\scp{\fRRK_i}{\fRRK_i}>0$ for an intermediate step $i$, where $\fRRK_i=\fRRK(\mathbf{Y}_i)$, 
  since they are non-zero by Taylor's expansion around $t_n$ \cite{ranocha2019relaxation}. 
\end{remark}

In order to guarantee $\gamma_n>0$ (Remark \ref{rem:hessian}), 
we restrict our choice of Lyapunov functions, $\Lyap$, that satisfy $\frac{\partial^2\Lyap}{\partial\bmU^2}(\URRK^n)>0$, \ie,
convex Lyapunov functions. This choice is obviously 
model dependent. Furthermore, as highlighted at the beginning of this
section,
the convexity of $\Lyap$ is also required to ensure a one-to-one map between the
state variables, $\bmU$, and the 
Lyapunov variables, $\W$.

\vspace{0.2cm}
\begin{theorem}\label{thrm:fixedpoints}
  If the weights of the RRK method, $b_{i}>0$, and $\gamma_n> 0$ then $E_{\Delta t}=E$.
\end{theorem}
\begin{proof}
 The conditions such that $E_{\Delta t}=E$ are proven in \cite{aljahdali2024brain}.
\end{proof}

\section{Numerical Results} \label{sec:numerical-results}

In this section, we numerically verify the accuracy and Lyapunov-consistency
properties of the new discretizations by solving a dimerization problem.
Dimerization is a fundamental concept in chemistry that plays a crucial role in various chemical processes and systems. It involves the combination of two identical or similar molecules to form a dimer through the formation of chemical bonds. A  dimer is a molecular entity consisting of two monomers held together by chemical bonds. 

Consider the reversible dimerization reaction
\begin{equation}
\label{eq:reversible-dimerization}
\schemestart
2 P
\arrow{<=>[$k_f$][$k_r$]}
Q,
\schemestop
\end{equation}
where monomers $P$ dimerize at rate $k_f$, and dimers $Q$ dissociate at rate
$k_r$. The convection-diffusion-reaction dynamics of the local concentrations of
monomers and dimers at position $x \in \Omega$ and time $t > 0$ is governed by

\begin{subequations}\label{eq:pde-system-monomers-dimers}
\begin{align}
& \frac{\partial}{\partial t} \, P = d_{P} \, \nabla^{2} P - \bm{\alpha} \cdot \nabla P -2 \, k_f P^2+2 \, k_r Q, 
\\
& \frac{\partial}{\partial t} \,Q = d_{Q} \, \nabla^{2} Q  -\bm{\beta} \cdot \nabla Q+k_f P^2-k_r Q, 
\end{align}
\end{subequations}
where the positive constants $d_{P}$ and $d_{Q}$ are
the  monomers and dimers diffusion coefficients, respectively. The quantities
$\bm{\alpha}$ and $\bm{\beta}$ are the velocity field vectors. The symbol
$\Omega \subset \mathbb{R}^{dim}$ indicates a bounded domain with
piecewise smooth boundary, $\Gamma$. At the boundary, periodic boundary conditions are
used. The general initial conditions are give
by $\left(P,Q \right)(x,0) = \left(P_0,Q_0 \right)(x)$ with $x \in \Omega$.
We assume the dispersal strategy of the monomers and dimers only differs in their reaction rate, i.e., $d_{P}=d_{Q}=d$, and $\bm{\alpha}=\bm{\beta}=\bm{a}$.

Shear \cite{shear1967analog} demonstrated the existence of the Lyapunov function
for homogeneous chemical reaction, such as the reversible dimerization reaction
\eqref{eq:reversible-dimerization}. He proved that the equilibrium point,
$\bmU_{eq} = [P_{eq}, Q_{eq}]\Tr$, is globally asymptotically stable by using the following Lyapunov function
\begin{equation}\label{eq:lyapunov-function}
    \Lyap(P,Q)=P \ln(P/P_{eq})-P+P_{eq}+Q \ln(Q/Q_{eq})-Q+Q_{eq}.
\end{equation}
The Lyapunov functional for system \eqref{eq:pde-system-monomers-dimers} is
obtained by integrating in space the Lyapunov function \eqref{eq:lyapunov-function},
\begin{equation}\label{eq:lyapunov-functional}
\tilde{\Lyap}\equiv\int_{\Omega}\Lyap\mr{d}\Omega.
\end{equation}

The Lyapunov consistency presented in Section \ref{sec:lyapunov-consistency}
requires defining a two-point flux function, $\bm{\mathfrak{F}}$. 
For system \eqref{eq:reversible-dimerization}, the components of the two-point flux function 
read
\begin{equation}\label{eq:2point-flux-function}
   \mat{F}^{lc}_{x_{l}}\left(U_1,U_2\right) = a_{x_l} \, \, \bmU_{eq} \,
    \left[\left[\frac{\bmU}{\bmU_{eq}}\right]\right]_{\log}, \quad l = 1,
    \ldots, dim, 
\end{equation}
with the logarithmic average given by
    \begin{equation}
    \left[\left[\frac{\bmU}{\bmU_{eq}}\right]\right]_{\log}=\frac{\frac{\bmU_{1}}{\bmU_{eq}}-\frac{\bmU_{2}}{\bmU_{eq}}}{\log{\frac{\bmU_{1}}{\bmU_{eq}}}-\log{\frac{\bmU_{2}}{\bmU_{eq}}}},
\end{equation}
where $\bm U_{1}=[P_1,Q_1]$, $\bm U_{2}=[P_2,Q_2]$, and the subscripts $1$
and $2$ indicate the two states needed for the two-point flux function.

\subsection{Convergence Study} 
In this section, we aim to assess the order of accuracy of the proposed fully-discrete schemes by applying the method of manufactured solution (MMS) \cite{roy2005}. 
Here, we use the MMS to study the order of convergence of our new schemes
applied to the convection-diffusion-reaction system \eqref{eq:pde-system-monomers-dimers} that models reversible dimerization chemical
reactions. 

The computational domain is a cube $\Omega = \left[0,2\pi\right]^{3}$, where the following manufactured solution is used for the convergence study:
 
\begin{equation}\label{pde-p-q-manu}
   P = Q = \prod_{i=1,2,3}f_{i}(x_{i},t) \, \, \text{where}\, \, f_{i}(x_{i},t) = 1.25 + 0.75 \sqrt{\cos{\left(-\pi t + \pi x_{i}\right)}}.
\end{equation}
The final time of the simulation is set to $T_{f} = 1$. We run the convergence study with initial local concentrations of monomers and dimers as $P_0=10$ for monomers and $Q_0=1$ for dimers, respectively. The dimerization rate of monomers is set to $k_f=10$, and the dissociation rate of dimers is set to $k_r=1$. The diffusion coefficient $d$ is fixed at 0.05, and the velocity field vector $\bm{a}$ is set to the unit vector $\bm{1}$.

The convergence rates for the discrete error norms $L^{1}$, $L^{2}$, and
$L^{\infty}$, defined in  Appendix \ref{sec:discrete-error-norms}, corresponding
to second- ($p=1$),  third- ($p=2$), fourth- ($p=3$), and fifth-order ($p=4$) accurate fully-discrete algorithms, are shown in Table \ref{tab:conv_relax_time}. 
The solution polynomial degree $p$ of the DG SBP-SAT  operator used in the spatial discretization and the number of cells $K$ in the computational mesh are shown in the table.
In addition, we listed the relaxation Runge--Kutta time integration methods that are used
to evolve the system of ODEs arising from the DG SBP-SAT spatial discretization \cite{ALJAHDALI2022111333}. 
We observe that numerical schemes with even-degree DG solution polynomials exhibit an $L^2$ convergence rate of approximately $p+1$. In contrast, those with odd-degree DG polynomials converge at a higher rate. A similar pattern is seen for the $L^1$ and $L^{\infty}$ error norms, except the $p=1$ scheme in the case of the $L^{\infty}$ norm.

\begin{table}[H]
\centering
\setlength\tabcolsep{0.5pt}
  \caption{Convergence study for the convection-diffusion-reaction dynamics of
    the dimerization chemical reaction model using Lyapunov consistent DG
    SBP-SAT schemes with different solution polynomial degrees, $p$, and
    relaxation Runge--Kutta methods.}
  \label{tab:conv_relax_time}
\begin{tabular*}{\linewidth}{@{\extracolsep{\fill}}*9c@{}}
    \toprule
    $p$ & RK Method
    &$K$& $L^1$ Error & $L^1$ Rate
    & $L^2$ Error & $L^2$ Rate
    & $L^\infty$ Error & $L^\infty$ Rate
    \\
    \midrule
      1 & ERK(8,2)  &
512      & 7.99854e-01 &  -    & 9.56101e-01 &   -   & 3.18092e+00 & - \\&& 
4,096    & 3.56431e-01 & 1.166 & 4.23122e-01 & 1.176 & 1.27381e+00 & 1.320  \\&& 
32,768   & 1.09019e-01 & 1.709 & 1.27223e-01 & 1.734 & 3.72503e-01 & 1.774  \\&& 
262,144  & 2.78390e-02 & 1.969 & 3.22426e-02 & 1.980 & 9.07434e-02 & 2.037 \\&& 
2,097,152& 6.87503e-03 & 2.018 & 8.05057e-03 & 2.002 & 2.20967e-02 & 2.038  \\
    \midrule
      2 & BSRK(4,3)-3(2) pair  &

  1,728         & 1.14981e-01   &    -   & 1.32592e-01 &   -    & 4.56714e-01 & - \\&& 
  13,824        & 8.51894e-03   & 3.755  & 8.85260e-03 &  3.905 & 4.23364e-02 & 3.431   \\&& 
  110,592       & 9.15459e-04   & 3.218  & 8.82595e-04 &  3.326 & 4.77381e-03 & 3.149   \\&& 
  884,736       & 1.09862e-04   & 3.059  & 1.07194e-04 &  3.042 & 5.92147e-04 & 3.011   \\&& 
  7,077,888     & 1.35060e-05   & 3.024  & 1.34799e-05 &  2.991 & 7.30730e-05 & 3.019   \\
    \midrule
      3 & RK(4,4)  &
 4,096  & 4.95896e-03 & -     & 5.64148e-03 & -     & 2.87246e-02 &  - \\&&  
 32,768  & 5.03612e-04 & 3.299 & 5.33964e-04 & 3.401 & 3.10829e-03 & 3.208   \\&& 
 262,144  & 2.93424e-05 & 4.101 & 2.83613e-05 & 4.235 & 1.81243e-04 & 4.100   \\&& 
 2,097,152  & 1.28045e-06 & 4.518 & 1.16181e-06 & 4.609 & 7.75138e-06 & 4.547   \\&& 
 16,777,216  & 4.61820e-08 & 4.793 & 4.12820e-08 & 4.815 & 2.72962e-07 & 4.828   \\

            \midrule
      4 & BSRK(8,5)-4(5) pair  &
8,000     &7.35452e-04  &   -    & 8.41095e-04 &  -    & 6.28798e-03 &  -\\&& 
64,000     &2.34012e-05  & 4.974 & 2.62451e-05 &5.002 & 2.10179e-04 & 4.903      \\&& 
512,000     &7.20832e-07  & 5.021 & 7.63498e-07 &5.103 & 6.37179e-06 & 5.044      \\&&  
4,096,000     &2.19966e-08  & 5.034 & 2.34475e-08 &5.025 & 2.07971e-07 & 4.937    \\&&
32,768,000     &6.80908e-10  & 5.014 & 7.33225e-10 &4.999 & 6.70438e-09 & 4.955     \\ 
    \bottomrule
  \end{tabular*}
\end{table}

Finally, we study the time-dependent coefficients for the diffusion terms in the
same PDE model \eqref{eq:pde-system-monomers-dimers}. For this purpose, the
monomers and dimers diffusion coefficients  $d_{P}$ and $d_{Q}$ are multiplied by $C_{fact} = 1 + 0.5 \sin(2 \pi t)$, where $t$ is the time.
The results are shown in Table \ref{tab:conv_iv_p3to5_relax_time}.
We observe similar orders of accuracy to those reported for the previous case, where the diffusion coefficients remained constant over time. 

\begin{table}[H]
\centering
\setlength\tabcolsep{0.5pt}
  \caption{Convergence study for the convection-diffusion-reaction dynamics of
    the dimerization chemical reaction model using Lyapunov consistent DG
    SBP-SAT schemes with different solution polynomial degrees, $p$, and relaxation
    Runge--Kutta methods. The diffusion coefficients $C_{l,m}$ are time-dependent.}
  \label{tab:conv_iv_p3to5_relax_time}
  \begin{tabular*}{\linewidth}{@{\extracolsep{\fill}}*9c@{}}
    \toprule
    $p$ & RK Method
    &$K$& $L^1$ Error & $L^1$ Rate
    & $L^2$ Error & $L^2$ Rate
    & $L^\infty$ Error & $L^\infty$ Rate
    \\
    \midrule
      1 & ERK(8,2)  &
 512             & 7.95096e-01 &  -    & 9.51325e-01&   -   & 3.20918e+00 & - \\&& 
 4,096           & 3.43179e-01 & 1.212 & 4.06416e-01& 1.227 & 1.20634e+00 & 1.412  \\&& 
 32,768          & 1.05092e-01 & 1.707 & 1.22263e-01& 1.733 & 3.58101e-01 & 1.752  \\&& 
 262,144         & 2.69018e-02 & 1.965 & 3.10639e-02& 1.977 & 8.77128e-02 & 2.030 \\&& 
 2,097,152       & 6.64138e-03 & 2.018 & 7.75579e-03& 2.009 & 2.13368e-02 & 2.039  \\
    \midrule
      2 & BSRK(4,3)-3(2) pair  &
  1,728     & 1.11522e-01 &    -   & 1.28893e-01 &   -    & 4.55491e-01 & - \\&& 
  13,824    & 8.43921e-03 & 3.724  & 8.75667e-03 &  3.880 & 4.17861e-02 & 3.446   \\&& 
  110,592   & 9.13462e-04 & 3.208  & 8.80679e-04 &  3.314 & 4.75647e-03 & 3.135   \\&& 
  884,736   & 1.09778e-04 & 3.057  & 1.07145e-04 &  3.039 & 5.91121e-04 & 3.008   \\&& 
  7,077,888 & 1.35018e-05 & 3.023  & 1.34784e-05 &  2.991 & 7.30127e-05 & 3.017   \\
 
    \midrule
      3 & RK(4,4)  &
  4,096      & 4.85403e-03 &   -   & 5.51374e-03 &    -  & 2.79356e-02 &  - \\&&  
  32,768     & 4.95293e-04 & 3.293 & 5.24052e-04 & 3.395 & 3.07503e-03 & 3.183   \\&& 
  262,144    & 2.91624e-05 & 4.086 & 2.81495e-05 & 4.219 & 1.80249e-04 & 4.093   \\&& 
  2,097,152  & 1.27850e-06 & 4.512 & 1.15982e-06 & 4.601 & 7.73758e-06 & 4.542   \\&& 
  16,777,216 & 4.61633e-08 & 4.792 & 4.12647e-08 & 4.813 & 2.72839e-07 & 4.826   \\

            \midrule
      4 & BSRK(8,5)-4(5) pair  &
8,000     & 7.33054e-04 &    -  & 8.39935e-04 &    -  & 6.29025e-03 & - \\&& 
64,000    & 2.34562e-05 & 4.966 & 2.62946e-05 & 4.997 & 2.10260e-04 & 4.903    \\&& 
512,000   & 7.21108e-07 & 5.024 & 7.63691e-07 & 5.106 & 6.37183e-06 & 5.044    \\&&  
4,096,000 & 2.19973e-08 & 5.035 & 2.34481e-08 & 5.025 & 2.07940e-07 & 4.937    \\&& 
32,768,000& 6.80907e-10 & 5.014 & 7.33227e-10 & 4.999 & 6.76091e-09 & 4.943    \\ 
    \bottomrule
  \end{tabular*}
\end{table}

\subsection{Chemical reaction for the reversible dimerization model}

To verify numerically the Lyapunov consistency properties of the algorithms, we
compute the spatio-temporal evolution of the solution of the system
\eqref{eq:pde-system-monomers-dimers} and monitor the convergence to the
equilibrium point. 
The simulations are carried out using a spatially fourth-order accurate
Lyapunov consistent DG SBP-SAT operator (i.e., the solution polynomial degree is $p=3$). The time integration is performed using the RK 3(2) method of Bogacki--Shampine with adaptive time stepping \cite{ALJAHDALI2022111333}.

To obtain the equilibrium condition of system
\eqref{eq:pde-system-monomers-dimers}, we use the rate equations that describe
the local concentrations of monomers and dimers, i.e.,
\begin{subequations}\label{eq:rate-equations}
\begin{align}
& P_t = -2 k_f P^2+2 \, k_r Q, \\
& Q_t = k_f P^2-k_r Q.
\end{align}
\end{subequations}
Then, the equilibrium condition reads 
\begin{equation}\label{eq:equilibrium-condition}
\frac{Q_{eq}}{P_{eq}^2} = \frac{k_f}{k_r}, 
\end{equation}
and the conservation relation $P+2Q=P_0+2Q_0$ can be used to find the equilibrium
point: 
 \begin{equation}\label{eq:rate-equations}
    \bmU_{eq} = \left[P_{eq}, Q_{eq} \right]= \left[\frac{-k_r+\sqrt{k_r(8
     k_f(P_0+2Q_0)+k_r)}}{4k_f}, \frac{P_0+2Q_0-P_{eq}}{2} \right]\Tr.
\end{equation}

Here, we show the time evolution of the model
\eqref{eq:pde-system-monomers-dimers} towards the equilibrium \eqref{eq:rate-equations}. 
The simulation is performed in a cube domain with a side length of  $1$ unit.
This domain is discretized into $64$ hexahedral cells in each coordinate
direction. To investigate the impact of the reaction rate on the local
concentrations of monomers and dimers, we distribute the monomers and dimers using a spherical harmonics shape
located at the center of the cube, which is expressed in spherical coordinates
as $HS =\tau \, (1+0.2 \, (\sin(7\theta)\sin(4\phi))$, with a radius of $\tau=0.30$.
Within this region, the local concentrations of monomers and dimers are assigned
specific values, $P_0=10$ for monomers and $Q_0=1$ for dimers. The local concentrations of monomers and dimers are assigned
different values in the surrounding area, $P_0=0.1$ and $Q_0=0.1$, respectively. The dimerization rate
of monomers is set to $k_f=10$, and the dissociation rate of dimers is set to
$k_r=1$. The diffusion coefficient value is chosen as $d=0.05$, and the
velocity field vector is $\bm{a}=\bm{1}$.

\begin{figure}[H]
\centering
  \begin{minipage}[b]{.9\textwidth}
      {\includegraphics[width=0.98\textwidth]{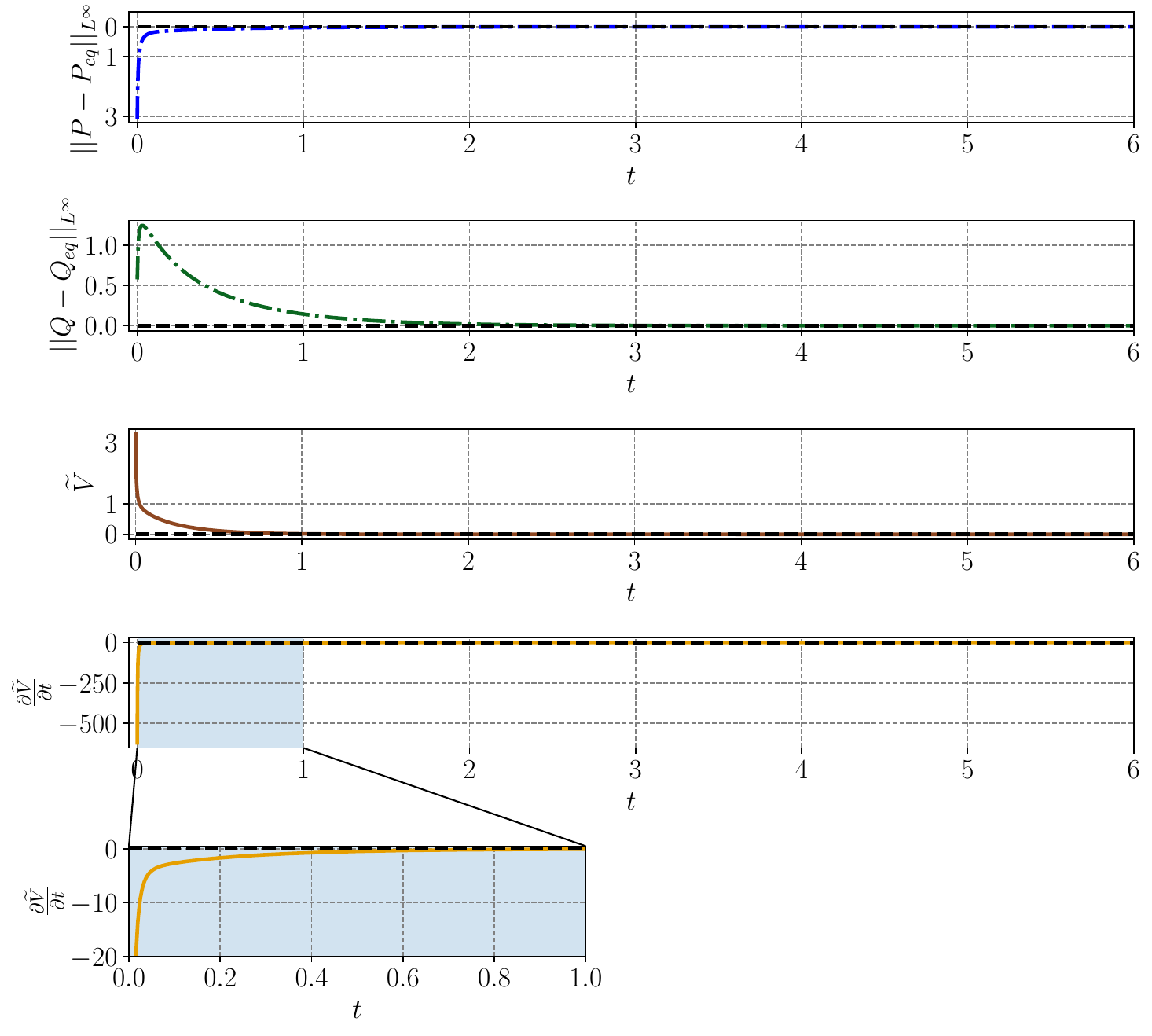}}
  \end{minipage}
\vspace*{3mm}
    \caption{Temporal evolution of the maximum norm of the difference of the solution
    $(P,Q)$ and the equilibrium point $(P_{eq}, Q_{eq})$, Lyapunov functional, $\widetilde{\Lyap}$, and time derivative of the Lyapunov functional, $\frac{d \widetilde{\Lyap}}{d t}$, for the chemical reaction model reversible dimerization.} 
\label{fig:p_q_model}
\end{figure}

The first two one-dimensional plots at the top of Figure  \ref{fig:p_q_model} illustrate the
convergence of the solution to the equilibrium point $\bmU_{eq}= (P_{eq}, Q_{eq})=(0.267912, 0.717767)$. 
Specifically, these plots show the temporal changes in the maximum norm of the difference 
between the solution $(P,Q)$, with respect to the equilibrium point $(P_{eq}, Q_{eq})$. 
We observe that the solution gradually converges toward the equilibrium value over time.
In addition, Figure  \ref{fig:p_q_model} shows the time evolution of the Lyapunov functional, 
$\widetilde{\Lyap}$, and its derivative, $d\widetilde{\Lyap}/dt$, in the third and fourth one-dimensional
plots, respectively. The spatiotemporal discretization preserves the convexity of  $\widetilde{\Lyap}$ and 
the negative definiteness of $d\widetilde{\Lyap}/dt$. These properties are crucial for establishing the global asymptotic stability of the equilibrium point.

 \begin{figure}[H]
    \centering
   \begin{minipage}{.43\textwidth}
            \begin{subfigure}{\textwidth}
            \centering
            \includegraphics[width=\textwidth,trim={10cm 1cm 12cm 0},clip]{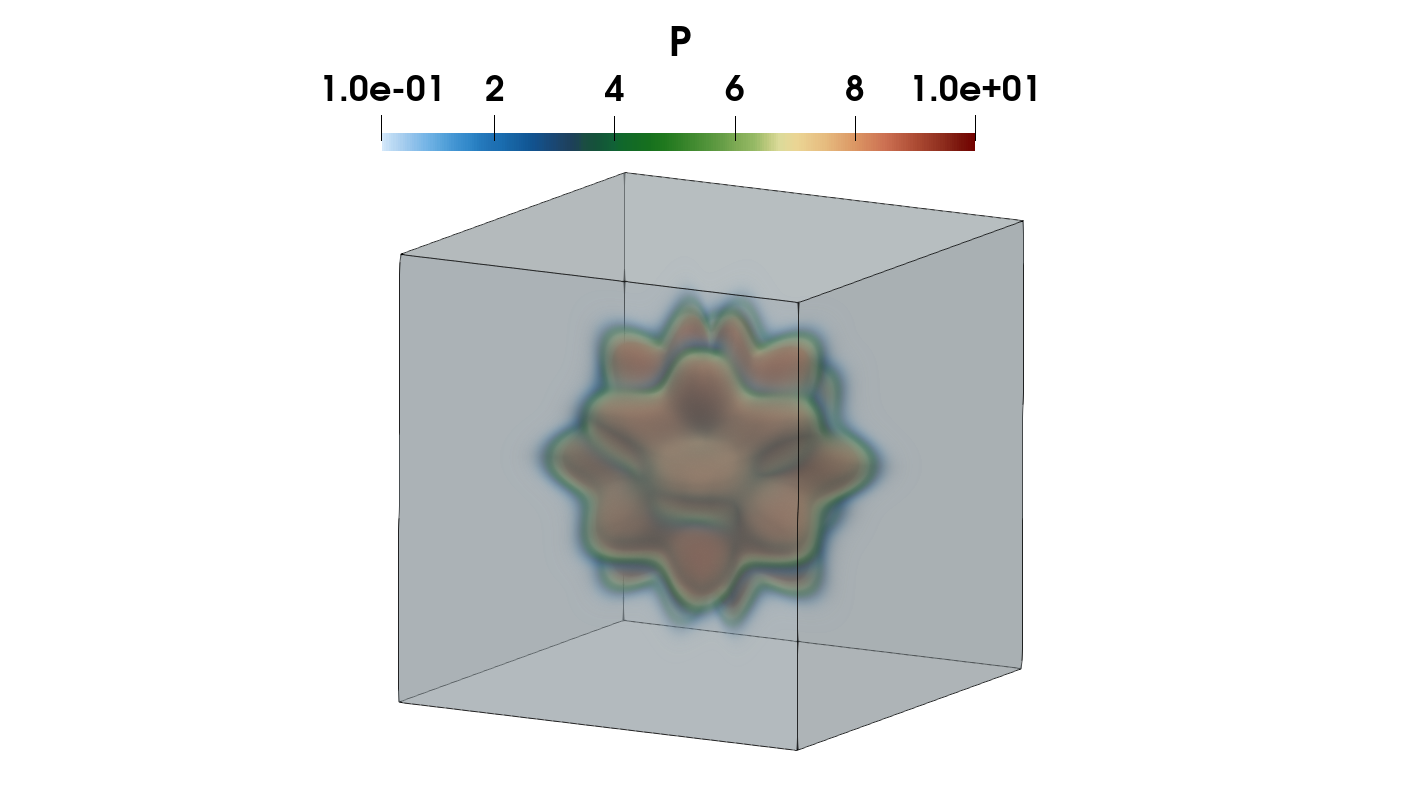}
            \end{subfigure}
            \subcaption{$t=0$}
        \end{minipage}%
        \hspace*{2.2cm}
   \begin{minipage}{.3\textwidth}
            \begin{subfigure}{\textwidth}
            \centering
            \includegraphics[width=\textwidth]{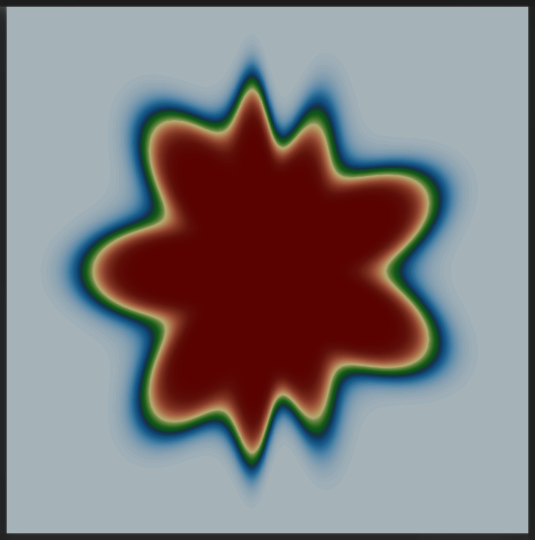}
            \end{subfigure}
            \subcaption*{\centerline{ Vertical cross-section of the domain at $t=0$}}
        \end{minipage}
   \begin{minipage}{.43\textwidth}
            \begin{subfigure}{\textwidth}
            \centering
            \includegraphics[width=\textwidth,trim={10cm 1cm 12cm 0},clip]{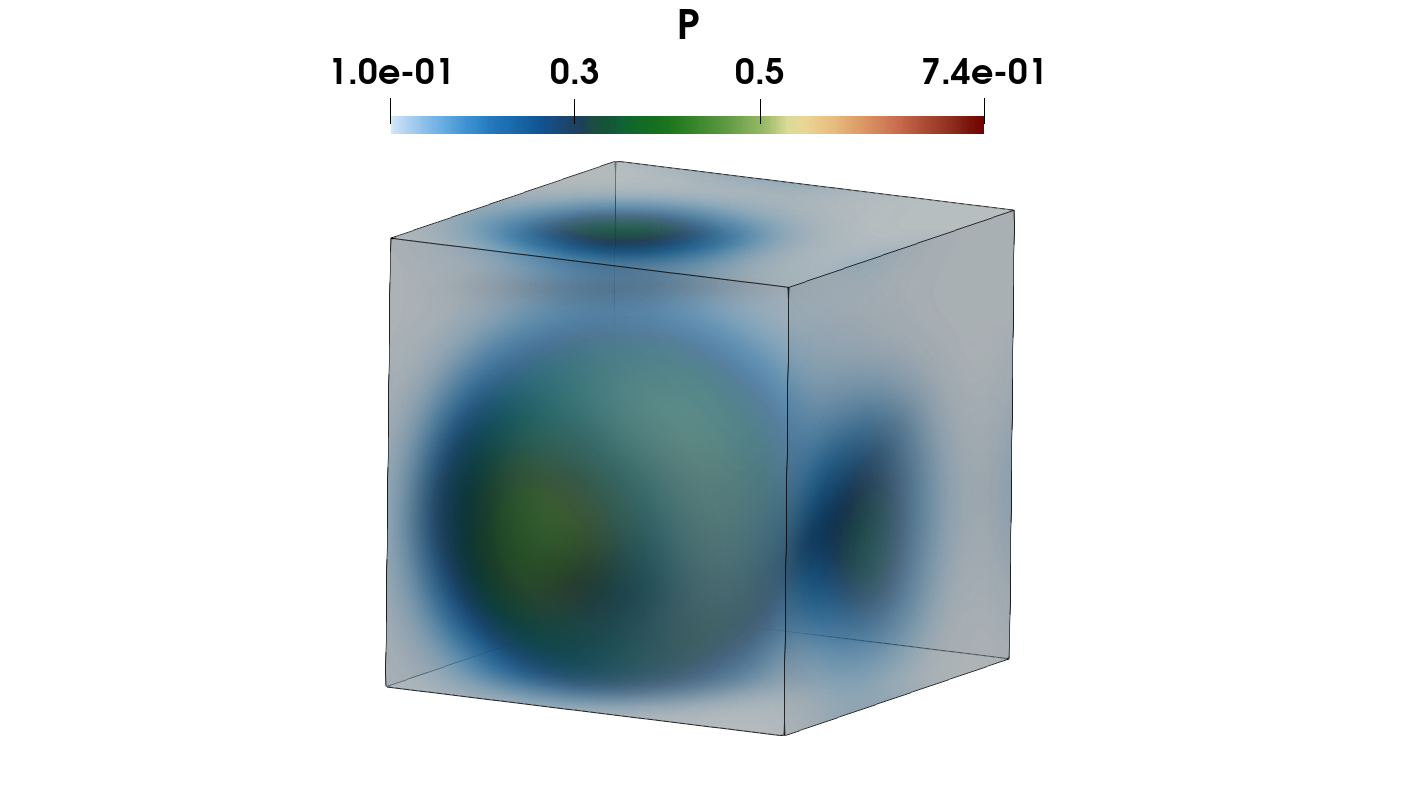}
           \end{subfigure}\
            \subcaption{$t=0.127$}
        \end{minipage}
                \hspace*{2.2cm}
        \begin{minipage}{.3\textwidth}
            \begin{subfigure}{\textwidth}
            \centering
            \includegraphics[width=\textwidth]{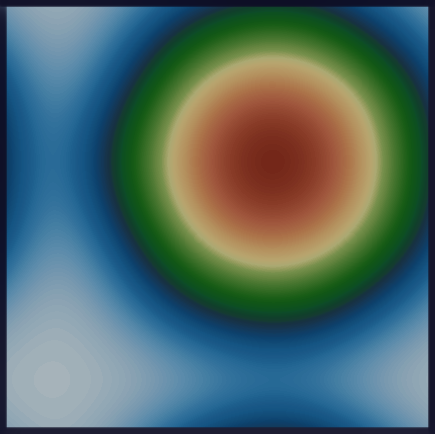}
            \end{subfigure}\
            \subcaption*{ \centerline{Vertical cross-section of the domain at $t=0.127$.}}
        \end{minipage}
           \begin{minipage}{.43\textwidth}
            \begin{subfigure}{\textwidth}
            \centering
            \includegraphics[width=\textwidth,trim={10cm 1cm 11cm 0},clip]{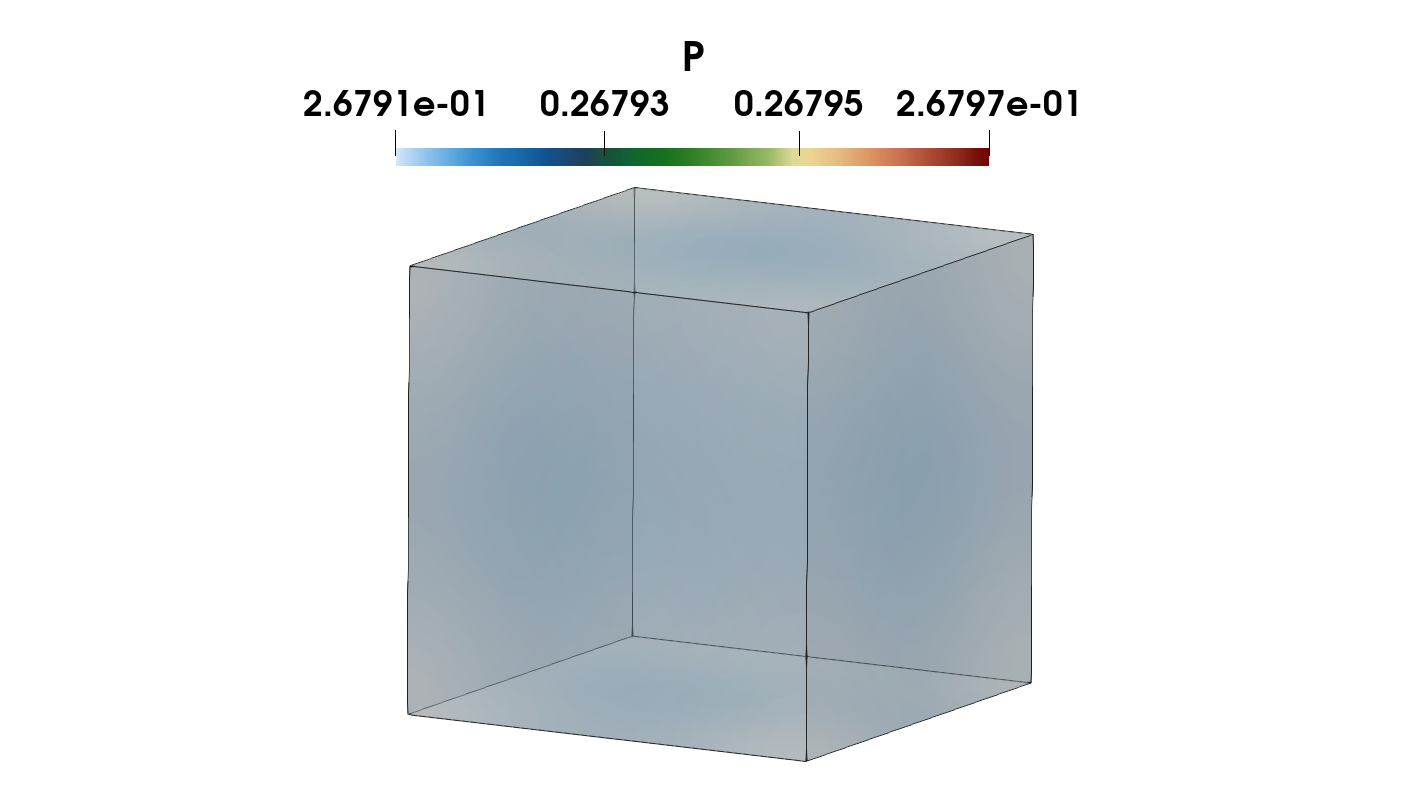}
           \end{subfigure}
            \subcaption{$t=6$}
        \end{minipage}%
                \hspace*{2.2cm}
        \begin{minipage}{.3\textwidth}
            \begin{subfigure}{\textwidth}
            \centering
            \includegraphics[width=\textwidth]{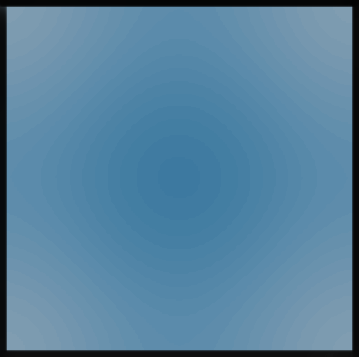}
            \end{subfigure}\\ 
            \subcaption*{ \centerline{Vertical cross-section of the domain at $t=6$.}}
        \end{minipage}%
           \caption{The local concentrations of monomers $P$ at the initial transition at $t=0$, the dampening stage at $t=0.127$, and the final transition to a uniform state at $t=6$. \label{fig:p}}
    \end{figure}
    
\begin{figure}[H]
    \centering
   \begin{minipage}{.43\textwidth}
            \begin{subfigure}{\textwidth}
            \centering
            \includegraphics[width=\textwidth,trim={10cm 1cm 12cm 0},clip]{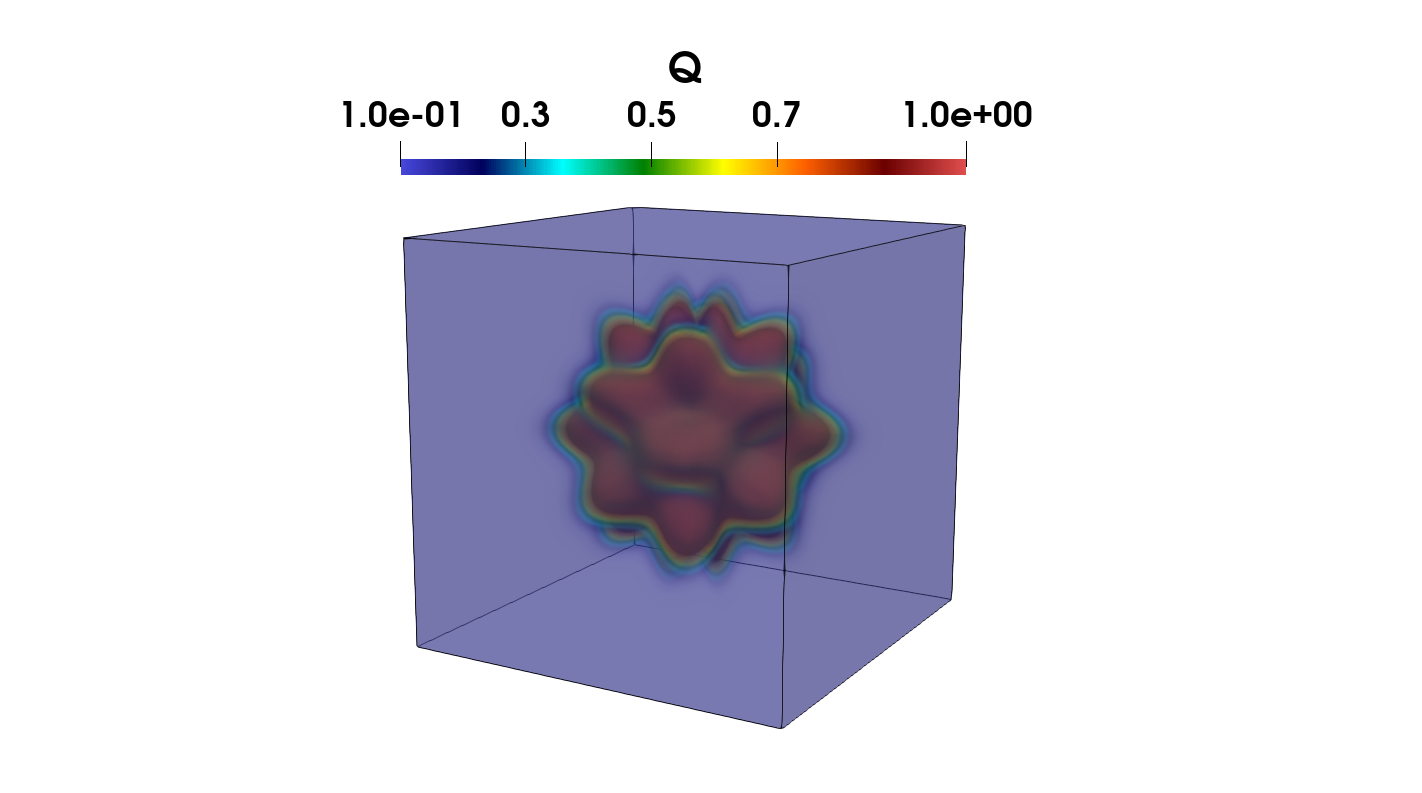}
            \end{subfigure}
            \subcaption{$t=0$}
        \end{minipage}%
        \hspace*{2.2cm}
   \begin{minipage}{.3\textwidth}
            \begin{subfigure}{\textwidth}
            \centering
            \includegraphics[width=\textwidth]{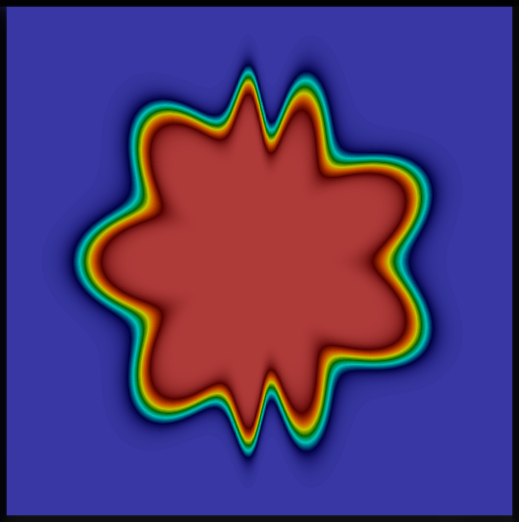}
            \end{subfigure}
            \subcaption*{\centerline{ Vertical cross-section of the domain at $t=0$}}
        \end{minipage}
   \begin{minipage}{.43\textwidth}
            \begin{subfigure}{\textwidth}
            \centering
            \includegraphics[width=\textwidth,trim={10cm 1cm 12cm 0},clip]{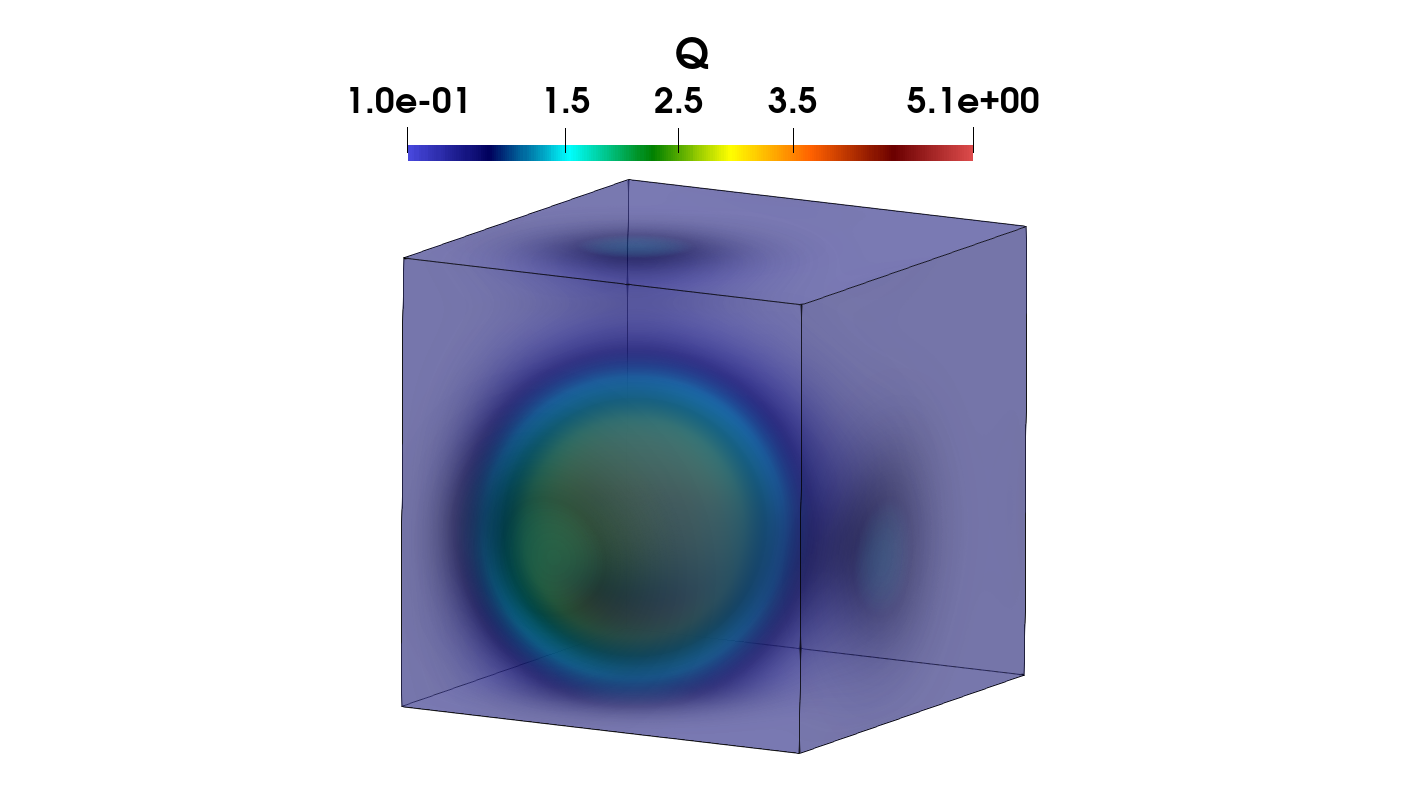}
           \end{subfigure}\
            \subcaption{$t=0.127$}
        \end{minipage}
                \hspace*{2.2cm}
        \begin{minipage}{.3\textwidth}
            \begin{subfigure}{\textwidth}
            \centering
            \includegraphics[width=\textwidth]{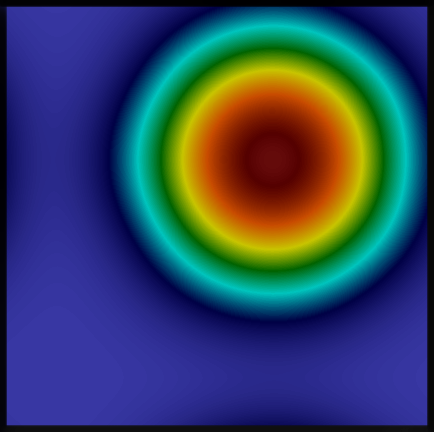}
            \end{subfigure}\
            \subcaption*{ \centerline{Vertical cross-section of the domain at $t=0.127$.}}
        \end{minipage}
           \begin{minipage}{.43\textwidth}
            \begin{subfigure}{\textwidth}
            \centering
            \includegraphics[width=\textwidth,trim={10cm 1cm 11cm 0},clip]{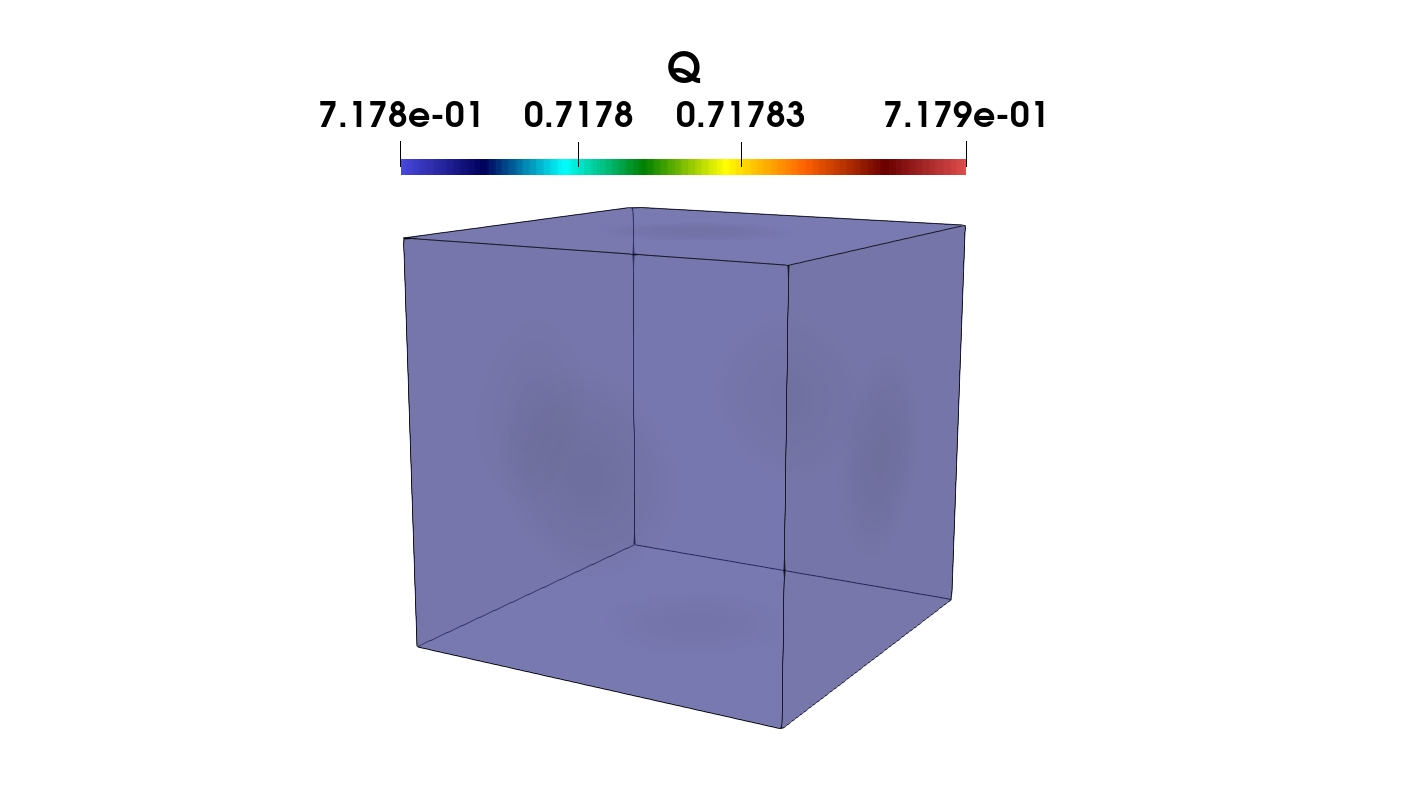}
           \end{subfigure}
            \subcaption{$t=6$}
        \end{minipage}%
                \hspace*{2.2cm}
        \begin{minipage}{.3\textwidth}
            \begin{subfigure}{\textwidth}
            \centering
            \includegraphics[width=\textwidth]{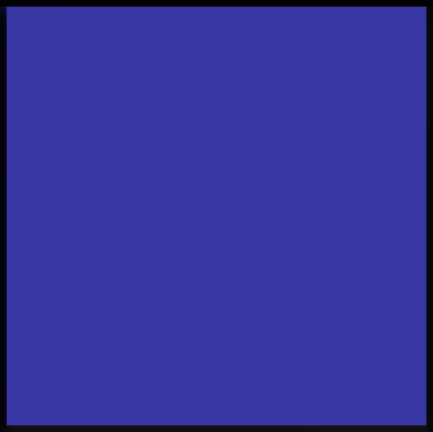}
            \end{subfigure}\\ 
            \subcaption*{ \centerline{Vertical cross-section of the domain at $t=6$.}}
        \end{minipage}%
           \caption{The local concentrations of monomers $Q$ at the initial transition at $t=0$, the dampening stage at $t=0.127$, and the final transition to a uniform state at $t=6$. \label{fig:q}}
    \end{figure}

Figures \ref{fig:p} and \ref{fig:q} show the three-dimensional contour plots and vertical cross-sections illustrating the local concentrations of monomers $P$ and dimers $Q$. The plots correspond to the initial time, $t=0$, the dampening stage at $t=0.127$, and the final uniform state at $t=6$. These results demonstrate that the solution exhibits a persistent and consistent pattern, undergoing gradual reaction-diffusion-convection processes until it converges to a state of equilibrium characterized by uniform values in the whole domain, $\Omega$. In addition, in Appendix \ref{sec:mesh-Convergence}, we show the influence of mesh refinement on the convergence rate to the equilibrium point by conducting numerical simulations of the dimerization model with various numbers of cells.

\subsection{Preserving Lyapunov stability with $hp$-nonconforming discretization} 

\begin{figure}[H]
\centering
  \begin{minipage}[b]{.93\textwidth}
      {\includegraphics[width=0.98\textwidth]{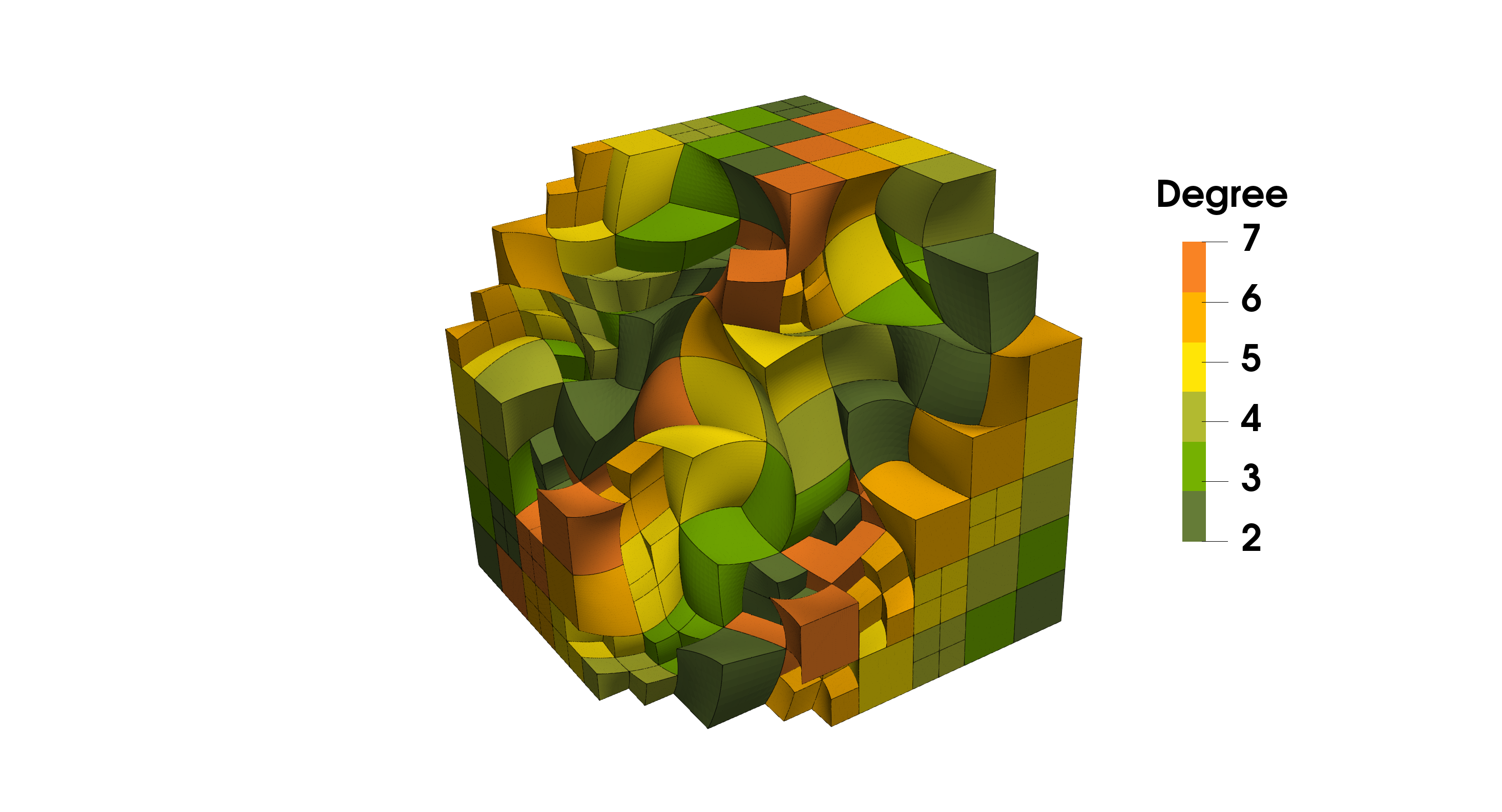}}
  \end{minipage}
\vspace*{0.1mm}
    \caption{Mesh structure with $hp$-adaptation. The colors indicating the solution
        polynomial degree distribution ($p=2$ to $p=7$) with one level of random mesh refinement.} 
\label{fig:hp-adaptation}
\end{figure}

To demonstrate that the algorithms maintain the Lyapunov stability of the continuous
model when using $hp$-adaptation \cite{fernandez_entropy_stable_hp_ref_snpdea_2019}, we conducted a numerical simulation of the dimerization
model \eqref{eq:pde-system-monomers-dimers}. A cubical computational domain with periodic boundary conditions on all six faces is used.
We simulated this model using
a grid with seven hexahedral elements in each coordinated direction, featuring a nonconforming 
interface, as shown in Figure \ref{fig:hp-adaptation}.
The grid is generated by assigning
a random integer solution polynomial degree to each element, selected from the set 
$\{2, 3, 4, 5\}$. A second-order degree polynomial representation is used for all the boundary cell interfaces 
except those with periodic boundary conditions.
All dissipation terms
related to interface coupling \cite{Parsani2015,fernandez_entropy_stable_hp_ref_snpdea_2019},
including upwind and interior-penalty SATs,  are turned off. The time integration is carried 
out using the BSRK(4,3)-3(2) pair method of Bogacki--Shampine with relaxation and adaptive time stepping based on error estimation \cite{ranocha2019relaxation,ALJAHDALI2022111333}. We maintain the
same parameter values as in the previous section.

\begin{figure}[H]
\centering
  \begin{minipage}[b]{.9\textwidth}
      {\includegraphics[width=0.98\textwidth]{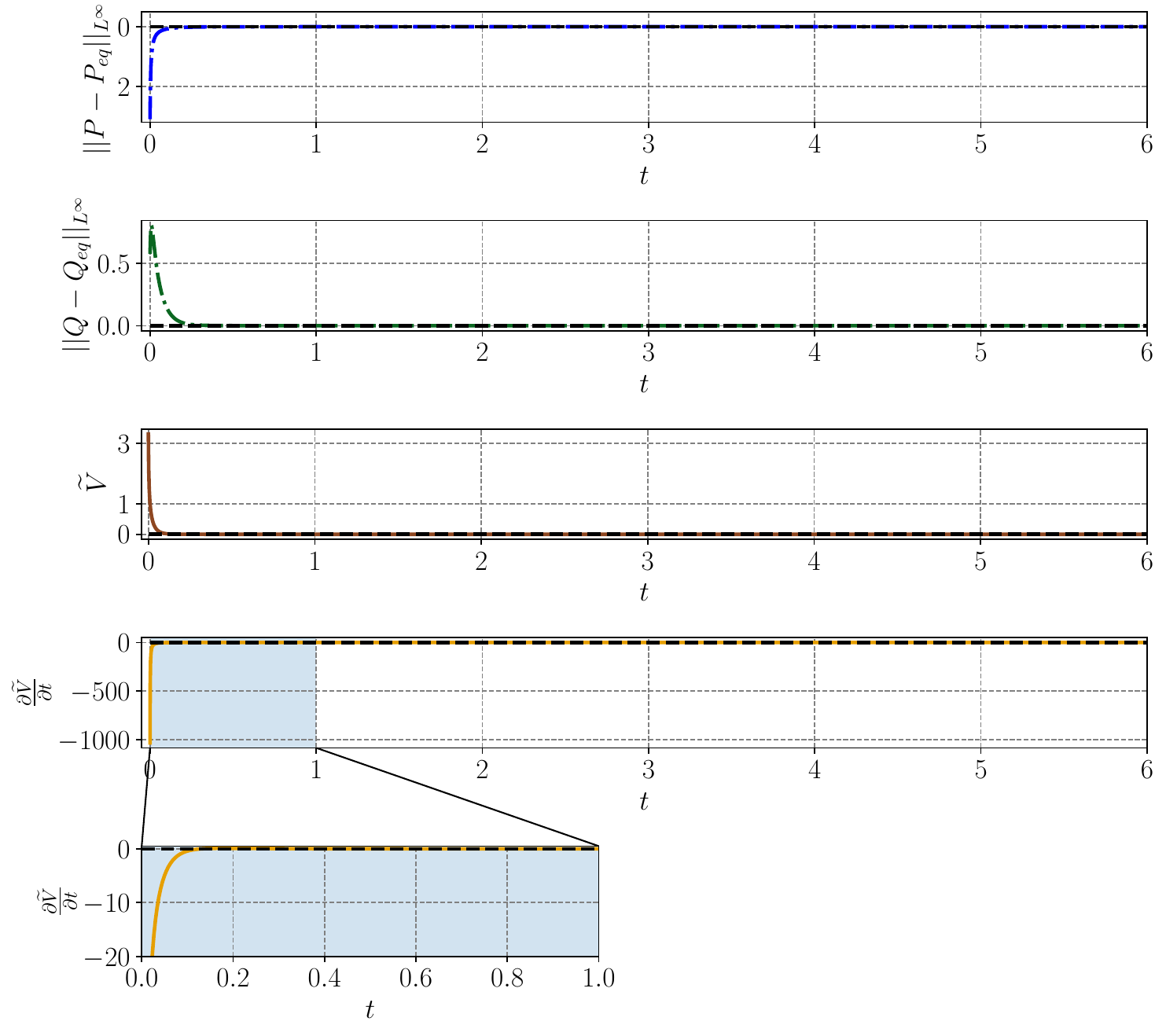}}
  \end{minipage}
\vspace*{0.1mm}
    \caption{Temporal evolution of the maximum norm of the difference of the solution
    $(P,Q)$ and the equilibrium point $(P_{eq}, Q_{eq})$, Lyapunov functional, $\widetilde{\Lyap}$, and time derivative of the Lyapunov functional, $\frac{d \widetilde{\Lyap}}{d t}$, for the chemical reaction model reversible dimerization with $hp$ -nonconforming discretization.} 
\label{fig:p_q_model_hp}
\end{figure}

We show the time evolution of the two-component model 
\eqref{eq:pde-system-monomers-dimers} and its Lyapunov functional in Figure \ref{fig:p_q_model_hp}. The top 
two one-dimensional plots illustrate the solution's convergence to the equilibrium point  $\bmU_{eq}= (P_{eq}, Q_{eq}) = (0.267912, 0.717767)$. These
plots show the temporal changes in the maximum norm of the difference between
the solution $(P,Q)$, with respect to the equilibrium point $(P_{eq}, Q_{eq})$.
We observe that the solution gradually converges toward the equilibrium point. Furthermore, the third and fourth plots show 
the time evolution of the Lyapunov functional, $\widetilde{\Lyap}$, and its derivative,
$d\widetilde{\Lyap}/dt$, respectively. Also, in this case, we observe that the spatiotemporal discretization preserves the convexity of  $\widetilde{\Lyap}$ and 
the negative definiteness of $d\widetilde{\Lyap}/dt$.

\begin{figure}[H]
\centering
  \begin{minipage}[b]{.93\textwidth}
      {\includegraphics[width=0.98\textwidth]{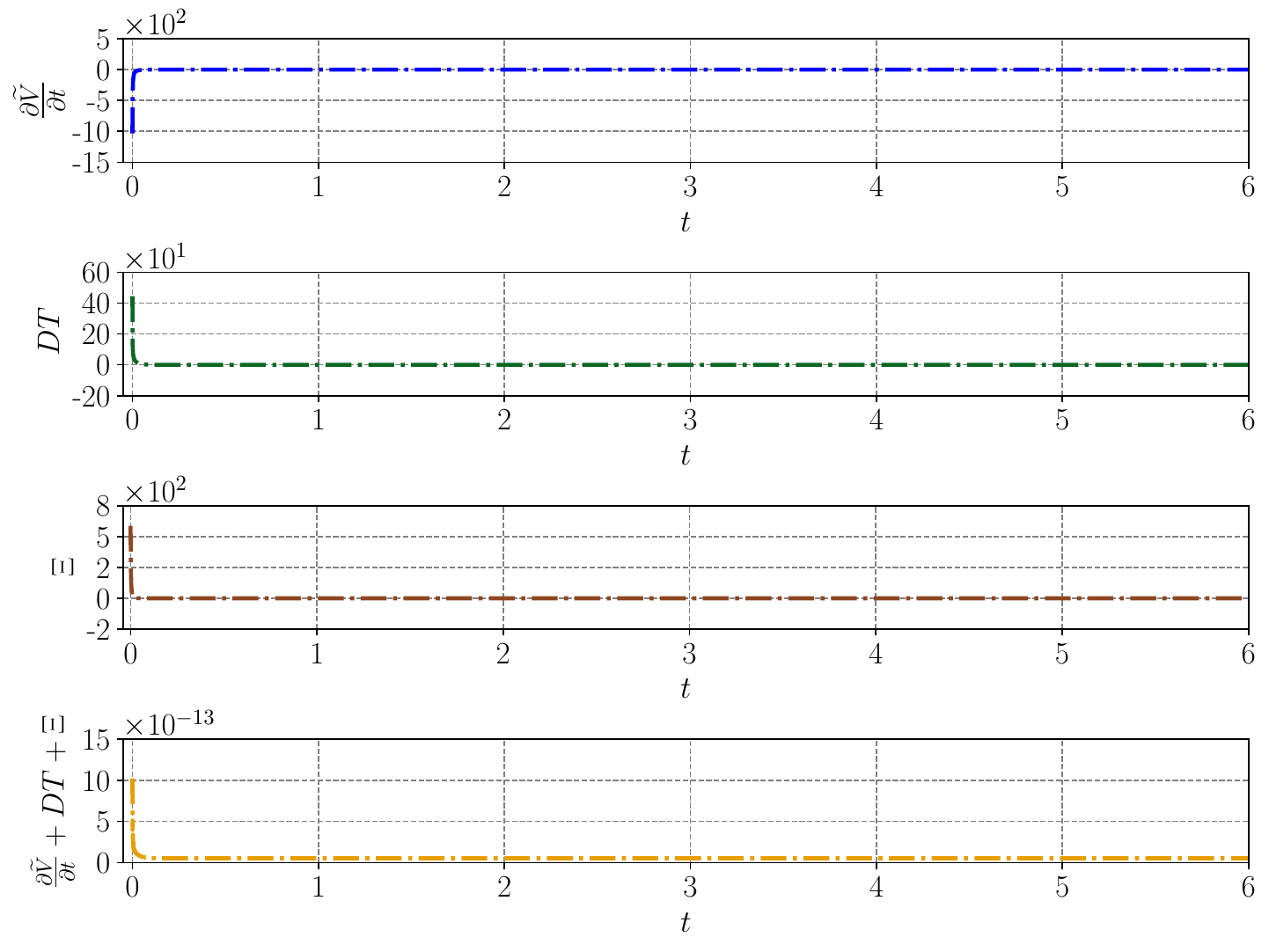}}
  \end{minipage}
\vspace*{0.1mm}
    \caption{Temporal evolution of the discrete-time rate of change of the Lyapunov functional, $d\widetilde{\Lyap}/dt$, the discrete dissipation term, $DT$, the discrete contribution in the time derivative of Lyapunov functional of reaction term, $\Xi$, and the instantaneous Lyapunov functional balance, $d\widetilde{\Lyap}/dt+DT+\Xi$. } 
\label{fig:hp-balance}
\end{figure}

To illustrate that the space and time discretizations and their coupling are 
indeed Lyapunov conservative, we compute the terms of the discrete Lyapunov functional 
balance \eqref{eq:discentropy2}. In Figure \ref{fig:hp-balance}, we show the discrete-time rate of 
change of the Lyapunov functional, $d\widetilde{\Lyap}/dt$, the discrete dissipation 
term, $DT$, and the discrete contribution of the reaction term, $\Xi = \left(\wk{g}\right)\Tr \,
      \MJ^{g} \, \fk{g}$. Additionally, we plot the instantaneous Lyapunov functional
balance, $d\widetilde{\Lyap}/dt+DT+\Xi$. The latter quantity is below machine's precision, numerically verifying the Lyapunov conservative property of the discretization.
\section{Discussion} \label{sec:discussion-conclusions}

This work presents a general framework for constructing fully discrete Lyapunov-consistent discretizations for convection-diffusion-reaction PDE models. 
The proposed methods mimic the Lyapunov stability analysis at the fully discrete level. 
The spatial discretization is based on collocated discontinuous Galerkin operators
with the summation-by-part property and the simultaneous approximation term technique for unstructured grids with 
tensor-product elements. In contrast, the temporal integration relies on relaxation Runge-Kutta methods. 
The proposed methodology is general and yields arbitrarily high-order accurate discretizations in space and time.
To illustrate the capabilities of the proposed fully discrete algorithms, we solve a convection-diffusion-reaction system 
that models the evolving dynamics between monomer and dimer concentrations during the dimerization process. Finally, 
it is essential to highlight that the algorithms apply to a wide range of convection-reaction-diffusion
equations, offering robust numerical tools that can provide valuable insights into the dynamics of complex systems.

\begin{appendices}

\section{Global SBP operators for the first and second derivatives}\label{app:globalSBP}
In this section, we will first introduce SBP operators 
constructed in computational space and their properties. These are then extended to 
Cartesian derivatives by using a curvilinear mapping and a skew-symmetric splitting 
of the Cartesian derivatives. Next, the local Cartesian SBP operators are extended 
to global SBP operators by adding SATs to weakly couple elements at interfaces. 
To mimic the nonlinear integration by parts that occurs for the convective terms, \ie, 
\begin{equation*}
  \int_{\Omega}\sum\limits_{l=1}^{d}\W\Tr\frac{\partial\bmU}{\partial\xm{l}}\mr{d}\Omega=
  \int_{\Omega}\sum\limits_{l=1}^{d}\frac{\partial\Lyap}{\partial\xm{m}}\mr{d}\Omega=
  \oint_{\Gamma}\sum\limits_{l=1}^{d}\Lyap\nxm{m}\mr{d}\Gamma,
\end{equation*}
we introduce the additional mechanics of Tadmor's two-point flux functions and its 
extension to SBP operators via the Hadammard formalism.
\subsection{Computational space SBP operators}
We start by constructing one-dimensional matrix difference SBP operators in a fixed computational reference 
space $\xil{}\in[-1,1]$. The specific class of SBP operators we consider are 
 of diagonal-norm, one-dimensional collocated SBP operators defined as follows:
\begin{definition}\label{def:SBP1D}
    One-dimensional summation-by-parts operator for the first derivative:
    A matrix operator, $\Dxil{}\in\mathbb{R}^{N_{1d} \times N_{1d}}$, is a degree $p$ SBP operator approximating 
    the derivative $\partial/\partial\xil{}$ on the $N_{1d}$ node nodal distribution 
    $\bm{\xi}\equiv\left[-1,\xi_{2},\dots,\xi_{N_{1d}-1},1\right]\Tr$, if 
    \begin{enumerate}
        \item $\Dxil{}\bm{\xi}^{k}=k\bm{\xi}^{k-1}$, $k=0,\dots,p$;
        \item $\Dxil{}\equiv\M_{1d}^{-1}\Qxil{}$, where the norm matrix, $\M_{1d}$, is diagonal and symmetric positive definite; 
        \item $\Qxil{}\equiv\Sxil{}+\frac{1}{2}\Exil{}$, where $\Sxil{}=-\Sxil{}\Tr$, $\Exil{}=\Exil{}\Tr$;
        \item $\Exil{}\equiv \mr{diag}\left(-1,0,\dots,0,1\right)=\el{N_{1d}}\Tr\el{N_{1d}}-\el{1}\Tr\el{1}$, 
        $\el{N_{1d}}\equiv\left[0,\dots,0,1\right]^{T}$, and $\el{1}\equiv\left[1,0,\dots,0\right]^{T}$.

            \vspace{0.5cm}
  \noindent Moreover, $\bm{\xi}^{k}\equiv\left[-1^{k},\xi_{2}^{k},\dots,\xi_{N_{1d}-1}^{k},1^{k}\right]\Tr$, 
            $\xi_{i}$ is the $\xil{}$ coordinate of the $i^{th}$ node, and the convention that 
        $\bm{\xi}^{k}=\bm{0}$ if $k<0$ is used.
    \end{enumerate}
\end{definition}
The critical property that we want these operators to mimic discretely is the IBP property, 
\ie, 
\begin{equation}\label{eq:IBP1D}
  \int_{-1}^{1}\fnc{V}\Tr\frac{\partial\fnc{U}}{\partial \xil{}}\mr{d}\xil{}
  \int_{-1}^{1}\fnc{V}\Tr\frac{\partial\fnc{U}}{\partial \xil{}}\mr{d}\xil{}=
\left.\fnc{V}\fnc{U}\right|_{-1}^{1}.
\end{equation}
Given two scalar functions $U, V  \in L^{2}(\xil{})$, $\xil{}\in[-1,1]$, we first note that the 
norm matrix, $\M_{1d}$, is a high-order approximation to the $L_2$ inner product
of the two functions evaluated at the $N_{1d}$ node nodal distribution, \ie, 
\begin{equation}\label{eq:P1d}
  \bm{v}\Tr \, \M_{1d} \, \bm{u}\approx\int_{-1}^{1}\fnc{V} \, \fnc{U}\mr{d}\xil{}.
\end{equation}
Using~\eqref{eq:P1d} and the SBP derivative operators to discretize the left-hand side 
of~\eqref{eq:IBP1D} results in the following equality:
\begin{equation}\label{eq:SBP1D}
  \bm{v}\Tr\M_{1d}\Dxil{1}\bm{u}+\bm{u}\Tr\M_{1d}\Dxil{}\bm{v}=\bm{v}\Tr\Exil{}\bm{u}=
  \bm{v}(N)\bm{u}(N)-\bm{v}(1)\bm{u}(1),
\end{equation}
where the notation $\bm{v}(i)$ means the $i^{th}$ entry of the vector $\bm{v}$. 
Note that the terms in~\eqref{eq:SBP1D} are high-order approximations individually 
to those in~\eqref{eq:IBP1D}. Moreover, the operators mimic the IBP rule in a telescoping 
sense in that the surface integral on the right-hand side of~\eqref{eq:IBP1D} is split into individual 
surface contributions and, in the multidimensional case, into node-by-node contributions at surfaces. 
This property is important for the imposition of boundary conditions and minimal dissipation 
at interior element surfaces.  

The one-dimensional SBP operators are extended to multiple dimensions and systems of equations 
using Kronecker products. Thus, for example, if $dim=2$, there are $N_{1d}$ nodes in each direction, 
and there are $r$ equations then the operators for the $\xil{1}$ direction are  
\begin{equation}
  \begin{split}
      &\Dxil{1}\equiv\Dxil{}\otimes\Imat{N_{1d}}\otimes\Imat{r},\quad
  \M\equiv\M_{1d}\otimes\M_{1d}\otimes\Imat{r},\quad
  \Qxil{1}\equiv\Qxil{}\otimes\M_{1d}\otimes\Imat{r},\\
  &\Exil{1}\equiv\left(\el{N_{1d}}\Tr\el{N_{1d}}-\el{1}\Tr\el{1}\right)\otimes\Imat{N_{1d}}\otimes\Imat{r}=
  \left(\Rxiln{1}{N_{1d}}\right)\Tr\Morthoxil{1}\Rxiln{1}{N_{1d}}
 - \left(\Rxiln{1}{1}\right)\Tr\Morthoxil{1}\Rxiln{1}{1},\\
  &\Rxiln{1}{N_{1d}}\equiv\el{N_{1d}}^{T}\otimes\Imat{N_{1d}}\otimes\Imat{r},\quad
   \Rxiln{1}{1}\equiv\el{1}^{T}\otimes\Imat{N_{1d}}\otimes\Imat{r},\quad
   \Morthoxil{1}\equiv\M_{1d}\otimes\Imat{r},
  \end{split}
\end{equation}
where the generic matrix $\Imat{M}$ is an $M\times M$ identity matrix. Extending the one-dimensional SBP operators 
preserves the SBP property. Consider the integration by parts rule in multiple dimensions, which now reads
\begin{equation}\label{eq:IBP}
\int_{\hat{\Omega}}\bfnc{V}\Tr\frac{\partial\bmU}{\partial\xil{l}}\mr{d}\hat{\Omega}+
\int_{\hat{\Omega}}\bfnc{U}\Tr\frac{\partial\bfnc{V}}{\partial\xil{l}}\mr{d}\hat{\Omega}
=\oint_{\hat{\Gamma}}\bfnc{V}\Tr\bfnc{U}\nxil{l}\mr{d}\hat{\Gamma},
\end{equation}
where the computational domain in three dimensions is defined as $\hat{\Omega}\equiv[-1,1]^{3}$ and $\nxil{l}$ 
is the $l^{th}$ component of the outward facing unit normal on the surface of the computational domain 
$\hat{\Gamma}$. 

As in the one-dimensional case, the norm matrix, $\M$, defines a high-order discrete 
approximation to the $L^{2}$ inner product, \ie, 
\begin{equation}\label{eq:P}
  \bm{v}\Tr\M\bm{u}\approx\int_{\hat{\Omega}}\bfnc{V}\Tr\bmU\mr{d}\hat{\Omega}.
\end{equation}
Discretizing the left-hand side of~\eqref{eq:IBP} using~\eqref{eq:P} and the 
multidimensional SBP operators results in the following equality,
\begin{equation}\label{eq:SBP}
  \bm{v}\Tr\M\Dxil{l}\bm{u}+\bm{u}\Tr\M\Dxil{l}\bm{v} = \bm{v}\Tr\Exil{}\bm{u} = 
  \bm{v}\Tr\left(\Rxiln{l}{N_{1d}}\right)\Tr\Porthoxil{l}\Rxiln{l}{N_{1d}}\bm{u}-
  \bm{v}\Tr\left(\Rxiln{l}{1}\right)\Tr\Porthoxil{l}\Rxiln{l}{1}\bm{u}.
\end{equation}
%
%
\subsection{Local Cartesian SBP operators}
In the last section, we introduced SBP operators for fixed computational coordinates. In 
this section, we will use those operators to construct approximations to the first 
and the second derivatives in physical Cartesian coordinates that retain the SBP property. To do so 
consider the following expansion of the Cartesian derivatives in terms of the computational 
derivatives:
\begin{equation}\label{eq:cordtransformation}
  \frac{\partial}{\partial\xm{l}}=\sum\limits_{i=1}^{dim}
  \frac{\partial\xil{i}}{\partial\xm{l}}\frac{\partial}{\partial\xil{i}},\quad
   \frac{\partial}{\partial\xm{l}}=\sum\limits_{i=1}^{dim}\fnc{J}^{-1}
  \frac{\partial}{\partial\xil{i}}\fnc{J}\frac{\partial\xil{i}}{\partial\xm{l}},
\end{equation}
where we denote the first form as the chain-rule form and the second as the 
strong-conservation form. The second form is derived from the first by multiplying through by the metric 
Jacobian $\fnc{J}$, using chain rule, and the metric identities 
\begin{equation}\label{eq:GCL}
  \sum\limits_{l=1}^{dim}\frac{\partial}{\partial\xil{l}}
  \left(\fnc{J}\frac{\partial\xil{i}}{\partial\xm{l}}\right) = 0.
\end{equation}
Averaging the two expansions~\eqref{eq:cordtransformation} of the Cartesian 
coordinates, we obtain a third expansion, which we denote as the skew-symmetric form, 
namely
\begin{equation}\label{eq:cordtransformationskew}
  \frac{\partial}{\partial\xm{l}}=\frac{1}{2}\sum\limits_{i=1}^{dim}\left(
  \frac{\partial\xil{i}}{\partial\xm{l}}\frac{\partial}{\partial\xil{i}}+
   \fnc{J}^{-1}
  \frac{\partial}{\partial\xil{i}}\fnc{J}\frac{\partial\xil{i}}{\partial\xm{l}}\right).
\end{equation}
Before constructing local Cartesian SBP approximations, we first define a norm matrix
\begin{equation}\label{eq:PJ}
  \MJ\equiv\M\matJk{},
\end{equation}
where $\matJk{}$ is a diagonal matrix with the metric Jacobian $\fnc{J}$ along its 
diagonal. Since $\M$ is a diagonal positive-definite matrix, and $\matJk{}$ is a 
diagonal positive-definite matrix so too is $\MJ$. 

To construct a Cartesian SBP operator, we use our computational SBP operators 
to discretize the skew-symmetric form of the Cartesian operators, that is 
\begin{equation}\label{eq:Dxms}
  \Dxms{l}\equiv\frac{1}{2}\sum\limits_{i=1}^{dim}\left(
    \matJk{}^{-1}\Dxil{i}\matJxilxmk{i}{l}{}+
    \matJk{}^{-1}\matJxilxmk{i}{l}{}\Dxil{i}
    \right),
\end{equation}
where the superscript $s$ reminds us that this is the skew-symmetric approximation 
(we will use a different form for the second derivative). Using the norm matrix,
$\MJ$, 
we can construct the various matrices of the skew-symmetric Cartesian SBP operator as 
\begin{equation}\label{eq:matSBPC-1}
  \begin{split}
  &\Qxms{l}\equiv\MJ\Dxms{l}=
\frac{1}{2}\sum\limits_{i=1}^{dim}\left(
    \Qxil{i}\matJxilxmk{i}{l}{}+
    \matJxilxmk{i}{l}{}\Qxil{i}
    \right).
  \end{split}
\end{equation}
Thus, 
\begin{equation}\label{eq:matSBPC-2}  
  \begin{split}
  \Exm{l}\equiv\Qxms{l}+\left(\Qxms{l}\right)\Tr=&
  \sum\limits_{i=1}^{dim}\Exil{i}\matJxilxmk{i}{l}{}=\\
 &\sum\limits_{i=1}^{dim}\left(
  \left(\Rxiln{1}{N}\right)\Tr\Morthoxil{1}\Rxiln{1}{N}
 -
 \left(\Rxiln{1}{1}\right)\Tr\Morthoxil{1}\Rxiln{1}{1}\right)\matJxilxmk{i}{l}{},
  \end{split}
\end{equation}
and we can construct the skew symmetric part of the $\Qxms{l}$ matrix as 
\begin{equation}\label{eq:matSBPC-3}
  \begin{split}
  &\Sxms{l}\equiv\frac{1}{2}\sum\limits_{i=1}^{dim}\left(
    \Sxil{i}\matJxilxmk{i}{l}{}+
    \matJxilxmk{i}{l}{}\Sxil{i}
    \right).
  \end{split}
\end{equation}

Note that $\Exm{l}$ does not have a superscript because it will appear for 
the second derivative operator as well. 

Now, we demonstrate that 
these operators satisfy the SBP property. Discretizing the Cartesian version of IBP~\eqref{eq:IBP} 
using the operators given in \eqref{eq:Dxms}-\eqref{eq:matSBPC-2}  yields the equality
\begin{equation}\label{eq:SBPC}
  \bm{v}\Tr\MJ\Dxms{l}\bm{u}+\bm{u}\Tr\MJ\Dxms{l}=\bm{v}\Tr\Exm{l}\bm{u}=
\sum\limits_{i=1}^{dim}\bm{v}\Tr\left(
  \left(\Rxiln{i}{N_{1d}}\right)\Tr\Morthoxil{i}\Rxiln{i}{N_{1d}}
 - \left(\Rxiln{i}{1}\right)\Tr\Morthoxil{i}\Rxiln{i}{1}\right)\matJxilxmk{i}{l}{}\bm{u},
\end{equation}
where we see that not only do we mimic IBP rule term-by-term and to high order, 
but the surface integrals are composed of individual contributions at surface nodes. 

Before presenting the discrete operator for the second derivative, it is instructive to 
examine the IBP property that we intend to preserve, namely 
\begin{equation}\label{eq:IBP2}
 \int_{\Omega}\bfnc{V}\Tr\frac{\partial}{\partial\xm{l}}\left(\Vhatmn{l}{m}\frac{\partial\bmU}{\partial\xm{m}}\right)\mr{d}\Omega = 
 \oint_{\Gamma}\bfnc{V}\Tr\Vhatmn{l}{m}\frac{\partial\bmU}{\partial\xm{m}}\nxm{l}\mr{d}\Gamma
 -\int_{\Omega}\frac{\partial\bfnc{V}}{\partial\xm{l}}\Tr\Vhatmn{l}{m}\frac{\partial\bmU}{\partial\xm{m}}\mr{d}\Omega.
\end{equation}
While we could use the skew-symmetric formulation, we instead discretize using a combination of 
the strong conservation form, which is discretized as
\begin{equation}\label{eq:Dxmstrong}
  \Dxmsc{l}\equiv\sum\limits_{i=1}^{dim}\matJk{}^{-1}\Dxil{i}\matJxilxmk{i}{l}{},
\end{equation}
and the chain rule form, which is discretized as  
\begin{equation}\label{eq:Dxmchain}
  \Dxmc{l}\equiv\sum\limits_{i=1}^{dim}\matJk{}^{-1}\matJxilxmk{i}{l}{}\Dxil{i},
\end{equation}
where the superscript $sc$ and $c$ are used to denote that these are discretizations for the 
strong conservation form and the chain-rule form, respectively. Note that these operators  
are not SBP by themselves. However, their combination for discretizing the
second derivative satisfies the SBP property. Thus, we discretize the second derivative as  
\begin{equation}\label{eq:D2V}
  \Dxmsc{l}\matVhatmnk{l}{m}{}\Dxmc{m}\bm{u}\approx
  \frac{\partial}{\partial\xm{l}}\left(\Vhatmn{l}{m}\frac{\partial\bmU}{\partial\xm{m}}\right),
\end{equation} 
where $\matVhatmnk{l}{m}{}$ is a block diagonal matrix with $\Vhatmn{l}{m}$ evaluated at 
each node. Discretizing the left-hand side of~\eqref{eq:IBP2} using the operator
\eqref{eq:D2V} and $\MJ$ defined in \eqref{eq:PJ}
gives the following equality:
\begin{equation}\label{eq:SBP2}
  \bm{v}\Tr\MJ\Dxmsc{l}\matVhatmnk{l}{m}{}\Dxmc{m}\bm{u} = \bm{v}\Tr\Exm{l}\bm{u}
  -\bm{v}\Tr\left(\Dxmc{l}\right)\Tr\matVhatmnk{l}{m}{}\Dxmc{m}\bm{u}.
\end{equation}
%
%
\subsection{Global Cartesian SBP operators}\label{appA:globalCartSBP}

In the previous section, we constructed Cartesian SBP operators that act at the element 
In the previous section, we constructed Cartesian SBP operators that act at the element 
level and preserve the SBP property. However, to construct stability proofs of the theorems presented in Section \ref{sec:dg-sbp-sat-global}, we need SBP operators over the entire set of 
elements used to discretize the domain $\Omega$ with boundary $\Gamma$.  We can construct such global SBP operators 
by using appropriate SATs. In this section, we demonstrate the procedure using simple 
two-element examples, i.e., left and right elements denoted with the superscripts $L$ and
$R$, respectively. 

For this purpose, we consider two elements with a common vertical surface. For simplicity 
of presentation, we consider that all coordinate systems are aligned, \ie, the Cartesian 
and computational coordinates are aligned within each element and are aligned between 
elements (we don't mean parallel, but rather that they generally point in the same direction).
For simplicity of explanation but without loss of generality, suppose that the 
SBP operator $\Dxmsk{l}{L}$ is applied to the vector vector $\bm{u}^{L}$, i.e., $\Dxmsk{l}{L}\bm{u}^{L}$. 
Thus, we 
modify the Cartesian derivative of the left element, $\Dxmsk{l}{L}\bm{u}^{L}$, with an interface SAT coupling it to 
the right element as follows:
\begin{equation}\label{eq:DxmL}
  \Dxmsk{l}{L}\bm{u}^{L}+\bm{SAT}_{L}\approx\frac{\partial\bmU_{L}}{\partial \xm{m}},
\end{equation}
where 
\begin{equation}\label{eq:SATL}
    \bm{SAT}^{L}\equiv-\frac{1}{2}\left(\MJ^{L}\right)^{-1}\left(\Rxiln{1}{N_{1d}}\right)\Tr\Porthoxil{1}
  \left(\Rxiln{1}{N_{1d}} \matJxilxmk{1}{l}{L}\bm{u}_{L}
  -\Rxiln{1}{1}\matJxilxmk{1}{l}{R}\bm{u}_{R}\right).
\end{equation}
We proceed simiarily for the right element and obtain
\begin{equation}\label{eq:DxmR}
  \Dxmsk{l}{R}\bm{u}_{R}+\bm{SAT}_{R}\approx\frac{\partial\bmU_{R}}{\partial \xm{m}},
\end{equation}
where 
\begin{equation}\label{eq:SATR}
  \bm{SAT}^{R}\equiv\frac{1}{2}\left(\MJ^{R}\right)^{-1}\left(\Rxiln{1}{1}\right)\Tr\Porthoxil{1}
  \left(\Rxiln{1}{1}\matJxilxmk{1}{l}{R}\bm{u}_{R}
  -\Rxiln{1}{N_{1d}}\matJxilxmk{1}{l}{L}\bm{u}_{L}\right).
\end{equation}

Now, we construct a global SBP operator over the two elements and obtain the following:
\begin{equation}\label{eq:Dxmg}
  \Dxmsk{l}{g}\equiv
  \overbrace{
  \left[
    \begin{array}{cc}
      \MJ^{L}\\
      &\MJ^{R}
    \end{array}
  \right]^{-1}
  }^{\equiv\MJ^{g}}
  \left(
    \overbrace{
  \left[
\begin{array}{cc}
  \Sxmsk{l}{L}&\bm{CT}_{L}\\
\bm{CT}_{R}
  &\Sxmsk{l}{R}
\end{array}
  \right]}^{\equiv\Sxmsk{l}{g}}+\frac{1}{2}
  \overbrace{
  \left[
\begin{array}{cc}
  \Exmk{l}{L}-\bm{IS}_{L}\\
  &\Exmk{l}{R}-\bm{IS}_{R}
\end{array}
  \right]}^{\equiv\Exmk{l}{g}}
  \right),
  \end{equation}
where $\Sxmsk{l}{L/R}$ are constructed using \eqref{eq:matSBPC-3}, and  
\begin{equation}\label{eq:CT}
  \begin{split}
  &\bm{CT}_{L}\equiv  \frac{1}{2}\left(\Rxiln{1}{N_{1d}}\right)\Tr\Porthoxil{1}
  \Rxiln{1}{1}\matJxilxmk{1}{l}{R},\\
  &\bm{CT}_{R}\equiv-\frac{1}{2}\left(\Rxiln{1}{1}\right)\Tr\Porthoxil{1}
 \Rxiln{1}{N_{1d}}\matJxilxmk{1}{l}{L}.
  \end{split}
\end{equation}
We note that for consistent approximations the metric terms $\matJxilxmk{1}{l}{L}$ 
and $\matJxilxmk{1}{l}{R}$ are the same at the joined boundary points. 
Then, $\bm{CT}_{R}=-\bm{CT}_{R}$ and therefore $\Sxmsck{l}{g}$ is skew-symmetric. 
Finally, 
\begin{equation}\label{eq:Exmgex-1}
  \begin{split}
      &\bm{IS}_{L}\equiv\frac{1}{2}\left(\Rxiln{1}{N_{1d}}\right)\Tr\Porthoxil{1}
      \Rxiln{1}{N_{1d}}\matJxilxmk{1}{l}{L},\\
  &\bm{IS}_{R}\equiv-\frac{1}{2}\left(\Rxiln{1}{1}\right)\Tr\Porthoxil{1}
  \Rxiln{1}{1}\matJxilxmk{1}{l}{R}.
  \end{split}
\end{equation}
We also note that the action of the $\bm{IS}$ terms is to remove the contribution of the
$\Exmk{}{}$ matrices at the interior interface so that the symmetric matrix $\Exmk{l}{g}$ 
only contains terms that contribute to the boundaries of the two-element mesh. 

This procedure can be applied at all interior interfaces and the resultant 
global operators preserve the SBP property and have the following matrix 
  properties:
  \begin{equation}\label{eq:Dxmscg}
    \Dxmsk{l}{g}=\left(\MJ^{g}\right)^{-1}\Qxmsk{l}{g},\quad
    \Qxmsk{l}{g}+\left(\Qxmsk{l}{g}\right)\Tr=\Exmk{l}{g}.
  \end{equation}
The matrix $\Exmk{l}{g}$ is a surface integration matrix such that for $M$ surfaces it has the following form, 
  \begin{equation}\label{eq:Exmgapp}
    \Exmk{l}{g}\equiv\sum\limits_{i=1}^{M}n_{\xi_{\Gamma_i}}\left(\Rxis{\Gamma_i}\right)\Tr\Porthoxil{\Gamma_i}\Rxis{\Gamma_i}
    \matJxilxmk{\Gamma_i}{m}{},
  \end{equation}
  where $\Gamma_{i}$ specifies the computational coordinate orthogonal to the $i^{\text{th}}$ surface (i.e., in three dimensions 
 $\Gamma_{i}$ can be $1$ or $2$ or $3$).  For example, suppose the surface is 
  is orthogonal to $\xil{1}$ and at the maximum of $\xil{1}$ then 
  $\Rxis{\Gamma_{i}}=\Rxiln{1}{N}$, $\Porthoxil{\Gamma_i}=\Porthoxil{1}$, and $\matJxilxmk{\Gamma_i}{}{}=\matJxilxmk{1}{}{}$.

  We proceed in an identical manner for the approximation to the second derivative and 
  in the same way as at the element level, we define two global matrices $\Dxmsck{l}{g}$ 
  and $\Dxmck{m}{g}$ and construct the global SBP approximation to the second derivative as 
  \begin{equation}\label{eq:D2g}
    \Dxmsck{l}{g}\matVhatmnk{l}{m}{g}\Dxmck{m}{g}=\left(\MJ^{g}\right)^{-1}
    \left(-\left(\Dxmck{l}{g}\right)\Tr\MJ^{g}\matVhatmnk{l}{m}{g}\Dxmck{m}{g}+
    \Exmk{l}{g}\matVhatmnk{l}{m}{k}\Dxmck{m}{g}\right),
  \end{equation}
  where the right-hand side decomposition demonstrates the SBP property of this operator.
%
%
\subsection{Convective term: Tadmor's two-point flux functions and the Hadammard formalism}\label{app:Hadammard}
To mimic the continuous proof at the semi-discrete level, we need approximations 
to the convective terms that reduce to discrete surface integrals when discretely integrated against 
the Lyapunov variables. In particular, we rely on combining SBP operators with 
two-point flux functions of Tadmor~\cite{tadmor2003entropy}. The required mechanics are involved, and we 
attempt to convey all the details clearly. However, the interested reader is referred to 
the literature on entropy-stable methods.  

To begin, we introduce the following scalar function, called Laypunov potential,
\begin{equation}\label{eq:Lappotential}
  \fnc{\psi}\equiv\W\Tr\bmU-\Lyap,
\end{equation}
where $\W$, $\bmU$, and $V$ are the Lyapunov variables, the conserved variables, and
the Lyapunov function, respectively. Next, we introduce the two-point flux functions, $\FSCempty$, that 
satisfies the following three conditions:
\begin{enumerate}
    \item Consistency: $\USC{L}{j}{L}{j}= \bm{F}^{(c)}\left(\bm{u}^{L}(j,:)\right)$,
  \item  Tadmor shuffle condition for the convective spatial flux terms~\cite{tadmor2003entropy},
\begin{equation}\label{eq:tadmor}
  \left(\wki{L}{i}-\wki{R}{j}\right)\Tr\USC{L}{i}{R}{j}=\psiki{L}{i}-\psiki{R}{j}
\end{equation}
\item Symmetry: $\USC{L}{i}{R}{j}=\USC{R}{j}{L}{i}$,
\end{enumerate}
where $\uk{L}$ and $\uk{R}$ are the discrete solution vectors for an element to the left 
and right of an interface, respectively, with similar interpretations for 
$\wk{L}$, $\wk{R}$, $\psik{L}$, and $\psik{R}$, where the latter are constructed from 
the $\uk{}$ vectors using \eqref{eq:Lappotential}. The notation, for example, $\uki{L}{i}$ means the $r\times1$ solution 
vector at the $i^{\text{th}}$ node where $r$ is the number of equations in the PDE.

To achieve Lyapunov consistency, the approximations to the spatial derivatives of
the convective flux are constructed such that the continuous Lyapunov analysis is
mimicked in a term-by-term in the semi-discrete context. To this end, the approximation
to the spatial derivatives has to be mimetic to a special case of the integration by parts, 
namely
\begin{equation}\label{eq:nonlinearIBP}
  \begin{split}
      \int_{D} \bm{Q}\Tr \frac{\partial \bm{G}}{\partial \chi} \, dD = \int_{D} \frac{\partial h}{\partial \chi} \, dD = \oint_{B} h \, n_{\chi}  dB,
  \end{split} 
\end{equation}
where where $\bm{Q}$ and $\bm{G}(\bm{U})$ are vector valued functions, $h(\bm{U})$ is a scalar function, 
$\chi$ is some independent variable over the domain $D$ with boundary $B$, and $n_{\chi}$ 
is the component of the outward facing unit normal in the $\chi$ direction.
In the present analysis, $\bm{Q}$, $\bm{G}(\bm{U})$, and $h$ play the role of the Lyapunov variables,
$\W$, the convective fluxes, $\bmfI{m}$, or the conservative variables, $\bmU$, and the Lyapunov fluxes, $\bmF{m}$, 
or the Lyapunov function, $\Lyap$, respectively. Equation \eqref{eq:nonlinearIBP} is referred to as a nonlinear IBP relation 
and the numerical operators that mimic this property discretely are called nonlinear SBP operators.

The nonlinear SBP operators constructed and used in this work to discretize the divergence
of the convective and diffusive fluxes for the system of PDEs \eqref{eq:SEIR_con_diff} satisfy 
the following approximations \cite{fernandez2020entropy}, 
\begin{equation}\label{eq:nonlinearSBP}
  \begin{split}
  &2\left(\Dxil{i}\matJxilxmk{i}{l}{k}\right)\circ\matUSC{k}{k}\ones\approx
\left.\left(\frac{\partial}{\partial\xil{i}}\left(\Jxilxm{i}{l}\bm{F}^{(c)}\right)+
      \bm{F}^{(c)}\frac{\partial}{\partial\xil{i}}\left(\Jxilxm{i}{l}\right)\right)\right|_{C_k},\\
  &2\left(\matJxilxmk{i}{l}{k}\Dxil{i}\right)\circ\matUSC{k}{k}\ones\approx
\left.\left(\frac{\partial}{\partial\xil{i}}\left(\Jxilxm{i}{l}\bm{F}^{(c)}\right)\right)\right|_{C_k},
  \end{split}
\end{equation}
where $\circ$ represents the Hadammard product (i.e., the entry-wise multiplication of two matrices), and the 
notation ``$\left( a \right)|_{c_{k}}$" indicates the vector constructed by evaluating the quantity $a$ at the 
set of nodes $C$ of the element $k$, i.e., $C_{k}$. The proofs of \eqref{eq:nonlinearSBP} is given in 
\cite{fernandez2020entropy} and follow in a straightforward manner 
from the Theorem 1 in \cite{Crean2018}.
The block diagonal matrix \matUSC{k}{k}, for the $k^{th}$ element having discrete 
solution vector $\uk{k}$, is constructed as 
\begin{equation}\label{eq:matUSC}
  \begin{split}
&\matUSC{k}{k}\equiv\\
&\left[
\begin{array}{ccc}
    \mr{diag}\left(\USC{k}{1}{k}{1}\right)&\dots&\mr{diag}\left(\USC{k}{1}{k}{N_{1d}}\right)\\
\vdots                                &     & \vdots\\
\mr{diag}\left(\USC{k}{N_{1d}}{k}{1}\right)&\dots&\mr{diag}\left(\USC{k}{N_{1d}}{k}{N_{1d}}\right)
\end{array}
\right],
\end{split}
\end{equation}
where $N_{1d}^{dim}$ (i.e., $N_{1d}$ to the power $dim$) is the total number of nodes in the element $k$. Note that via the properties of 
$\FSCempty$ 
the matrix $\matUSC{k}{k}$ is symmetric. 

Our construction of a global nonlinear SBP operator follows identically as before, and we obtain the 
approximation \cite{Fisher2013b,fernandez2020entropy,fernandez_entropy_stable_hp_ref_snpdea_2019,Crean2018}
\begin{equation}
  2\Dxmsk{l}{g}\circ\matUSC{g}{g}\bm{1}\approx\frac{\partial\bm{F}^{(c)}}{\partial \xm{l}}\left(\bm{x}_{g}\right),
\end{equation}
where the notation $\frac{\partial\bm{F}^{(c)}}{\partial \xm{l}}\left(\bm{x}_{g}\right)$ means the vector 
that results from evaluating $\frac{\partial\bm{F}^{(c)}}{\partial \xm{l}}$ at the mesh nodes. In addition, 
the resulting (global) nonlinear SBP operator satisfies the discrete counterpart of the 
nonlinear IBP rule \eqref{eq:nonlinearIBP}; see \cite{Fisher2013b,fernandez2020entropy,fernandez_entropy_stable_hp_ref_snpdea_2019} and also
Theorem 4.2 in \cite{yamaleev2019entropy}.

For our discretization to result in a Lyapunov-consistent formulation, we need 
our global SBP operator to discretely satisfy the metric invariants  (similarly to the entropy conservative and 
stable discretizations proposed in \cite{fernandez_entropy_stable_hp_ref_snpdea_2019}), which are 
\begin{equation}\label{eq:GCL}
  \sum\limits_{i=1}^{dim}\frac{\partial}{\partial\xil{i}}\left(\fnc{J}\frac{\partial\xil{i}}{\partial\xm{l}}\right) = 0,\quad l=1,\dots,dim.
\end{equation}
The discrete counter-part to~\eqref{eq:GCL} turns out to be \cite{fernandez_entropy_stable_hp_ref_snpdea_2019,friedrich2018entropy} 
\begin{equation}\label{eq:DGCLxm}
 \Dxmsk{l}{g}\bm{1} = 0.
\end{equation}
To understand why, note that the SATs used to construct $\Dxmsk{l}{g}$ in Appendix \ref{appA:globalCartSBP} cancel for the 
constant vector. Thus, we have the following local statements
\begin{equation}\label{eq:DGCL}
  \sum\limits_{i=1}^{dim}\Dxil{i}\matJxilxmk{i}{l}{}\bm{1}+\matJxilxmk{i}{l}{}\Dxil{i}=
  \sum\limits_{i=1}^{dim}\Dxil{i}\matJxilxmk{i}{l}{}\bm{1}=0.
\end{equation}
In this work, the metric terms for conforming and non-conforming elements with 
mesh- ($h$-) and solution polynomial degree ($p$-) refinements are constructed to satisfy~\eqref{eq:DGCL} as 
reported in \cite{fernandez_entropy_stable_hp_ref_snpdea_2019}. Therefore, \eqref{eq:DGCLxm} holds.  

We need one last theorem before we can show why the involved discretization 
of the convective terms mimics the continuous proof of the Lyapunov consistency. This theorem is given next.

\vspace{0.5cm}

\begin{theorem}\label{thrm:shuffel}
  Consider an arbitrary matrix $\bar{\mat{A}}$ of size $N^{dim}\times N^{dim}$ and the matrix 
  $\matUSC{L}{R}$ constructed from a two-point flux function $\USC{L}{i}{R}{j}$, 
  that is symmetric and satisfies the Tadmor's shuffle condition \cite{tadmor2003entropy}. Then, for 
  $\mat{A}\equiv\bar{\mat{A}}\otimes\Imat{r}$
  \begin{equation}\label{eq:shuffle}
      \left(\wk{L}\right)\Tr\left(\mat{A}\circ\matUSC{L}{R}\right)\ones-\ones\Tr\left(\mat{A}\circ\matUSC{L}{R}\right)\wk{R} = 
      \left({\psik{L}}\right)\Tr\bar{\mat{A}}\onesbar-\onesbar\Tr\bar{\mat{A}}\psik{R}.
  \end{equation}
\end{theorem}
\begin{proof}
    The proof is given in \cite{fernandez2020entropy} (see also Ref. \cite{Crean2018} Lemmas 2and 3).
\end{proof}

We can finally prove Theorem~\ref{thrm:SBPnon}. To begin, since $\MJ^{g}$ is diagonal we obtain
  \begin{equation}\label{eq:ent1}
      2\left(\wk{g}\right)\Tr\MJ^{g}\Qxmsk{l}{g}\circ\matUSC{g}{g}\bm{1}= 
2\left(\wk{g}\right)\Tr\Qxmsk{l}{g}\circ\matUSC{g}{g}\bm{1}.
  \end{equation}
  Next splitting the right-hand side of~\eqref{eq:ent1} into two terms, taking the transpose 
  of the second term and using the symmetry of the matrix $\matUSC{g}{g}$, we obtain
  \begin{equation}\label{eq:ent2}
\left(\wk{g}\right)\Tr\MJ^{g}\Dxmsk{l}{g}\circ\matUSC{g}{g}\bm{1}= 
\left(\wk{g}\right)\Tr\Qxmsk{l}{g}\circ\matUSC{g}{g}\bm{1}+\bm{1}\Tr\left(\Qxmsk{l}{g}\right)\Tr\circ\matUSC{g}{g}\wk{g}.
  \end{equation}
Using the SBP property, $\Qxmsk{l}{g}=-\left(\Qxmsk{l}{g}\right)\Tr+\Exmk{l}{g}$, yields
  \begin{equation}\label{eq:ent2}
    \begin{split}
2\left(\wk{g}\right)\Tr\MJ^{g}\Qxmsk{l}{g}\circ\matUSC{g}{g}\bm{1}=& 
\left(\wk{g}\right)\Tr\Qxmsk{l}{g}\circ\matUSC{g}{g}\bm{1}-\bm{1}\Tr\Qxmsk{l}{g}\circ\matUSC{g}{g}\wk{g}\\
&+\bm{1}\Tr\Exmk{l}{g}\circ\matUSC{g}{g}\wk{g}.
    \end{split}
  \end{equation}
  Applying Theorem~\ref{thrm:shuffel} to the first two terms on the right-hand side
  of~\eqref{eq:ent2} leads to
  \begin{equation}\label{eq:ent3}
    \begin{split}
2\left(\wk{g}\right)\Tr\MJ^{g}\Qxmsk{l}{g}\circ\matUSC{g}{g}\bm{1}=& 
        \left({\psik{g}}\right)\Tr\Qxmskbar{l}{g}\bar{\bm{1}}-\bar{\bm{1}}\Tr\Qxmskbar{l}{g}\psik{g}\\
&+\bm{1}\Tr\Exmk{l}{g}\circ\matUSC{g}{g}\wk{g}.
    \end{split}
  \end{equation}
 If~\eqref{eq:DGCLxm} holds, then, $\Qxmskbar{l}{g}\bar{\bm{1}}=\bm{0}$. Furthermore, 
  by a similar reasoning we have
  \begin{equation*}
  \bar{\bm{1}}\Tr\Qxmskbar{l}{g}=
  \bar{\bm{1}}\Tr\left(-\left(\Qxmskbar{l}{g}\right)\Tr+\Exmkbar{l}{g}\right)
  =\bar{\bm{1}}\Tr\Exmkbar{l}{g}.
  \end{equation*}
  Thus, Equation~\eqref{eq:ent3} reduces to 
  \begin{equation}\label{eq:ent4}
    \begin{split}
2\left(\wk{g}\right)\Tr\MJ^{g} \, \Qxmsk{l}{g}\circ\matUSC{g}{g}\bm{1}=&\bm{1}\Tr\Exmk{l}{g}\circ\matUSC{g}{g}\wk{g}-\bar{\bm{1}}\Tr
\, \Exmkbar{l}{g}\psik{g}.
    \end{split}
  \end{equation}
The right-hand side of Equation \eqref{eq:ent4} contains only boundary terms arising 
from the convective term of the PDE system.
For the particular class of SBP operators, we are examining, $\Exmk{l}{g}$ is diagonal 
  and has nonzeros associated with boundary nodes. Moreover, since $\matUSC{g}{g}\wk{g}$ is 
  consistent, the block diagonals associated with boundary nodes are simply the diagonal 
  matrix $\left[\bm{F}^{(c)}\right]$.  
%
%
\section{Curvilinear element-wise discretization}\label{app:discelement}

The discretization on the element with index $k$ is given as 
  \begin{equation}\label{eq:element}
    \begin{split}
    \matJk{k}\frac{\mr{d}\uk{k}}{\mr{d}t} =& 
    \fk{k}
    -\sum\limits_{l,m=1}^{dim}\left(\Dxil{m}\matJxilxmk{m}{l}{k}+\matJxilxmk{m}{l}{k}\Dxil{m}\right)\circ\matUSC{k}{k}\ones
    +\sum\limits_{i,j=1}^{dim}\Dxil{i}\matVtildelmk{i}{j}{k}\Thetajk{j}{k}\\
        &+\bm{SAT}^{(c)}+\bm{SAT}^{(d)}+\bm{SAT}^{(bc)}+\bm{diss}^{(c)}+\bm{diss}^{(d)},
    \end{split}
  \end{equation}
  where 
  \begin{equation}
\matVtildelmk{i}{j}{k}\equiv\sum\limits_{l,m=1}^{dim}\matJxilxmk{i}{l}{k}\matVhatmnk{l}{m}{k}\matJk{k}^{-1}\matJxilxmk{i}{l}{k}.
  \end{equation}
  Moreover, $\Thetajk{j}{k}$ is defined as 
  \begin{equation}\label{eq:theta}
      \Thetajk{j}{k}=\Dxil{j}\wk{k}+\bm{SAT}^{(w)}.
  \end{equation}

  The interface dissiaption terms, $\bm{diss}^{(c)}$ and $\bm{diss}^{(d)}$, are discussed in Appendix~\ref{app:diss}.  
  Below, we describe the remaining SAT terms, where the SATs are constructed for the surface perpendicular to $\xil{1}$ and located where $\xil{1}$ is at a maximum. 
  The $\bm{SAT}^{(c)}$ is defined as  
  \begin{equation*}
      \begin{split}
          \bm{SAT}^{(c)}\equiv\M^{-1}\sum\limits_{m=1}^{dim}&\left\{
      \left(\left(\Rxiln{1}{N_{1d}}\right)\Tr\Morthoxil{1}\Rxiln{1}{N_{1d}}\matJxilxmk{1}{m}{k}\right)\circ\matUSC{k}{k}\ones \right.\\
          &\left.-\left(\left(\Rxiln{1}{N_{1d}}\right)\Tr\Morthoxil{1}\Rxiln{N_{1d}}{1}\matJxilxmk{1}{m}{a}\right)\circ\matUSC{k}{a}\ones\right\},
      \end{split}
  \end{equation*}
  where $\uk{a}$ and $\matJxilxmk{1}{m}{a}$ are the solution and the metrics terms of the adjoining element, respectively.

The $\bm{SAT}^{(d)}$ term reads 
  \begin{equation*}
      \bm{SAT}^{(d)}\equiv-\frac{1}{2}\M^{-1}\sum\limits_{j=1}^{dim}\left(
      \left(\Rxiln{1}{N_{1d}}\right)\Tr\Morthoxil{1}\Rxiln{1}{N_{1d}}\matVtildelmk{1}{j}{k}\Thetajk{j}{k}\ones
      -\left(\Rxiln{1}{N_{1d}}\right)\Tr\Morthoxil{1}\Rxiln{N_{1d}}{1}\matVtildelmk{1}{j}{a}\Thetajk{j}{a}\right),
  \end{equation*}
  where $\Thetajk{j}{a}$ and $\matVtildelmk{1}{j}{a}$ are from the adjoining element. 

Finally, $\bm{SAT}^{(w)}$ indicates the SAT term applied to the gradient of the discrete entropy variables. 
For example, the $\bm{SAT}^{(w)}$ for $\Thetajk{1}{k}$ is constructed as 
\begin{equation}\label{eq:SATWj}
    \bm{SAT}^{(w)} \equiv-\frac{1}{2}\M^{-1}\left(\Rxiln{1}{N_{1d}}\right)\Tr\Morthoxil{1}
    \left(\Rxiln{1}{N_{1d}}\wk{k}-\Rxiln{}{1}\wk{a}\right),
\end{equation}
where $\wk{a}$ is from the adjoining element. 
%
%
\section{Interface dissipation}\label{app:diss}
The $\bm{diss}^{(d)}$ term adds interface dissipation for the viscous terms and is constructed as 
\begin{equation}\label{eq:SATIP}
    \bm{diss}^{(d)}\equiv-\M^{-1}\left(\Rxiln{1}{N_{1d}}\right)\Tr\Morthoxil{1}
  \widetilde{\mat{Co}}
    \left(\Rxiln{1}{N_{1d}}\wk{k}-\Rxiln{1}{1}\wk{a}\right),
\end{equation}
where $\widetilde{\mat{Co}}$ is an appropriate combination of interface values of 
$\matVtildelmk{1}{j}{k}$ and $\matVtildelmk{1}{j}{a}$ (for example the average) and 
the terms $\wk{a}$ and $\matVtildelmk{1}{j}{a}$ originate from the adjoining element.

We continue to use the same interface for the analysis as in the previous sections. 
The $\bm{diss}^{(d)}$ term from the $L$ element results in  
\begin{equation}\label{eq:STLIP}
    diss^{(d),L}\equiv-\left(\wk{L}\right)\Tr\left(\Rxiln{1}{N_{1d}}\right)\Tr\Porthoxil{1}\widetilde{\mat{Co}}
    \left(\Rxiln{1}{N_{1d}}\wk{L}-\Rxiln{1}{1}\wk{R}\right),
\end{equation}
while for the $R$ element we obtain
\begin{equation}\label{eq:STRIP}
    diss^{(d),R}\equiv-\left(\wk{R}\right)\Tr\left(\Rxiln{1}{1}\right)\Tr\Porthoxil{1}\widetilde{\mat{Co}}
    \left(\Rxiln{1}{1}\wk{R}-\Rxiln{1}{N_{1d}}\wk{L}\right).
\end{equation}
Summing~\eqref{eq:STLIP} and~\eqref{eq:STRIP} gives
\begin{equation*}
    diss^{(d),L}+diss^{(d),R} = -\left(\Rxiln{1}{N_{1d}}\wk{L}-\Rxiln{1}{1}\wk{R}\right)\Tr\Porthoxil{1}\widetilde{\mat{Co}}
    \left(\Rxiln{1}{N_{1d}}\wk{L}-\Rxiln{1}{1}\wk{R}\right),
\end{equation*}
which is a negative semi-definite term and therefore does not impact the
Lyapunov consistency and stability statement.

Also the Lyapunov consistent SATs terms, $SAT^{(c)}$, are augmented with dissipative
terms inspired by the upwinding used in the Roe approximate Riemann solver, which has the
following generic form:
\begin{equation}\label{eq:roe-upwind-generic}
\boldsymbol{F}^{*}=\frac{\boldsymbol{F}^{+}+\boldsymbol{F}^{-}}{2}-\frac{1}{2} \mathrm{Y}|\Lambda| \mathrm{Y}^{-1}\left(\boldsymbol{U}^{+}-\boldsymbol{U}^{+}\right)
\end{equation}
where the superscript ``+" denotes quantities evaluated on the side of an interface 
in the positive direction of the unit normal, and the superscript ``-" denotes 
quantities in the negative direction. Since the Roe flux is not entropy consistent, the central flux in
\eqref{eq:roe-upwind-generic} (i.e., the first term on the right-hand side) is replaced by the numerical fluxes given in the above Lyapunov conservative $SAT^{(c)}$ 
terms. To provide entropy dissipation at element interfaces, the dissipative term, i.e., the second term
on the right-hand side of \eqref{eq:roe-upwind-generic} also has to be modified. The approach is constructing a dissipation term that enables Lyapunov stability 
while remaining accurate and conservative in design order. This last step is accomplished by 
using the flux Jacobian with respect to the Lyapunov variables rather than the conservative variables \cite{Fisher2013b,Carpenter2014}.
The eigenvectors of the conservative variable flux Jacobian can be scaled such that \cite{merriam1989entropy}, 
$$
\frac{\partial \boldsymbol{U}}{\partial \boldsymbol{W}}=Y Y^{\mathrm{T}}.
$$
Using this approach, after some algebraic manipulation and including the metric terms for curvilinear grids, 
we arrive at the following expression for the $\bm{diss}^{(c)}$,
\begin{equation*}
\bm{diss}^{(c)}\equiv
    -\frac{1}{2}\M^{-1}\left(\Rxiln{1}{N_{1d}}\right)\Tr\Porthoxil{1}\matdutildedwjk{1}{}
    \left(\Rxiln{1}{N_{1d}}\wk{k}-\Rxiln{1}{1}\wk{a}\right),
\end{equation*}
where $\matdutildedwjk{1}{}$ is a block diagonal matrix constructed from an appropriate average of 
\begin{equation*}
  \frac{\partial\tilde{\bmU}_{m}}{\partial \W}\equiv\sum\limits_{l=1}^{dim}\fnc{J}\frac{\partial\xil{m}}{\partial\xm{m}}\frac{\partial\bmU}{\partial \W},
\end{equation*}
from element $k$ and the adjecent element denoted by $a$. 
We highlight that, the same approach can be used for dynamic, unstructured grids. 
We refer the reader to the detailed and clear analysis reported 
by Yamaleev and co-authors \cite{yamaleev2019entropy}.

Now, we compute the resulting term of $\bm{diss}^{(c)}$ for an interface shared by two elements denoted by the usual
superscripts ``$L$" and ``$R$".
The contribution from the element $L$ is given by
\begin{equation}\label{eq:dissconL}
    diss^{(c),L} \equiv -\wk{L}\left(\Rxiln{1}{N_{1d}}\right)\Tr\Porthoxil{1}\matdutildedwjk{1}{}
    \left(\Rxiln{1}{N_{1d}}\wk{L}-\Rxiln{1}{1}\wk{R}\right),
\end{equation}
while for the $R$ element the contribution is
\begin{equation}\label{eq:dissconR}
    diss^{(c),R} \equiv -\wk{R}\left(\Rxiln{1}{1}\right)\Tr\Porthoxil{1}\matdutildedwjk{1}{}
    \left(\Rxiln{1}{1}\wk{R}-\Rxiln{1}{N_{1d}}\wk{L}\right).
\end{equation}
Adding~\eqref{eq:dissconL} to~\eqref{eq:dissconR}  gives
\begin{equation*}
    diss^{(c),L} + diss^{(c),R} = -\left(\Rxiln{1}{N_{1d}}\wk{L}-\Rxiln{1}{1}\wk{R}\right)\Tr\Porthoxil{1}\matdutildedwjk{1}{}
    \left(\Rxiln{1}{N_{1d}}\wk{L}-\Rxiln{1}{1}\wk{R}\right).
\end{equation*}

\section{Semi-discrete stability analysis: convective terms}\label{app:stabcondisc}

The element-wise Lyapunov consistency (stability) analysis of the diffusive and reaction terms is 
reported in \cite{aljahdali2024brain}, and it is decoupled from the analysis of the convective term.  
Thus, in this appendix, we need only to show that the addition of the convective term preserves the Lyapunov consistency of the 
fully discrete algorithms for reaction-diffusion PDE systems presented in \cite{aljahdali2024brain}.
That is, for the generic element with index $k$, we need to prove that
\begin{equation}\label{eq:condiscstab}
  \begin{split}
C\equiv&-\left(\wk{k}\right)\Tr\M\sum\limits_{i,l=1}^{dim}\left(\Dxil{i}\matJxilxmk{i}{l}{k}+\matJxilxmk{i}{l}{k}\Dxil{i}\right)
\circ\matUSC{k}{k}\ones\\
      &+\left(\wk{k}\right)\Tr\M\left(\bm{SAT}^{(c)}+\bm{SAT}^{(bc)}\right)\leq 0.
  \end{split}
\end{equation}

As before, we will first drop the boundary condition 
SAT terms, $\bm{SAT}^{(bc)}$, as well as the interface dissipation in the $\bm{SAT}^{(c)}$, 
i.e., the dissipation term $\bm{diss}^{(c)}$, presented in the previous section. 
Thus, after using the 
definition $\Dxil{i}=\M^{-1}\Qxil{i}$, Equation \eqref{eq:condiscstab} reduces to
 \begin{equation}\label{eq:condiscstab2}
  \begin{split}
      C=&-\left(\wk{k}\right)\Tr\sum\limits_{i,l=1}^{dim}\left(\Qxil{i}\matJxilxmk{i}{l}{k}+\matJxilxmk{i}{l}{k}\Qxil{i}\right)
\circ\matUSC{k}{k}\ones\\
      &+\left(\wk{k}\right)\Tr\M\bm{SAT}^{(c)}.
  \end{split}
\end{equation}
Taking the transpose of the second set of terms on the right-hand side of~\eqref{eq:condiscstab2} (first line), 
using the SBP property $\Qxil{i}\Tr=-\Qxil{i}+\Exil{i}$, and the symmetry of $\matUSC{}{}$, we obtain 
 \begin{equation}\label{eq:condiscstab3}
  \begin{split}
C=&-\sum\limits_{i,l=1}^{dim}\left( 
      \left(\wk{k}\right)\Tr\left(\Qxil{i}\matJxilxmk{i}{l}{k}\right)\circ\matUSC{k}{k}\ones-
\ones\Tr\left(\Qxil{i}\matJxilxmk{i}{l}{k}\right)\circ\matUSC{k}{k}\wk{k}
\right)\\
      &+\sum\limits_{i,l=1}^{dim}\left(\wk{k}\right)\Tr\left(\Exil{i}\matJxilxmk{i}{l}{k}\right)\circ\matUSC{k}{k}\ones
      +\left(\wk{k}\right)\Tr\M\bm{SAT}^{(c)}.
  \end{split}
\end{equation}
Using Theorem~\ref{thrm:shuffel} on the first sum on the right-hand side of~\eqref{eq:condiscstab3} gives 
 \begin{equation}\label{eq:condiscstab4}
  \begin{split}
C=&-\sum\limits_{i,l=1}^{dim}\left( 
      \left({\psik{k}}\right)\Tr\Qxilbar{i}\matJxilxmkbar{i}{l}{k}\onesbar-
\onesbar\Tr\Qxilbar{i}\matJxilxmkbar{i}{l}{k}\psik{k}
\right)\\
      &-\sum\limits_{i,l=1}^{dim}\left(\wk{k}\right)\Tr\left(\Exil{i}\matJxilxmk{i}{l}{k}\right)\circ\matUSC{k}{k}\ones
      +\left(\wk{k}\right)\Tr\M \, \bm{SAT}^{(c)}.
  \end{split}
\end{equation}
The first term on the right-hand side of~\eqref{eq:condiscstab4} is zero via the discrete metric identities 
assumption~\eqref{eq:DGCLxm}, \ie, 
\begin{equation*}
  \begin{split}
\sum\limits_{i,l=1}^{dim} 
      \left({\psik{k}}\right)\Tr\Qxilbar{i}\matJxilxmkbar{i}{l}{k}\onesbar&= 
      \sum\limits_{i=1}^{dim}\left({\psik{k}}\right)\Tr\sum\limits_{i=1}^{dim} 
  \Qxilbar{i}\matJxilxmkbar{i}{l}{k}\onesbar\\  
  &= 
      \sum\limits_{i=1}^{dim}\left({\psik{k}}\right)\Tr\M\sum\limits_{i=1}^{dim} 
  \Dxilbar{i}\matJxilxmkbar{i}{l}{k}\onesbar=
  0.
  \end{split}
\end{equation*}
Thus, Equation~\eqref{eq:condiscstab4} reduces to
 \begin{equation}\label{eq:condiscstab5}
  \begin{split}
C=&\sum\limits_{i,l=1}^{dim} 
\onesbar\Tr\Qxilbar{i}\matJxilxmkbar{i}{l}{k}\psik{k}
      -\sum\limits_{i,l=1}^{dim}\left(\wk{k}\right)\Tr\left(\Exil{i}\matJxilxmk{i}{l}{k}\right)\circ\matUSC{k}{k}\ones
      +\left(\wk{k}\right)\Tr\M \, \bm{SAT}^{(c)}.
  \end{split}
\end{equation}
Finally,  using consistency of the SBP operator (i.e., $\Dxil{i}\ones=0\rightarrow\ones\Tr\Qxil{i}=\ones\Tr\Exil{i}$), 
Equation~\eqref{eq:condiscstab5} reduces to  
 \begin{equation}\label{eq:condiscstab6}
  \begin{split}
C=&\sum\limits_{i,l=1}^{dim} 
\onesbar\Tr\Exilbar{i}\matJxilxmkbar{i}{l}{k}\psik{k}
      -\sum\limits_{i,l=1}^{dim}\left(\wk{k}\right)\Tr\left(\Exil{i}\matJxilxmk{i}{l}{k}\right)\circ\matUSC{k}{k}\ones
      +\left(\wk{k}\right)\Tr\M \, \bm{SAT}^{(c)}.
  \end{split}
\end{equation}
All the terms in~\eqref{eq:condiscstab6} are surface terms. We continue the analysis  
by considering the contributions at a particular surface. For the $L$ element, the surface terms are
\begin{equation}\label{eq:conL}
  \begin{split}
      C_{L} &\equiv \sum\limits_{l=1}^{dim}\onesbar\Tr\left(\Rxilnbar{1}{N_{1d}}\right)\Tr\Porthoxilbar{1}\Rxilnbar{1}{N_{1d}}\psik{L}\\
      &-\sum\limits_{l=1}^{dim}\left(\wk{L}\right)\Tr\left(\left(\Rxiln{1}{N_{1d}}\right)\Tr\Porthoxil{1}\Rxiln{1}{N_{1d}}\matJxilxmk{i}{l}{L}\right)\circ\matUSC{L}{L}\ones\\
      &+\sum\limits_{l=1}^{dim}\left(\wk{L}\right)\Tr\left(\left(\Rxiln{1}{N_{1d}}\right)\Tr\Porthoxil{1}\Rxiln{1}{N_{1d}}\matJxilxmk{i}{l}{L}\right)\circ\matUSC{L}{L}\ones\\
      &-\sum\limits_{l=1}^{dim}\left(\wk{L}\right)\Tr\left(\left(\Rxiln{1}{N_{1d}}\right)\Tr\Porthoxil{1}\Rxiln{1}{1}\matJxilxmk{i}{l}{R}\right)\circ\matUSC{L}{R}\ones.
  \end{split}
\end{equation}
Thus, $C_{L}$ reduces to
\begin{equation}\label{eq:conL1}
  \begin{split}
      C_{L}&\equiv\sum\limits_{l=1}^{dim}\onesbar\Tr\left(\Rxilnbar{1}{N_{1d}}\right)\Tr\Porthoxilbar{1}\Rxilnbar{1}{N_{1d}}\psik{L}\\
      &-\sum\limits_{l=1}^{dim}\left(\wk{L}\right)\Tr\left(\left(\Rxiln{1}{N_{1d}}\right)\Tr\Porthoxil{1}\Rxiln{1}{1}\matJxilxmk{i}{l}{R}\right)\circ\matUSC{L}{R}\ones.
  \end{split}
\end{equation}
An identical analysis for the contributions from the $R$ element results in 
\begin{equation}\label{eq:conR1}
  \begin{split}
      C_{R}\equiv&-\sum\limits_{l=1}^{dim}\onesbar\Tr\left(\Rxilnbar{1}{1}\right)\Tr\Porthoxilbar{1}\Rxilnbar{1}{1}\psik{R}\\
      &+\sum\limits_{l=1}^{dim}\left(\wk{R}\right)\Tr\left(\left(\Rxiln{1}{1}\right)\Tr\Porthoxil{1}\Rxiln{1}{N_{1d}}\matJxilxmk{i}{l}{L}\right)\circ\matUSC{R}{L}\ones.
  \end{split}
\end{equation}
Suming~\eqref{eq:conL1} and~\eqref{eq:conR1} yields
\begin{equation}\label{eq:conC1}
  \begin{split}
      C_{L}+C_{R} &=
      \sum\limits_{l=1}^{dim}\onesbar\Tr\left(\Rxilnbar{1}{N_{1d}}\right)\Tr\Porthoxilbar{1}\Rxilnbar{1}{N_{1d}}\psik{L}
-\sum\limits_{l=1}^{dim}\onesbar\Tr\left(\Rxilnbar{1}{1}\right)\Tr\Porthoxilbar{1}\Rxilnbar{1}{1}\psik{R}\\
      &-\sum\limits_{l=1}^{dim}\left(\wk{L}\right)\Tr\left(\left(\Rxiln{1}{N_{1d}}\right)\Tr\Porthoxil{1}\Rxiln{1}{1}\matJxilxmk{i}{l}{R}\right)\circ\matUSC{L}{R}\ones\\
      &+\sum\limits_{l=1}^{dim}\left(\wk{R}\right)\Tr\left(\left(\Rxiln{1}{1}\right)\Tr\Porthoxil{1}\Rxiln{1}{N_{1d}}\matJxilxmk{i}{l}{L}\right)\circ\matUSC{R}{L}\ones.
  \end{split}
\end{equation}
Taking the transpose of the terms in the last sum on the right-hand side of~\eqref{eq:conC1} and 
using the symmetry of $\matUSC{}{}$ (\ie, $\matUSC{L}{R}=\matUSC{R}{L}\Tr$), Equation~\eqref{eq:conC1} reduces 
to
\begin{equation}\label{eq:conC2}
  \begin{split}
      C_{L}+C_{R} &=
      \sum\limits_{l=1}^{dim}\onesbar\Tr\left(\Rxilnbar{1}{N_{1d}}\right)\Tr\Porthoxilbar{1}\Rxilnbar{1}{N_{1d}}\psik{L}
-\sum\limits_{l=1}^{dim}\onesbar\Tr\left(\Rxilnbar{1}{1}\right)\Tr\Porthoxilbar{1}\Rxilnbar{1}{1}\psik{R}\\
      &-\sum\limits_{l=1}^{dim}\left(\wk{L}\right)\Tr\left(\left(\Rxiln{1}{N_{1d}}\right)\Tr\Porthoxil{1}\Rxiln{1}{1}\matJxilxmk{i}{l}{R}\right)\circ\matUSC{L}{R}\ones\\
      &+\sum\limits_{l=1}^{dim}\ones\Tr\left(\matJxilxmk{i}{l}{L}\left(\Rxiln{1}{N_{1d}}\right)\Tr\Porthoxil{1}\Rxiln{1}{N_{1d}}\right)\circ\matUSC{L}{R}\left(\wk{R}\right)\Tr.
  \end{split}
\end{equation}

Now, considering that $\left(\Rxiln{1}{N_{1d}}\right)\Tr\Porthoxil{1}\Rxiln{1}{1}\matJxilxmk{i}{l}{R}$ 
is diagonal and picks off the values of the metrics at the surface nodes and that at the surface nodes 
the metrics are the same, we have
\begin{equation}\label{eq:Rmetrics}
    \matJxilxmk{i}{l}{L}\left(\Rxiln{1}{N_{1d}}\right)\Tr\Porthoxil{1}\Rxiln{1}{N_{1d}} = 
    \left(\Rxiln{1}{N_{1d}}\right)\Tr\Porthoxil{1}\Rxiln{1}{1}\matJxilxmk{i}{l}{R}.
\end{equation}
Using identity \eqref{eq:Rmetrics}, Equation~\eqref{eq:conC2} becomes
\begin{equation}\label{eq:conC3}
  \begin{split}
      C_{L}+C_{R} &=
      \sum\limits_{l=1}^{dim}\onesbar\Tr\left(\Rxilnbar{1}{N_{1d}}\right)\Tr\Porthoxilbar{1}\Rxilnbar{1}{N_{1d}}\psik{L}
-\sum\limits_{l=1}^{dim}\onesbar\Tr\left(\Rxilnbar{1}{1}\right)\Tr\Porthoxilbar{1}\Rxilnbar{1}{1}\psik{R}\\
      &-\sum\limits_{l=1}^{dim}\left(\wk{L}\right)\Tr\left(\left(\Rxiln{1}{N_{1d}}\right)\Tr\Porthoxil{1}\Rxiln{1}{1}\matJxilxmk{i}{l}{R}\right)\circ\matUSC{L}{R}\ones\\
      &+\sum\limits_{l=1}^{dim}\ones\Tr\left(\left(\Rxiln{1}{N_{1d}}\right)\Tr\Porthoxil{1}\Rxiln{1}{N_{1d}}\matJxilxmk{i}{l}{R}\right)\circ\matUSC{L}{R}\left(\wk{R}\right)\Tr.
  \end{split}
\end{equation}
Now, we apply Theorem~\ref{thrm:shuffel} on the last two sums on the right-hand side of~\eqref{eq:conC3} 
which gives 
\begin{equation}\label{eq:conC4}
  \begin{split}
      C_{L}+C_{R} &=
      \sum\limits_{l=1}^{dim}\onesbar\Tr\left(\Rxilnbar{1}{N_{1d}}\right)\Tr\Porthoxilbar{1}\Rxilnbar{1}{N_{1d}}\psik{L}
-\sum\limits_{l=1}^{dim}\onesbar\Tr\left(\Rxilnbar{1}{1}\right)\Tr\Porthoxilbar{1}\Rxilnbar{1}{1}\psik{R}\\
      &-\sum\limits_{l=1}^{d}\left({\psik{L}}\right)\Tr\left(\left(\Rxilnbar{1}{N_{1d}}\right)\Tr\Porthoxilbar{1}\Rxilnbar{1}{1}\matJxilxmkbar{i}{l}{R}\right)\onesbar\\
      &+\sum\limits_{l=1}^{dim}\onesbar\Tr\left(\left(\Rxilnbar{1}{N_{1d}}\right)\Tr\Porthoxilbar{1}\Rxiln{1}{N_{1d}}\matJxilxmkbar{i}{l}{R}\right)\psik{R}=0,
  \end{split}
\end{equation}
where the final equality follows from using identity~\eqref{eq:Rmetrics}. 
Thus, the interface 
terms do not destroy the result generated from the discrete diffusion operator analysis reported in \cite{aljahdali2024brain}. 

Next we analyze the effect of the additional dissipation. The contribution from element $L$ is 
\begin{equation}\label{eq:dissconL}
    diss^{(c),L} \equiv -\wk{L}\left(\Rxiln{1}{N_{1d}}\right)\Tr\Porthoxil{1}\matdutildedwjk{1}{}
    \left(\Rxiln{1}{N_{1d}}\wk{L}-\Rxiln{1}{1}\wk{R}\right),
\end{equation}
while for the $R$ element the contribution is
\begin{equation}\label{eq:dissconR}
    diss^{(c),R} \equiv -\wk{R}\left(\Rxiln{1}{1}\right)\Tr\Porthoxil{1}\matdutildedwjk{1}{}
    \left(\Rxiln{1}{1}\wk{R}-\Rxiln{1}{N_{1d}}\wk{L}\right).
\end{equation}
Adding~\eqref{eq:dissconL} to~\eqref{eq:dissconR}  gives
\begin{equation*}
    diss^{(c),L}+ diss^{(c),R} = -\left(\Rxiln{1}{N_{1d}}\wk{L}-\Rxiln{1}{1}\wk{R}\right)\Tr\Porthoxil{1}\matdutildedwjk{1}{}
    \left(\Rxiln{1}{N_{1d}}\wk{L}-\Rxiln{1}{1}\wk{R}\right),
\end{equation*}
which is a negative semi-definite term since $\matdutildedwjk{1}{}$ is symmetric 
positive-semi definite. Therefore, it does not impact the Lyapunov consistency statement.

\section{Mesh Convergence Study}\label{sec:mesh-Convergence}

\begin{figure}[H]
    \centering 
             \begin{minipage}{.52\textwidth}
             \begin{subfigure}{\textwidth}
            \centering
            \includegraphics[width=\textwidth]{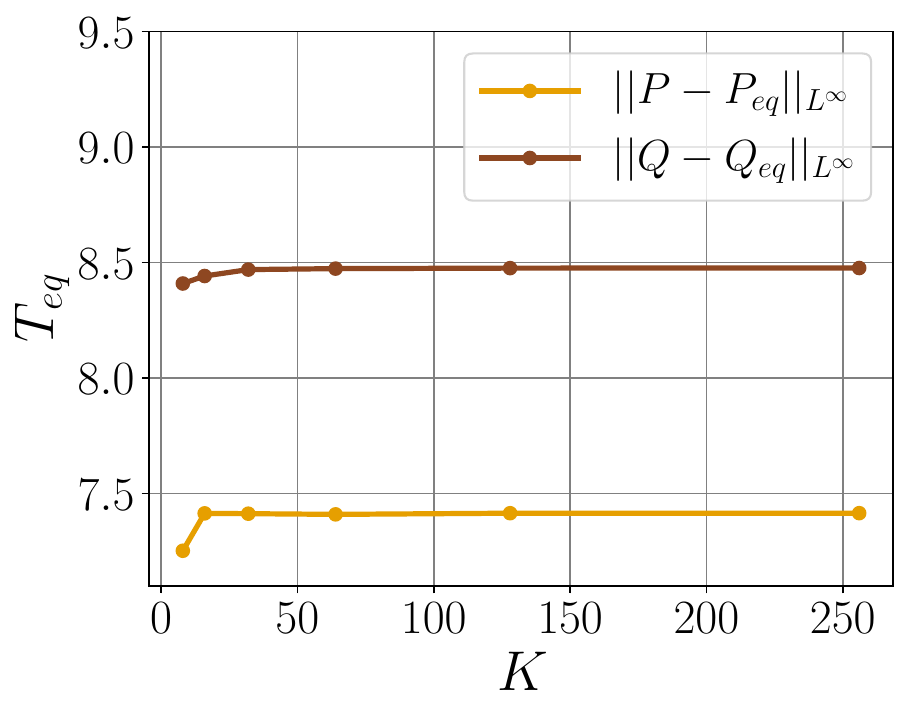}
            \end{subfigure}
              \end{minipage}
           \caption{
           Effect of the mesh refinement on the convergence to the equilibrium for the dimerization model \eqref{eq:pde-system-monomers-dimers}. \label{fig:mesh_convergence_study}}
    \end{figure}

To evaluate the impact of mesh refinement on the convergence to the equilibrium 
point of the proposed algorithms, we conduct numerical simulations of the 
dimerization model \eqref{eq:pde-system-monomers-dimers} using identical parameter 
settings, boundary conditions, and initial conditions. The simulations are performed
within a cubic domain with a side length of $1$ unit,  discretized with 
$8$, $16$, $32$, $64$, $128$, and $256$ quadrilateral cells in each coordinate direction.
We determine the time $T_{eq}$ at which the maximum norm of the difference between 
the numerical solution and the equilibrium value becomes less than or equal 
to $10^{-8}$. The $T_{eq}$ results are presented in Figure \ref{fig:mesh_convergence_study}, 
plotted against the number of cells, $K$. As the cell density increases, the solution 
progressively approaches the equilibrium point. Notably, the reported times for the grids 
of $128$ and $256$ cells demonstrate negligible differences and are practically indistinguishable.

\section{Discrete Error Norms}\label{sec:discrete-error-norms}
In this work, we used the following formula of $L^{1}$, $L^{2}$, $L^{\infty} $ norms to compute the error:
\begin{equation}
\begin{aligned}
& \|\bmU\|_{L^{1}}=\Omega_{c}^{-1}\sum
    \limits_{j=1}^{N_{el}}\onek\Tr \, \Mk \, \matJkmod \,
    \textrm{abs}\left(\bmU_{j}\right),\\
&\|\bmU\|_{L^{2}}^2=\Omega_{c}^{-1}\sum\limits_{j=1}^{N_{el}}\bmU_{j}\Tr \,
    \Mk \, \matJkmod \, \bmU_{j},\\
& \|\bmU\|_{L^{\infty}}=\max\limits_{j=1\dots
    N_{el}}\textrm{abs}\left(\bmU_{j}\right),
\end{aligned}
\end{equation}
where  $\matJkmod$  represents the metric Jacobian of the curvilinear transformation mapping the physical space to the computational space of the $j^{th}$  hexahedral element. The mesh's total number of hexahedral elements is denoted by $N_{el}$. In addition, $\Omega_{c}$ corresponds to the volume of the domain of interest $\Omega$ associated with the considered PDE or system of PDEs. The volume $\Omega_{c}$ is determined as 
\begin{equation}
 \Omega_{c}\equiv\sum\limits_{\kappa=1}^{K}\onek_{\kappa}\Tr\M^{\kappa}\matJkmod\onek_{\kappa}, 
 \end{equation}
where $\onek_{\kappa}$ denotes a vector of ones with a length equivalent to the number of nodes on the $\kappa^{th}$ element.

\end{appendices}

\section*{Acknowledgments}
The work described in this paper was supported by King Abdullah University of Science 
and Technology through the grant number BAS/1/1663-01-01.
The authors are also thankful for the computing resources of the Supercomputing Laboratory 
at King Abdullah University of Science and Technology.

\section*{Declarations}

\bmhead{Conflict of interest}
The authors declare that they have no conflict of interest.

\bmhead{Availability of code, data, and materials}

\end{document}